\documentclass[11pt,centertags,twoside,reqno]{amsart}

\usepackage{amsmath,amstext,amsthm,amscd,typearea,hyperref}
\usepackage{amssymb}
\usepackage{algorithm}
\usepackage{algpseudocode}
\usepackage{a4wide}
\usepackage[mathscr]{eucal}
\usepackage{mathrsfs}
\usepackage{typearea}
\usepackage{CJKutf8}
\usepackage{dsfont}
\usepackage{hyperref}
\usepackage{charter}
\usepackage{pdfsync}
\usepackage[a4paper,width=16.2cm,top=2.3cm,bottom=2.5cm]{geometry}

\numberwithin{equation}{section}

\usepackage{cleveref}

\crefformat{section}{\S#2#1#3} 
\crefformat{subsection}{\S#2#1#3}
\crefformat{subsubsection}{\S#2#1#3}

\usepackage[
sortcites=true,
backend=biber,
style=alphabetic,
sorting=nyt
]{biblatex}
\renewbibmacro{in:}{\ifentrytype{article}{}{\printtext{\bibstring{in}\intitlepunct}}}

\usepackage[dvipsnames]{xcolor}

\hypersetup{
    colorlinks=true,
    linkcolor=black,
    filecolor=magenta,
    urlcolor=NavyBlue,
    citecolor = NavyBlue,
}

\newenvironment{conjecturedequations}
  {%
    \let\originaltheequation\theequation
    \renewcommand{\theequation}{\originaltheequation*}%
    \ignorespaces
  }
  {%
    \ignorespacesafterend
  }

\usepackage{yhmath}

\usepackage{datetime2}
\usepackage{graphicx}
\usepackage{verbatim}
\usepackage{subcaption}
\usepackage{blindtext}
\usepackage{ragged2e}
\usepackage{xcolor}
\usepackage{array}
\usepackage{xypic}
\usepackage{tikz}
\usepackage{nameref}
\usepackage{extarrows}
\usepackage{enumitem}
\usepackage{tabularx}
\usepackage{booktabs}
\usepackage[labelfont=bf,format=plain,justification=raggedright,singlelinecheck=false]{caption}

\usetikzlibrary{backgrounds}
\usetikzlibrary{arrows, decorations.pathmorphing}
\usetikzlibrary{quotes,angles}

\usepackage{extarrows}

\allowdisplaybreaks

\newtheorem{theorem}{Theorem}[section]

\newtheorem{definition}[theorem]{Definition}
\newtheorem{proposition}[theorem]{Proposition}

\newtheorem{lemma}[theorem]{Lemma}

\newtheorem{example}[theorem]{Example}

\newcommand{\bkt}[1]{\biggl( #1 \biggr) }
\newcommand{\Bkt}[1]{\Biggl( #1 \Biggr) }
\newcommand{\sbkt}[1]{\bigl( #1 \bigr) }

\newcommand{\bigbkt}[1]{\biggl[ #1 \biggr] }

\newcommand{\norm}[1]{\bigl|\bigl| #1 \bigr| \bigr|}

\newcommand{\Norm}[1]{\biggl|\biggl| #1 \biggr| \biggr|}

\newcommand{\abs}[1]{\biggl| #1 \biggr| }

\newcommand{\set}[1]{\biggl \{ #1 \biggr\} }

\emergencystretch=\maxdimen
\AtBeginDocument{%
  \setlength{\evensidemargin}{\oddsidemargin}%
}

\title[Near-optimal bounds for the SHF and the delta-Bose gas]{%
  Near-optimal bounds for the critical 2d SHF and the delta-Bose gas}

\author{Xue-Mei Li}
\address{Institute of Mathematics, École polytechnique fédérale de Lausanne (EPFL),
1015 Lausanne, Switzerland}
\email{xue-mei.li@epfl.ch}

\author{Te-Chun Wang}
\address{Institute of Mathematics, École polytechnique fédérale de Lausanne (EPFL),
1015 Lausanne, Switzerland}
\email{te-chun.wang@epfl.ch}

\date{\today}

\begin{document}

\begin{abstract} \centering \justifying 
We establish a near-optimal upper bound $\exp(\exp(CN^{1+o(1)}))$ for the $N$-th moment of the critical two-dimensional stochastic heat flow, improving on the bound $\exp(\exp(CN^{2}))$ of Gu, Quastel, and Tsai \cite{Gu_2d}. Near-optimality follows from the lower bound $\exp(\exp(cN))$ of Ganguly and Nam \cite{ganguly2025sharp}. These estimates establish the exponential scale of the ground-state energy of the two-dimensional $N$-body delta-Bose gas predicted by Rajeev’s conjecture \cite{rajeev1999condensation}, up to a vanishing error in the exponent.  As a further step toward this conjecture, we derive, for each fixed $N$, a formula for the ground-state energy in terms of the bottom of the spectrum of a principal operator. The derivation relies on an analytic continuation of the resolvent formula in \cite{Gu_2d}.
\end{abstract}

\pagestyle{headings}
\markboth
  {X.-M. Li, T.-C, Wang}
  {Near-optimal bounds for the SHF and the delta-Bose gas}

\maketitle

\thispagestyle{empty}

\medskip

\begin{minipage}{0.9\textwidth}
\noindent {\centering \justifying   Classification AMS 2020}: 60H15, 81Q10, 81V73. 
\end{minipage}

\medskip

\begin{minipage}{0.9\textwidth}
\noindent {\centering \justifying  Keywords: Kardar-Parisi-Zhang equation; stochastic heat equation; continuum directed polymer; critical 2d stochastic heat flow; delta-Bose gas.} 
\end{minipage}

\medskip


\tableofcontents

\section{Introduction}

The critical two-dimensional stochastic heat flow (SHF) is a universal measure-valued process $\mathscr{Z}^{\theta} = (\mathscr{Z}_{s,t}^{\theta}(dx,dy))_{0\leq s<t<\infty}$ that arises as the scaling limit of directed polymers in a random environment when the disorder strength is tuned within a critical window parametrized by $\theta$ \cite{Heat_flow}. The directed polymer model was introduced in the physics literature by Huse and Henley \cite{DDP2} as a simplified description of interfaces in Ising systems with random couplings. A similar discrete model was later studied rigorously by Imbrie and Spencer \cite{DDP}. Directed polymer models have since attracted considerable mathematical interest; see the review \cite{comets2004probabilistic} and the references therein.

The SHF also arises as a renormalized limit of the solution to the two-dimensional stochastic heat equation (SHE) \cite{Heat_flow,tsai2024stochastic}, formally written as
\begin{align}\label{equation:SHE}
\partial_t \mathbf{X}(x,t)
=
\frac{1}{2}\Delta_x \mathbf{X}(x,t)
+
\beta \mathbf{X}(x,t)\xi(x,t),
\end{align}
where $\xi$ is space-time white noise and the coupling constant $\beta$ indicates the noise strength. The SHF has also been obtained as a scaling limit of other stochastic PDEs and directed polymer models; see \cite{surendranath2026continuoustimedirectedpolymerscritical,long_SHE,drillick2026universalityresultcritical2d}.

The SHE is a central example of a stochastic PDE, with rich physical and probabilistic interpretations. In particular, under the Cole–Hopf transformation, the logarithm of its solution formally satisfies the Kardar–Parisi–Zhang (KPZ) equation \cite{KPZ}, which has been studied extensively in one spatial dimension; see \cite{solving,Burgers,ferrari2010random} and see also \cite{quastel} for a review. The higher moments of the SHE are related to the delta-Bose gas, which has a long history in the physics and functional analysis literature; see \cite{GalloneMichelangeli2023,albeverio2012solvable} for reviews. The SHE is also closely related to the continuum directed random polymer, with its solution serving as the polymer’s partition function. In one spatial dimension, this polymer was rigorously constructed in \cite{Alberts_2013}.

Due to the roughness of the noise, the SHE is well-defined within classical Itô theory \cite{walsh1986introduction} only in spatial dimension one. The two-dimensional case is critical: the equation is formally invariant under parabolic scaling. However, the general theories of singular stochastic PDEs in \cite{RS,Para} only apply to  subcritical equations. A standard approach in the critical dimension is to consider the regularized SHE and study the limit as the regularization is removed:
\begin{align}\label{equation:mollified_SHE}
    \partial_t \mathbf{X}_{\varepsilon}^{\theta}(x,t)
    =
    \frac{1}{2}\Delta_x \mathbf{X}_{\varepsilon}^{\theta}(x,t)
    +
    \beta_{\varepsilon}(\theta)\mathbf{X}_{\varepsilon}^{\theta}(x,t)\xi_{\varepsilon}(x,t),
\end{align}
where the regularized noise is
\begin{align*}
    \xi_{\varepsilon}(x,t)
    :=
    \int_{\mathbb{R}^{2}} dy \phi_{\varepsilon}(x-y)\xi(y,t).
\end{align*}
Here $\phi_{\varepsilon}(x):=\varepsilon^{-2}\phi(x/\varepsilon)$, with $\phi\in C_c^\infty(\mathbb{R}^{2};\mathbb{R}_{\geq0})$ and $\int_{\mathbb{R}^{2}}\phi(x)\,dx=1$. Its spatial covariance~is
\begin{equation}\label{eq:R}
 R_{\varepsilon}(x)
    :=
    \int_{\mathbb{R}^{2}}dy \phi_{\varepsilon}(y+x)\phi_{\varepsilon}(y).
\end{equation}
A nontrivial limit is obtained by tuning the noise strength within the critical window parametrized by $\theta\in\mathbb{R}$:
\begin{align}\label{equation:critical_scaling}
    \beta_{\varepsilon}(\theta)^2
    &:=
    \frac{2\pi}{\log\varepsilon^{-1}}
    +
    \frac{\pi}{(\log\varepsilon^{-1})^2}
    \bkt{
        \theta-2\log2+2\gamma_{\textup{EM}}
        +2\int_{\mathbb{R}^{2}\times\mathbb{R}^{2}} dxdx'
        R_1(x)\log|x-x'|R_1(x')
    },
\end{align}
where $\gamma_{\textup{EM}}$ is the Euler–Mascheroni constant. Under this scaling, for every fixed $s<t$, the fundamental solution satisfies
\begin{align}\label{equation:SHF}
    \mathbf{X}_{\varepsilon}^{\theta}((x',t);(x,s))\,dx\,dx'
    \overset{d}{\longrightarrow}
    \mathscr{Z}^{\theta}_{s,t}(dx,dx')
    \quad\text{in }\mathcal{M}_{+}(\mathbb{R}^{4})
\end{align}
as $\varepsilon\downarrow0$. Here $\mathbf{X}_{\varepsilon}^{\theta}(\cdot;(x,s))$ solves \eqref{equation:mollified_SHE} with initial condition $\delta_x$ at time $s$, and $\mathcal{M}_{+}(\mathbb{R}^{4})$ is the space of locally finite positive Borel measures on $\mathbb{R}^{4}$, equipped with the vague topology: the coarsest topology making $\mu\mapsto\mu(f)$ continuous for every $f\in C_c(\mathbb{R}^{4})$. See \cite{Heat_flow} and \cite[Proposition 1.5]{tsai2024stochastic} for details.

We establish upper bounds for the higher moments of the SHF, resolving the open problem on moment growth posed in \cite[Section 11.2, item 5]{caravenna2026critical2dstochasticheat}, up to an error quantified in \eqref{equation:rewrite_exponent}. Moment estimates played an important role in the construction of the SHF. Under the critical scaling \eqref{equation:critical_scaling}, second-moment convergence for the regularized SHE \eqref{equation:mollified_SHE} was established in \cite{Bertini_1998}, using the theory of two-dimensional point interactions developed in \cite[Chapter I.5]{albeverio2012solvable}. Third-moment convergence was subsequently proved in \cite{caravenna2019moments}, while arbitrary higher moments were treated in various settings in \cite{Gu_2d,chen}. These results yield tightness and hence subsequential limits. However, owing to the rapid growth of the moments \cite{caravenna2023critical,ganguly2025sharp}, moment convergence alone does not identify the limiting law. The fundamental problem of uniqueness of the scaling limit was resolved in \cite{Heat_flow}. The properties of the SHF have since been studied extensively; see \cite{berger2026critical2dkpzprocess} and see also
\cite{caravenna2026disordered} for a review.

\section{Main results}

We first present moment and tail probability estimates for the SHF in Section \ref{section:main_result_moment_bound}. In Section \ref{section:main_result_ground_state_energy}, we give bounds on the ground-state energy of the delta-Bose gas and a formula expressing this energy in terms of a principal operator.

\subsection{Part I: Stochastic heat flow} \label{section:main_result_moment_bound}

We denote by
\begin{align} \label{equation:density_Brownian_motion}
    G_{t}(x) := \frac{1}{2\pi t} \exp(-|x|^{2}/(2t))
\end{align}
the transition density of two-dimensional Brownian motion. We write
\begin{align*}
    \mathscr{Z}^{\theta}_{t}(dx,dy) = \mathscr{Z}^{\theta}_{0,t}(dx,dy),
\end{align*}
and let $\mathbb{E}$ denote expectation with respect to the underlying probability measure.

\begin{definition}
Let $f,h$ be non-negative Borel-measurable functions on $\mathbb{R}^{2}$. We say that a pair $(h,f)$ satisfies \textit{the upper-bound hypotheses} if for some $\nu_{f},C_{f},D_{h}>0$
\begin{equation} \label{eq:UBH} \tag{UBH}
    h \leq D_{h}, \quad f \leq C_{f} G_{\nu_{f}}.
\end{equation}
This class of pairs is denoted by \eqref{eq:UBH}.
\end{definition}

The upper-bound hypotheses hold for every pair of non-negative continuous functions with compact support. Our first main result gives a near-optimal upper bound on the higher moments of the SHF. Its proof is given in Section \ref{section:proof_idea_theorem:main_result_higher_moment}.
\begin{theorem} \label{theorem:main_result_higher_moment}
Suppose that $(h,f) \in \eqref{eq:UBH}$ and $\theta\in \mathbb{R}$. There exists a positive constant
$C_{\eqref{equation:main_result_higher_moment}}$, depending only on $\theta,h,f$, such that, for all $t\geq 1$
and integers $N\geq 3$,
\begin{align} \label{equation:main_result_higher_moment}
    \mathbb{E}[\mathscr{Z}^{\theta}_{t}(h,f)^{N}] \leq \exp(t\exp(C_{\eqref{equation:main_result_higher_moment}}N \log N)).
\end{align}
\end{theorem}

Together with the identity \eqref{equation:connection_SHF_DBG} and Theorem \ref{theorem:2d_DBG} below, this yields the following bound:
\begin{align} \label{equation:main_result_higher_moment_L2} 
    \sup_{||f||_{L^{2}(\mathbb{R}^{2})} = ||h||_{L^{2}(\mathbb{R}^{2})} = 1}\mathbb{E}[|\mathscr{Z}^{\theta}_{t}(h,f)|^{N}] \leq \exp(t\exp(C_{\eqref{equation:main_result_higher_moment_L2}}N \log N)) \quad \forall N \geq 3,
\end{align}
where the constant $C_{\eqref{equation:main_result_higher_moment_L2}}$ depends only on $\theta$.

For comparison with earlier results, we set $t=1$ in the following
discussion. The previously known upper bound of the form
$\exp(\exp(CN^{2}))$ follows from the results of Gu, Quastel, and Tsai
\cite{Gu_2d}. Their proof adapts the resolvent method of Dimock and
Rajeev \cite{JDimock_2004}, which builds on earlier ideas of Rajeev
\cite{rajeev1999condensation}.
Dimock and Rajeev constructed the renormalized Hamiltonian of the
two-dimensional delta-Bose gas, on the symmetric subspace of $L^{2}(\mathbb{R}^{2N})$, as a resolvent limit of momentum-cutoff
Hamiltonians and obtained a lower bound on its spectrum. Gu, Quastel, and Tsai instead considered the spatial regularization arising from the critically scaled mollified SHE \eqref{equation:mollified_SHE}. They established strong resolvent
convergence of the regularized Hamiltonians, on the full $L^{2}(\mathbb{R}^{2N})$, to the limiting
Hamiltonian $\mathscr{H}_{N}^{\theta}$  and derived an explicit diagrammatic representation of the semigroup $\exp(-t \mathscr{H}_{N}^{\theta})$; see Section \ref{section:main_result_ground_state_energy} for the definition of the Hamiltonian. 
This spatial regularization yields the connection between the limiting moments of the
mollified SHE \eqref{equation:mollified_SHE} and the semigroup $\exp(-t \mathscr{H}_{N}^{\theta})$. In terms of the SHF, this connection is expressed by
the following identity:
\begin{align} \label{equation:connection_SHF_DBG}
    \mathbb{E}[\mathscr{Z}^{\theta}_{t}(h,f)^{N}] = \langle h^{\otimes N}, \exp(-t \mathscr{H}_{N}^{\theta}) f^{\otimes N} \rangle_{L^{2}(\mathbb{R}^{2N})}, 
\end{align}
which holds for $t>0$, $N\geq 2$, and real-valued functions
$h,f\in L^{2}(\mathbb{R}^{2})$. Combining this identity with the spectral theorem and the resolvent
estimates in \cite{Gu_2d} yields the upper moment bound
$\exp(\exp(CN^{2}))$ at $t=1$.

From \eqref{equation:main_result_higher_moment} and
\begin{align} \label{equation:rewrite_exponent}
    N\log N = N^{1+o(1)}, \quad o(1) = \frac{\log \log N}{\log N},
\end{align}
we have the following bound 
\begin{align} \label{equation:main_result_higher_moment_1+o(1)}
    \mathbb{E}[\mathscr{Z}^{\theta}_{t}(h,f)^{N}] \leq \exp(t\exp(C_{\eqref{equation:main_result_higher_moment}}N^{1+o(1)})),
\end{align}
which can be compared directly with the expected optimal upper moment bound
\begin{align*}
    \exp(\exp(CN))
\end{align*}
and the lower moment bound $\exp(\exp(cN))$  established by Ganguly and Nam \cite{ganguly2025sharp}.
The  upper bound $\exp(\exp(CN))$ is supported both by Rajeev's conjecture \eqref{conjecture:Rajeev1999} for the ground-state energy of the delta-Bose gas and by the heuristic argument underlying our proof which is presented in Section \ref{section:informal_explanation}.  The lower bound in \cite{ganguly2025sharp}
strengthens the earlier lower moment bound $\exp(cN^{2})$ obtained by Caravenna, Sun, and Zygouras \cite{caravenna2023critical}. We refer to Section \ref{section:discussion_Rajeev_ conjecture} for a discussion of the optimal constant $C$ appearing in the conjectured upper bound $\exp(\exp(CN))$.

Our improvement from $\exp(\exp (CN^2))$ to $\exp(t\exp(CN \log N)) $ comes from extracting decay from the normalization factors of the Brownian transition densities appearing in the collision representation \eqref{definition:semigroup} for the higher moments.  Heuristically, in these formulas, the denominators contain sums of order $N$ independent, uniformly distributed random variables, suggesting a decay factor of $1/N$.

The key idea  is to use a collision-time Markov chain \eqref{equation:coulping_process} in the analysis of the higher moments. The contribution of diagrams of length $m$ is bounded in terms of a path functional involving the $(m-1)$-fold iterate of a weighted transition operator $\mathscr{Q}_N$; see Proposition \ref{proposition:FK_representation}. A key step is then to establish a weighted-Lyapunov inequality; see Proposition \ref{proposition:main_contractiv_inequality}. In their proof of the lower moment bound, Ganguly and Nam used the arithmetic mean-geometric mean (AM-GM) inequality to extract the required decay; see \cite[Lemma 3.5]{ganguly2025sharp}. In Section \ref{section:trade_off}, we explain the difficulties that arise when applying this inequality to obtain an upper moment bound.

Understanding the statistical properties of the SHF is a fundamental problem. As a measure-valued stochastic process, its probability law was characterized by Tsai \cite{tsai2024stochastic}, with the associated martingale formulations subsequently established by Chen  \cite{chen2025martingale} and Nakashima \cite{nakashima2025martingalemeasureassociatedcritical}. Alongside these process-level descriptions, estimates for the tail probability provide quantitative information about the SHF at fixed times. In this direction, Ganguly and Nam \cite[Theorem 1.4]{ganguly2025sharp} established two-sided bounds describing the SHF at large scale.  We improve their bounds in the following theorem.

\begin{theorem}
\label{corollary:two_side_bound_tail_probability}
The following holds.
\begin{enumerate} [label = (\roman*).]
    \item Let $(h,f)\in \eqref{eq:UBH}$, $t\geq 1$, and $\theta\in \mathbb{R}$. Then there exists  $z_{0}>0$ such that for all $z>z_{0}$,
    \begin{align} \label{equation:upper_bound_tail_probability}
    \mathbb{P}(\mathscr{Z}^{\theta}_{t}(h,f) \geq z) \leq \exp\bkt{-\frac{C_{\eqref{equation:upper_bound_tail_probability}}t}{\log \log ((\log z)/t)} \bkt{\log \frac{\log z}{t} } \bkt{\frac{\log z}{t} } }.
    \end{align}
    \item Let $t\geq 1$ and $\theta\in \mathbb{R}$, and let $ h,f \in C_{c}(\mathbb{R}^{2};\mathbb{R}_{\geq 0})$, not identically zero. Then there exists a constant $z_{0}$ such that for each $z>z_{0}$,
    \begin{align} \label{equation:lower_bound_tail_probability}
    \exp\bkt{-t\bkt{\frac{\log z}{t}}^{c_{\eqref{equation:lower_bound_tail_probability}} \log\log ((\log z)/t)} } \leq \mathbb{P}(\mathscr{Z}^{\theta}_{t}(h,f) \geq z).
    \end{align}
\end{enumerate}
Here $C_{\eqref{equation:upper_bound_tail_probability}}$ and $c_{\eqref{equation:lower_bound_tail_probability}}$ are positive constants depending only on $f,h,\theta$, independent of $t$, while $z_{0}$ may depend on $t$.
\end{theorem}

At $t=1$, the improvements over \cite[Theorem 1.4]{ganguly2025sharp}
are illustrated by the following comparisons of the bounds for large $z$:
\begin{align*}
    (\log z)^{(\log\log z)^{1+o(1)}}
    \gtrsim (\log z)^{c\log\log\log z},\quad \frac{(\log z)\log\log z}{\log\log\log z}
    \gtrsim (\log z)\sqrt{\log\log z}.
\end{align*}
Our estimates also allow both endpoint test functions to be in $C_{c}(\mathbb{R}^{2};\mathbb{R}_{\geq 0})$, whereas \cite[Theorem 1.4]{ganguly2025sharp} concerns the case in which one endpoint test function is constant.
To obtain the lower bound in this localized setting, we use
Ganguly and Nam's lower moment bound \cite[Theorem 1.2]{ganguly2025sharp} as the starting point
for a localization argument. The proof of these tail probability estimates is presented in
Section \ref{section:tail_probability_SHF}.

\subsection{Part II: Delta-Bose gas} \label{section:main_result_ground_state_energy}

Throughout this paper, we write
\begin{align*}
    x = (x_{1},\ldots,x_{N}) \in \mathbb{R}^{2N}, \quad N\geq 3,
\end{align*}
where $x_{i} \in \mathbb{R}^{2}$ denotes the position of the $i$-th particle. We now consider a nonrelativistic quantum system of $N$ particles in two dimensions, formally described by the Hamiltonian
\begin{align} \label{definition:formal_Hamiltonian}
    \mathscr{H}^{\theta}_{N} := -\frac{1}{2} \sum_{i = 1}^{N}\Delta_{x_{i}} -\beta(\theta)^{2} \sum_{(\mathfrak{i},\mathfrak{i}')\in \mathcal{E}_{N}} \delta_{0}(x_{\mathfrak{i}'} - x_{\mathfrak{i}}).
\end{align}
Here $\beta(\theta)$ formally represents the interaction strength and 
\begin{align}  \label{equation:index_set}
    \mathcal{E}_{N} := \set{(\mathfrak{i},\mathfrak{i}'): 1\leq \mathfrak{i}<\mathfrak{i}'\leq N}
\end{align}
is the collision-pair index set. In \eqref{definition:formal_Hamiltonian}, the first term represents the free kinetic energy, while the second describes an attractive pairwise contact interaction.

The formal Hamiltonian \eqref{definition:formal_Hamiltonian} does not directly define a self-adjoint operator on $L^2(\mathbb{R}^{2N})$. A rigorous realization was obtained in \cite{Teta} through a quadratic-form construction based on a momentum-cutoff approximation. Resolvent constructions were subsequently developed in \cite{JDimock_2004,Gu_2d}, yielding the corresponding self-adjoint Hamiltonian.

We use the realization obtained by regularizing the two-body interaction and tuning its strength within a critical window as the regularization is removed. For sufficiently small $\varepsilon>0$, consider the regularized Hamiltonian
\begin{align}  \label{H_varepsilon}
    \mathscr{H}^{\theta}_{N;\varepsilon} := -\frac{1}{2}\sum_{i=1}^{N}\Delta_{x_{i}}  -\beta_{\varepsilon}(\theta)^{2} \sum_{(\mathfrak{i},\mathfrak{i}')\in \mathcal{E}_{N}} R_{\varepsilon}(x_{\mathfrak{i}'} - x_{\mathfrak{i}}) \quad \text{on} \quad L^{2}(\mathbb{R}^{2N}),
\end{align}
where $\beta_{\varepsilon}(\theta)$ is given in \eqref{equation:critical_scaling} and $R_{\varepsilon}$ is defined in \eqref{eq:R}. Our parametrization in \eqref{equation:critical_scaling} differs from that in \cite[(1-7) and (1-8)]{Gu_2d}. Under our convention, their parameter $\beta_{\star}$ is identified with
\begin{align} \label{definition:beta_star}
    \beta_{\star} = \beta_{\star}(\theta) := \theta.
\end{align}
With this identification, Gu, Quastel, and Tsai proved that, for each fixed $N$ and $\theta$, $\mathscr{H}^{\theta}_{N;\varepsilon}$ converges in the strong resolvent sense as $\varepsilon\to0$ to a self-adjoint operator on $L^{2}(\mathbb{R}^{2N})$ \cite[Theorem 1.6 and Remark 1.8]{Gu_2d}. We denote this limit by $\mathscr{H}^{\theta}_{N}$, using the same notation as for the formal Hamiltonian \eqref{definition:formal_Hamiltonian}.

They also established that, for each $N\geq 2$ and $t>0$,
\begin{align} \label{equation:L_2_convergence}
    \exp(-t \mathscr{H}^{\theta}_{N;\varepsilon}) \longrightarrow\exp(-t \mathscr{H}^{\theta}_{N}) \quad \text{strongly on} \quad L^{2}(\mathbb{R}^{2N}) \quad \text{as } \varepsilon \to 0.
\end{align}
See \cite[(7-5)]{Gu_2d} for details. A diagrammatic path-integral representation of the limiting semigroup is given in Section \ref{section:path_integral_representation}. We refer to \cite{chen} for related convergence results in a different setting and to \cite{10.1214/25-AAP2179,10.1214/26-ECP771,chen2025stochasticmotionstwodimensionalmanybody_1,chen2025stochasticmotionstwodimensionalmanybody_2,chen2025stochasticmotionstwodimensionalmanybody_3} for studies of Feynman-Kac representations of this semigroup.

We now turn to the bosonic case. Let $L^{2}_{\textup{sym}}(\mathbb{R}^{2N})$ denote the subspace of $L^{2}(\mathbb{R}^{2N})$ consisting of functions invariant under permutations of the $N$ particle coordinates. Denote the restriction of $\mathscr{H}_{N;\varepsilon}^{\theta}$ to this subspace by $\mathscr{H}_{N;\varepsilon}^{\theta;\textup{sym}}$. Its strong-resolvent limit, denoted by $\mathscr{H}_{N}^{\theta;\textup{sym}}$, is self-adjoint on $L^{2}_{\textup{sym}}(\mathbb{R}^{2N})$ and agrees with the restriction of $\mathscr{H}_{N}^{\theta}$ to this subspace. Throughout this paper, we call $\mathscr{H}_{N}^{\theta;\textup{sym}}$ the Hamiltonian of the \textit{$N$-body delta-Bose gas}.

\subsubsection{\textbf{Ground-state energy estimates}} \label{section:main_result_ground_state_estimate}
The large-$N$ behavior of the ground-state energy is  fundamental in quantum many-body theory. For the two-dimensional delta-Bose gas, the best previously known lower bound was of the form $-\exp(CN^{2})$, established in \cite{Teta,JDimock_2004,Gu_2d}. We obtain near-sharp two-sided bounds for the ground-state energies of both the unrestricted and bosonic Hamiltonians. For a self-adjoint operator $\mathscr{H}$, write
\begin{align} \label{bottom_spectrum}
    \inf\mathscr{H}:=\inf\sigma(\mathscr{H}).
\end{align}
The following bounds are proved in Section \ref{section:2d_DBG}.
\begin{theorem} \label{theorem:2d_DBG}
For each $\theta \in \mathbb{R}$, there exist positive constants $C_{\eqref{equation:2d_DBG}},C'_{\eqref{equation:2d_DBG}}$ and an integer $N_{0} \geq 3$, all depending only on $\theta$, such that
\begin{align} \label{equation:2d_DBG}
    -\exp(C_{\eqref{equation:2d_DBG}} N \log N) \leq \inf \mathscr{H}_{N}^{\theta} \leq -\exp(C'_{\eqref{equation:2d_DBG}} N) \quad \forall N\geq N_{0}.
\end{align}
The same bounds hold for $\mathscr{H}_{N}^{\theta;\textup{sym}}$. In both cases, the lower bound holds for every $N\geq 3$.
\end{theorem}

Using a mean-field approximation, Rajeev \cite{rajeev1999condensation} conjectured that the ground-state energy of the bosonic Hamiltonian has the asymptotic form
\begin{conjecturedequations}
\begin{align} \label{conjecture:Rajeev1999}
    \inf \mathscr{H}_{N}^{\theta;\textup{sym}} = -\frac{\exp(4\pi C_{4-\textup{particles}}N)}{N}  \bkt{C_{3-\textup{particles}}+O(1/N) } \quad \text{as} \quad N \to \infty,
\end{align}
\end{conjecturedequations}
where $C_{3-\textup{particles}}$ and $C_{4-\textup{particles}}$ are positive constants independent of $N$. Throughout this paper, \textit{an asterisk in an equation number} indicates a conjectural statement. We refer to Section \ref{section:discussion_Rajeev_ conjecture} for further discussion of this conjecture, and to \cite{PhysRevLett.93.250408,Bazak_2018} for related predictions and numerical evidence concerning the exponential growth rate of the magnitude of $\inf \mathscr{H}_{N}^{\theta;\textup{sym}}$.

Using \eqref{equation:rewrite_exponent}, Theorem \ref{theorem:2d_DBG} implies
\begin{align*}
    -\exp(C_{\eqref{equation:2d_DBG}} N^{1+o(1)}) \leq \inf \mathscr{H}_{N}^{\theta;\textup{sym}} \leq -\exp(C'_{\eqref{equation:2d_DBG}} N).
\end{align*}
Thus, our estimate matches the conjectured power of $N$ in the exponent up to an $o(1)$ error as $N\to\infty$. We expect that the logarithmic factor in our bound can be removed, yielding a lower bound of the form
\begin{align*}
    -\exp(CN).
\end{align*}
 We refer to Section \ref{section:discussion_Rajeev_ conjecture} for a discussion of the conjectured optimal value of $C$. In one dimension, by contrast, the magnitude of the ground-state energy grows cubically in $N$; see \cite{10.1063/1.3663431} and, for the corresponding SHE moment asymptotics, \cite{AIHPB_2015__51_4_1486_0}.
Our model is related to the attractive mean-field models studied in \cite{JUNGE2026103800,10.1063/1.5131320,lewin:hal-01700668}, although we do not discuss this connection in detail.

The proof of Theorem \ref{theorem:2d_DBG} combines the semigroup-moment identity \eqref{equation:connection_SHF_DBG} with Theorem \ref{theorem:main_result_higher_moment} and the localized lower moment bound in Theorem \ref{theorem:localization}. To prove the lower bound in \eqref{equation:2d_DBG}, we use the
spectral theorem to recover the bottom of the spectrum from the large-time behavior of the semigroup tested against a heat kernel; see Lemma \ref{lemma:spectrum_theorem_1}. This circumvents the difficulty that our upper moment bounds hold only for a limited class of test functions and therefore do not directly control the operator norm of
$\exp(-t\mathscr{H}^{\theta}_{N})$. The large-time argument is our main reason for studying the dependence on $t$ in \eqref{equation:main_result_higher_moment}; see Section \ref{section:lower_bounds}.

\subsubsection{\textbf{Ground-state energy via the principal operator}} \label{section:discussion_Rajeev_ conjecture}
To investigate the constant $C_{4-\textup{particles}}$ in \eqref{conjecture:Rajeev1999}, we derive in
Section \ref{section:Ground_state_energy_and_mean_field_approximation} a formula for the ground-state energy
of the bosonic Hamiltonian in terms of the bottom of the spectrum of the
principal operator $\mathcal{A}_{-1}^{N;\theta}$ defined in
\eqref{equation:resolvent_representation_1}.

\begin{theorem} \label{theorem:informal_formula_ground_state_energy}
For each $\theta\in \mathbb{R}$ and $N\geq 3$, $\mathcal{A}_{-1}^{N;\theta}$ is a self-adjoint operator bounded from below and
\begin{align} \label{equation:informal_formula_ground_state_energy}
    \inf \mathscr{H}_{N}^{\theta;\textup{sym}} = -\exp(4\pi (-\inf \mathcal{A}_{-1}^{N;\theta})).
\end{align}
\end{theorem}

Combining Theorem \ref{theorem:2d_DBG} and Theorem \ref{theorem:informal_formula_ground_state_energy} gives
\begin{align} \label{A_-1_near_linear}
    c N\leq -\inf \mathcal{A}_{-1}^{N;\theta} \leq C N^{1+o(1)}
\end{align}
for all sufficiently large $N$.  The proof of \eqref{A_-1_near_linear} relies on the diagrammatic path-integral representation  \eqref{definition:semigroup} for the semigroup.  However, this representation  does not appear to provide a direct way to identify the constant $C_{4-\textup{particles}}$ in~\eqref{conjecture:Rajeev1999}. By contrast, the formulation \eqref{equation:informal_formula_ground_state_energy} provides a natural setting for investigating this coefficient.

The key idea in deriving
\eqref{equation:informal_formula_ground_state_energy} is to analytically
continue the resolvent representation \eqref{equation:resolvent_representation}
of Gu, Quastel, and Tsai \cite[Proposition 1.4]{Gu_2d} to a domain containing
$(-\infty,\inf \mathscr{H}_{N}^{\theta;\textup{sym}})$;
see Theorem \ref{theorem:resolvent_extension} for the precise statement
and Figure \ref{Figure:1} for an illustration of the extended domain.
In particular, the extended formula shows that the divergence of the
resolvent norm as $z$ approaches
$\inf\mathscr{H}_{N}^{\theta;\textup{sym}}$ from below forces the
operator norm of $(\mathcal{A}_{z}^{N;\theta})^{-1}$ to diverge as well.
Using this observation, we show that the bottom of the spectrum of
$\mathcal{A}_{z}^{N;\theta}$ tends to zero in this limit, yielding
\eqref{equation:equation_for_ground_state_energy}.
Rescaling then gives the ground-state energy formula
\eqref{equation:informal_formula_ground_state_energy}. An overview of the argument is given in
Section \ref{section:AC_GSE_formula};
see also \cite[Eq. (90)]{JDimock_2004} for an earlier resolvent formulation.

For the unrestricted Hamiltonian $\mathscr{H}_{N}^{\theta}$, the resolvent representation \eqref{equation:resolvent_representation} does not have the same simplified form. We do not pursue the corresponding ground-state energy formula here. 
\smallskip

We next study Rajeev's mean-field approximation for
$\inf\mathcal{A}_{-1}^{N;\theta}$, which is expected to capture
its leading-order behavior; see
\cite[pp.~9--10]{rajeev1999condensation}.
In our notation, this expectation takes the form
\begin{conjecturedequations}
\begin{align} \label{equation:Raj_idea}
    \inf\mathcal{A}_{-1}^{N;\theta}
    =
    E_{N,\theta}^{\textup{M-F}}+o(N)
    \qquad\text{as }N\to\infty.
\end{align}
\end{conjecturedequations}
Here $E_{N,\theta}^{\textup{M-F}}$ is obtained by restricting the
variational representation of $\inf\mathcal{A}_{-1}^{N;\theta}$
to a class of product functions; see \eqref{equation:Raj_energy}. Theorem \ref{proposition:mean_field_two_side}, stated in Section \ref{section:mean_field_approximation} and proved in Section \ref{section:mean_field_exponent}, shows that, for each fixed $\theta\in\mathbb{R}$,
\begin{align} \label{equation:lare_N_limit_E_M_F}
    \lim_{N\to\infty}\frac{E_{N,\theta}^{\textup{M-F}}}{N}
    =-C_{\textup{G-N}},
\end{align}
where $C_{\textup{G-N}}$ is the optimal constant in the
two-dimensional Gagliardo--Nirenberg inequality
\begin{align} \label{equation:Gagliardo_Nirenberg}
    \|f\|_{L^{4}(\mathbb{R}^{2})}^{4}
    \leq
    C_{\textup{G-N}}
    \|f\|_{L^{2}(\mathbb{R}^{2})}^{2}
    \|\nabla f\|_{L^{2}(\mathbb{R}^{2})}^{2}
    \qquad\forall f\in H^{1}(\mathbb{R}^{2}).
\end{align}
Together with the ground-state energy formula
\eqref{equation:informal_formula_ground_state_energy},
the conjectured approximation \eqref{equation:Raj_idea}
and the limit \eqref{equation:lare_N_limit_E_M_F} would imply
\begin{conjecturedequations}
\begin{align} \label{conjecture:Rajeev1999_1}
    \inf\mathscr{H}_{N}^{\theta;\textup{sym}}
    =
    -\exp\bigl(4\pi C_{\textup{G-N}}N+o(N)\bigr)
    \qquad\text{as }N\to\infty.
\end{align}
\end{conjecturedequations}
This suggests that $C_{4-\textup{particles}}=C_{\textup{G-N}}$. The numerical value of the sharp constant in \cite{cmp/1103922134} gives $\exp(4\pi C_{\textup{G-N}}) \approx 8.567$. 
This value matches the large-$N$ prediction of Hammer and Son
\cite{PhysRevLett.93.250408} for the ratio of successive ground-state
energies and is consistent with the numerical results of Bazak and
Petrov \cite{Bazak_2018}.

\subsection{Outline and notation} \label{section:notations}

The proofs of the main results presented in Section \ref{section:main_result_moment_bound} and Sections \ref{section:main_result_ground_state_estimate}-\ref{section:discussion_Rajeev_ conjecture} are organised as follows.
\begin{enumerate} [label = (\roman*).]
    \item  Section \ref{section:preliminaries} introduces the diagrammatic and spanning-tree representations, together with the spectral properties used below. Section \ref{section:proof_idea_theorem:main_result_higher_moment} presents the proof strategy for Theorem \ref{theorem:main_result_higher_moment}, and Sections \ref{section:SRW_formulation_proof_ideas}-\ref{section:estimate_for_path_functional} establish the required probabilistic and analytic bounds. Section \ref{section:tail_probability_SHF} proves the tail probability bounds in Theorem \ref{corollary:two_side_bound_tail_probability}.
    \item Section \ref{section:2d_DBG} proves the ground-state energy bounds in Theorem \ref{theorem:2d_DBG}.
   
    \item     Section \ref{section:Ground_state_energy_and_mean_field_approximation} presents the precise statements of Theorem \ref{theorem:informal_formula_ground_state_energy} and \eqref{equation:lare_N_limit_E_M_F}, and introduces their proof ideas and three further supporting results: Theorem \ref{theorem:resolvent_extension}, Proposition \ref{theorem:formula_ground_state_energy}, and Theorem \ref{proposition:mean_field_two_side}. These results are proved in  Sections \ref{section:proof_analytic_continuation}-\ref{section:mean_field_exponent}, respectively.
\end{enumerate}

\begin{table}[h]
\centering 
\caption{ \justifying  Frequently used notation}
\begin{tabularx}{\textwidth}{p{0.25\textwidth}Xc} 
\toprule
Symbol & Description & Def. \\ 
\midrule
$G_{t}(x)$ & The transition density of $2d$ Brownian motion & \eqref{equation:density_Brownian_motion} \\ 
$\mathcal{E}_{N}$ & The index set of two-body interactions & \eqref{equation:index_set} \\ 
$\mathcal{E}_{N}^{\otimes' m}$ & Sequences of collision labels with no consecutive repetition & \eqref{equation:non_consecutive_set} \\
$f^{\otimes N}$ & $N$-fold tensor product of $f$ & \eqref{definition:tensor_product}\\
$\inf \mathscr{H}$ &the bottom of the spectrum of $\mathscr{H}$ &\eqref{bottom_spectrum}\\
$\mathbb{C}_{<0}$ & $:= \{z \in \mathbb{C}:\Re(z) < 0\}$ &\\
\bottomrule
\end{tabularx}
\label{table_0}
\end{table}

\section{Preliminaries} \label{section:preliminaries}

We introduce the diagrammatic path-integral and spanning-tree representations and recall the spectral properties used below.

\subsection{Diagrammatic path-integral representation} \label{section:path_integral_representation}

The basic tool to prove Theorem \ref{theorem:main_result_higher_moment} is the diagrammatic path-integral representation of the semigroup $\exp(-t\mathscr{H}_{N}^{\theta})$. Theorem \ref{theorem:semigroup} follows by reformulating the operators in \cite[(1-16)]{Gu_2d} and using the identification of the limiting semigroup established in \cite[p. 217, proof of Theorem 1.1(b)]{Gu_2d}.

\begin{theorem} [Gu, Quastel, and Tsai \cite{Gu_2d}] \label{theorem:semigroup}
Assume that $\theta \in \mathbb{R}$ and $N\geq 3$. For each $\psi,\varphi \in L^{2}(\mathbb{R}^{2N})$ and all $t>0$, it holds that
\begin{align} \label{definition:semigroup}
    \langle \psi, \exp(-t \mathscr{H}_{N}^{\theta}) \varphi \rangle_{L^{2}(\mathbb{R}^{2N})} &=  \langle \psi,G_{t}^{\otimes N}\varphi \rangle_{L^{2}(\mathbb{R}^{2 N})} + \sum_{m=1}^{\infty} \sum_{ \mathbf{i}_{m} \in \mathcal{E}_{N}^{\otimes' m}} \langle \psi, \mathbf{L}_{t}^{\theta;\mathbf{i}_{m}}\varphi \rangle_{L^{2}(\mathbb{R}^{2N})},
\end{align}
The same representation holds for the restricted semigroup $\exp(-t\mathscr{H}_{N}^{\theta;\textup{sym}})$ when $\psi,\varphi \in L^{2}_{\textup{sym}}(\mathbb{R}^{2N})$.
\end{theorem}

We now explain the notation in \eqref{definition:semigroup}. For a measurable function $f:\mathbb{R}^{2} \mapsto \mathbb{R}$, define its $N$-fold tensor product by
\begin{align} \label{definition:tensor_product}
    f^{\otimes N}(z) := \prod_{i=1}^{N} f(z_{i}), \quad z =(z_{1},\ldots,z_{N}) \in \mathbb{R}^{2N}.
\end{align}
Define the set of sequences of collision-pair labels with no consecutive repetition by
\begin{align} \label{equation:non_consecutive_set}
    \mathcal{E}_{N}^{\otimes' m} &:= \set{\mathbf{i}_{m}=((\mathfrak{i}_{n},\mathfrak{i}'_{n}))_{1\leq n\leq m} \in (\mathcal{E}_{N})^{m}: (\mathfrak{i}_{n+1},\mathfrak{i}'_{n+1}) \neq (\mathfrak{i}_{n},\mathfrak{i}'_{n})\quad \forall 1\leq n\leq m-1 },
\end{align}
where $m\geq 1$. For convenience, we use the notation
\begin{align*}
    \begin{bmatrix}
        c_{1} \\
        \vdots \\
        c_{n}
    \end{bmatrix}_{\times}
    := \prod_{j=1}^{n} c_{j}.
\end{align*}
For each $\mathbf{i}_{m}\in\mathcal{E}^{\otimes' m}_{N}$, the path-integral operator $\mathbf{L}^{\theta;\mathbf{i}_{m}}_{t}$ is the integral operator with kernel
\begin{align} \label{equation:path_integral}
    \mathbf{L}_{t}^{\theta;\mathbf{i}_{m}}(z,z')
    &= \int_{b_{0}=0<a_{1}<b_{1}<\ldots<a_{m}<b_{m} < t=a_{m+1}} \prod_{j=1}^{m} da_{j} db_{j} \int_{\mathbb{R}^{2\times 2 m}} \prod_{j=1}^{m} dx_{j} dy_{j} \nonumber\\
    &\prod_{j=1}^{m} \begin{bmatrix}
        G_{a_{j}-b_{\mathrm{p}^{\mathbf{i}_{m}}(j,\mathfrak{i}'_{j})}}(x_{j} - q_{j,\mathfrak{i}'_{j}}) \\
        G_{a_{j}-b_{\mathrm{p}^{\mathbf{i}_{m}}(j,\mathfrak{i}_{j})}}(x_{j} - q_{j,\mathfrak{i}_{j}}) \\
    \end{bmatrix}_{\times}
    \begin{bmatrix}
        G_{(b_{j}-a_{j})/2}(y_{j} - x_{j}) \\
        \mathfrak{s}^{\theta}(b_{j}-a_{j}) \\
    \end{bmatrix}_{\times} \nonumber\\
    &\prod_{k=1}^{N} G_{a_{m+1} - b_{\mathrm{p}^{\mathbf{i}_{m}}(m+1,k)}} (z'_{k} - q_{m+1,k}), 
\end{align}
where
\begin{align*}
    z = (z_{1},\ldots,z_{N})\in \mathbb{R}^{2N}, \quad z'=(z'_{1},\ldots,z'_{N}) \in \mathbb{R}^{2N}
\end{align*}
are the initial and terminal particle positions, respectively. The ancestor map $\mathrm{p}^{\mathbf{i}_{m}}$, the positions $q_{j,\mathfrak{i}}$, and the function $\mathfrak{s}^{\theta}$ are defined below.

The kernel \eqref{equation:path_integral} admits the following diagrammatic interpretation along the timeline:
\begin{align} \label{definition:time}
    0<a_{1}<b_{1}<\ldots<a_{m}<b_{m}<t.
\end{align}
Starting from the initial configuration $z$, particles move freely between collisions. At time $a_{j}$, particles $\mathfrak{i}'_{j}$ and $\mathfrak{i}_{j}$ meet at $x_{j}$, move together during $[a_{j},b_{j}]$, and separate at $y_{j}$ at time $b_{j}$. The two incoming heat kernels in the second line of \eqref{equation:path_integral} describe their free motion to $x_{j}$. Their common motion from $x_{j}$ to $y_{j}$ is represented by $G_{(b_{j}-a_{j})/2}(y_{j}-x_{j})$, with the additional weight $\mathfrak{s}^{\theta}(b_{j}-a_{j})$.

To specify the starting points in the two incoming heat kernels, we record the most recent collision involving the corresponding particle.
For each  $1\leq j\leq m+1$ and $\mathfrak{i} \in [N]$,  define the ancestor map by
\begin{align} \label{definition:map_p}
    \mathrm{p}^{\mathbf{i}_{m}}(j,\mathfrak{i}) := \max (\{0\}\cup \{1\leq k\leq j-1:\mathfrak{i} \in  \{\mathfrak{i}'_{k},\mathfrak{i}_{k}\}  \}).
\end{align}
Thus, $\mathrm{p}^{\mathbf{i}_{m}}(j,\mathfrak{i})$ is the index of the most recent collision involving particle $\mathfrak{i}$ before the $j$-th collision, with value $0$ if no such collision has occurred. Set
\begin{align*}
    q_{j,\mathfrak{i}} := 
    \begin{cases}
        y_{\mathrm{p}^{\mathbf{i}_{m}}(j,\mathfrak{i})}, \quad &\text{if } \mathrm{p}^{\mathbf{i}_{m}}(j,\mathfrak{i}) \geq 1;\\
        z_{\mathfrak{i}}, \quad &\text{if } \mathrm{p}^{\mathbf{i}_{m}}(j,\mathfrak{i}) = 0.
    \end{cases}
\end{align*}
When $\mathrm{p}^{\mathbf{i}_{m}}(j,\mathfrak{i})\geq1$, the pair $(q_{j,\mathfrak{i}},b_{\mathrm{p}^{\mathbf{i}_{m}}(j,\mathfrak{i})})$ records the position and time at which that previous collision ends. Otherwise, it is the initial space-time pair $(z_{\mathfrak{i}},0)$. 

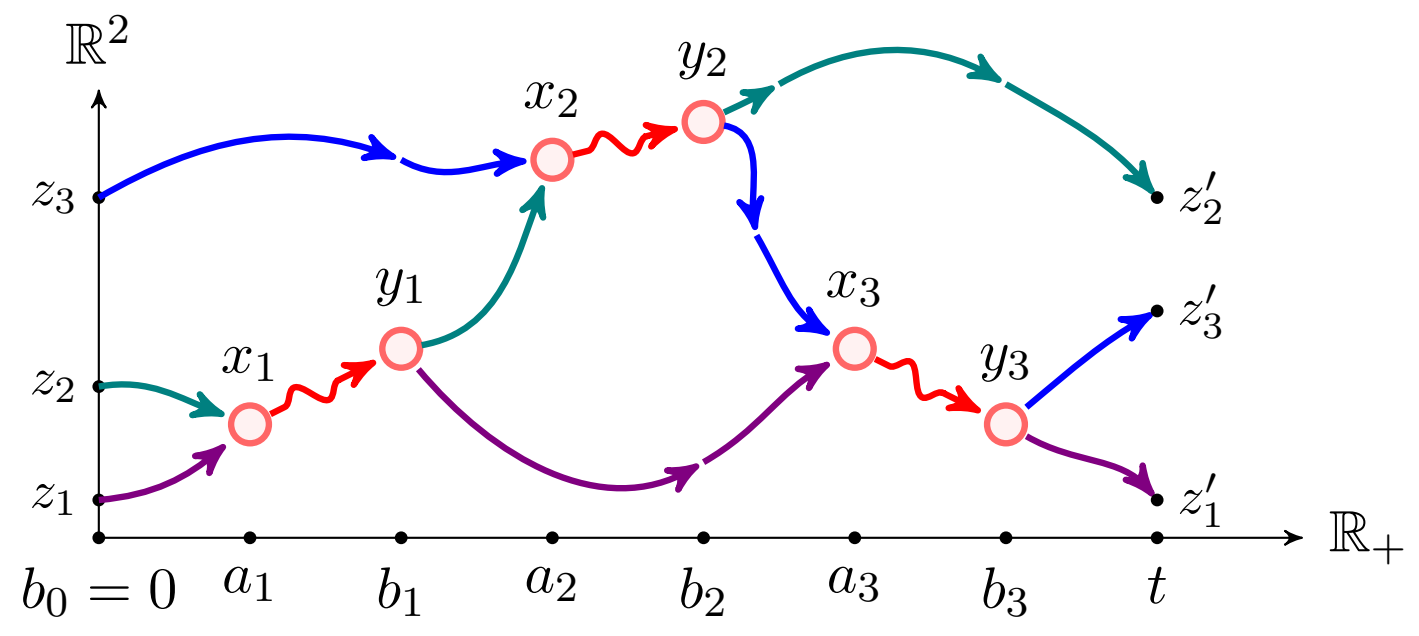
\begin{figure}
\centering
\begin{tikzpicture} [
        scale = 0.5,
        ->,    > = stealth',
       shorten > = 1pt,auto,
   node distance = 3cm,
      decoration = {snake,   
                    pre length=3pt,post length=7pt,
                    }
                    ]

\path[black] (-2,-0.12) node[anchor=north]{$b_{0} = 0$};

\path[black] (0,-0.12) node[anchor=north]{$a_{1}$};
\path[black] (2,-0.12) node[anchor=north]{$b_{1}$};
\path[black] (4,-0.12) node[anchor=north]{$a_{2}$};
\path[black] (6,-0.12) node[anchor=north]{$b_{2}$};
\path[black] (8,-0.12) node[anchor=north]{$a_{3}$};
\path[black] (10,-0.12) node[anchor=north]{$b_{3}$};
\path[black] (12,-0.12) node[anchor=north]{$t$};

\filldraw[black] (-2,0) circle (2pt);
\filldraw[black] (0,0) circle (2pt);
\filldraw[black] (2,0) circle (2pt);
\filldraw[black] (4,0) circle (2pt);
\filldraw[black] (6,0) circle (2pt);
\filldraw[black] (8,0) circle (2pt);
\filldraw[black] (10,0) circle (2pt);
\filldraw[black] (12,0) circle (2pt);

\draw[->] (-2,0) -- (14,0) node[anchor=west]{$\mathbb{R}_{+}$};
\draw[->] (-2,0) -- (-2,6) node[anchor=south]{$\mathbb{R}^{2}$};

\filldraw[black] (-2,0.5) circle (2pt) ;
\filldraw[black] (-2,2) circle (2pt);
\filldraw[black] (-2,4.5) circle (2pt);

\node (b) at (-2,0.5) {};
\node (c) at (-2,2) {};
\node (d) at (-2,4.5) {};

\node at (b)[anchor=east]{\small $z_{1}$};
\node at (c)[anchor=east]{\small $z_{2}$};
\node at (d)[anchor=east]{\small $z_{3}$};

\node (A) at (0,1.5) {};
\filldraw[color=red!60, fill=red!5, very thick](0,1.5) circle (0.25);
\node (B) at (2,2.5) {};
\filldraw[color=red!60, fill=red!5, very thick](2,2.5) circle (0.25);
\node (C) at (4,5) {};
\filldraw[color=red!60, fill=red!5, very thick](4,5) circle (0.25);
\node (D) at (8,2.5) {};
\filldraw[color=red!60, fill=red!5, very thick](8,2.5) circle (0.25);
\node (E) at (10,1.5) {};
\filldraw[color=red!60, fill=red!5, very thick](10,1.5) circle (0.25);
\node (CC) at (6,5.5) {};
\filldraw[color=red!60, fill=red!5, very thick](6,5.5) circle (0.25);

\draw[->] [violet,very thick] (-2,0.5) to[out=5,in=-140] (A);
\draw[->] [teal,very thick] (-2,2) to[out=10,in=160] (A);

\path[draw=red, very thick, decorate] (A) -- (B);
\path[draw=red, very thick, decorate] (C) -- (CC);
\path[draw=red, very thick, decorate] (D) -- (E);

\draw[->] [blue,very thick] (-2,4.5) to[out=30,in=160] (2,5);

\draw[->] [blue,very thick] (2,5) to[out=-30,in=-177] (C);
\draw[->] [teal,very thick] (B) to[out=10,in=-110] (C);
\draw[->] [violet,very thick] (B) to[out=-50,in=-150] (6,1);
\draw[->] [violet,very thick] (6,1) to[out=30,in=-150] (D);
\draw[->] [blue,very thick] (CC) to[out=-10,in=100] (6.7,4);
\draw[->] [blue,very thick] (6.7,4) to[out=-60,in=150] (D);
\draw[->] [teal,very thick] (CC) to[out=25,in=-150] (7,6);
\draw[->] [teal,very thick] (7,6) to[out=30,in=150] (10,6);
\draw[->] [teal,very thick] (10,6) to[out=-30,in=130] (12,4.5);
\draw[->] [blue,very thick] (E) to[out=40,in=-150] (12,3);
\draw[->] [violet,very thick] (E) to[out=-30,in=130] (12,0.5);

\node at (12,0.5) [anchor=west]{\small $z'_{1}$};
\node at (12,3) [anchor=west]{\small $z'_{3}$};
\node at (12,4.5) [anchor=west]{\small $z'_{2}$};

\filldraw[black] (12,0.5) circle (2pt) ;
\filldraw[black] (12,3) circle (2pt);

\filldraw[black] (12,4.5) circle (2pt);

\node at (0,1.8)[anchor=south]{$x_{1}$};

\node at (2,2.8)[anchor=south]{$y_{1}$};

\node at (4,5.3)[anchor=south]{$x_{2}$};

\node at (8,2.8)[anchor=south]{$x_{3}$};

\node at (10,1.8)[anchor=south]{$y_{3}$};

\node at (6,5.8)[anchor=south]{$y_{2}$};

\end{tikzpicture}
\begin{minipage}{0.9\textwidth}
\caption{
\justifying \footnotesize
The Feynman diagram corresponds to $\mathbf{L}_{t}^{\theta;\mathbf{i}_{3}}(z,z')$, where $(\mathfrak{i}_{1},\mathfrak{i}'_{1}) = (1,2)$, $(\mathfrak{i}_{2},\mathfrak{i}'_{2}) = (2,3)$, $(\mathfrak{i}_{3},\mathfrak{i}'_{3}) = (1,3)$, and $N = 3$. Here the first, second, and third particles are colored purple, green, and blue, respectively. In this example, the ancestor map is given by $\mathrm{p}^{\mathbf{i}_{3}}(3,3) = 2$, $\mathrm{p}^{\mathbf{i}_{3}}(3,1) = 1$, $\mathrm{p}^{\mathbf{i}_{3}}(2,3) = 0$, $\mathrm{p}^{\mathbf{i}_{3}}(2,2) = 1$, and $\mathrm{p}^{\mathbf{i}_{3}}(1,2) = \mathrm{p}^{\mathbf{i}_{3}}(1,1) = 0$. Hence, $q_{3,3} = y_{2}$, $q_{3,1} = y_{1}$, $q_{2,3} = z_{3}$, $q_{2,2} = y_{1}$, $q_{1,2} = z_{2}$, and $q_{1,1} = z_{1}$, respectively. 
}  \label{Figure:Feynman_0}
\end{minipage}
\end{figure}

The Volterra-type function $\mathfrak{s}^{\theta}: \mathbb{R}_{>0} \mapsto \mathbb{R}_{>0}$ appearing in \eqref{definition:semigroup} is defined by
\begin{align*}
    \mathfrak{s}^{\theta}(t) := (4\pi) \int_{0}^{\infty} \frac{t^{u-1} \exp(\theta u)}{\Gamma(u)} du.
\end{align*}
The product of heat kernels in the third line of \eqref{equation:path_integral} describes the final free motion of all particles to their terminal positions $z'$ at time $t$. Each sequence $\mathbf{i}_{m}\in\mathcal{E}^{\otimes' m}_{N}$ can be encoded by a Feynman diagram; see Figure \ref{Figure:Feynman_0} for an example.

\subsection{Spanning-tree representation} \label{section:matrix_tree_theorem}

We now recall the spanning-tree representation introduced by Ganguly and Nam \cite{ganguly2025sharp} for the higher moments of the SHF with constant initial datum and a Gaussian test function. More precisely, we consider
\begin{align} \label{equation:Gaussian_test_function}
    \mathbb{E}[\mathscr{Z}^{\theta}_{t}(G_{\nu})^{N}], \quad \mathscr{Z}^{\theta}_{t}(G_{\nu}):= \mathscr{Z}^{\theta}_{t}(\mathds{1},G_{\nu}),
\end{align}
where $\nu > 0$. We first specialize the diagrammatic path-integral representation to this choice of initial datum and test function. Lemma \ref{lemma:Gaussian_text_function} follows directly from \cite[Theorem 2.3]{chen}, so we omit the proof.
\begin{lemma} \label{lemma:Gaussian_text_function}
Assume that $\theta \in \mathbb{R}$ and $N\geq 3$. For each $t,\nu > 0$, 
\begin{align*}
    \exp(-t\mathscr{H}_{N}^{\theta})G_{\nu}^{\otimes N} \in L^{1}(\mathbb{R}^{2N})
\end{align*}
and
\begin{align} \label{equation:path_integral_2}
    \mathbb{E}[\mathscr{Z}^{\theta}_{t}(G_{\nu})^{N}] = \int_{\mathbb{R}^{2N}} dx' \exp(-t\mathscr{H}_{N}^{\theta})G_{\nu}^{\otimes N}(x')  =  1 + \sum_{m=1}^{\infty} \sum_{ \mathbf{i}_{m} \in \mathcal{E}_{N}^{\otimes' m}} \int_{\mathbb{R}^{2N}} dx' \mathbf{L}_{t}^{\theta;\textbf{i}_{m}}G_{\nu}^{\otimes N}(x').
\end{align}
For each $m \geq 1$ and $\mathbf{i}_{m} \in \mathcal{E}_{N}^{\otimes' m}$, 
\begin{align} \label{equation:path_integral_1}
    \int_{\mathbb{R}^{2N}} dx' \mathbf{L}_{t}^{\theta;\textbf{i}_{m}}G_{\nu}^{\otimes N}(x')
    &= \int_{\nu <a_{1}<b_{1}<\ldots<a_{m}<b_{m} < t+\nu} \prod_{j=1}^{m} da_{j} db_{j} \int_{\mathbb{R}^{2\times 2 m}} \prod_{j=1}^{m} dx_{j} dy_{j} \nonumber\\
    &\prod_{j=1}^{m} \begin{bmatrix}
        G_{a_{j}-b_{\mathrm{p}^{\mathbf{i}_{m}}(j,\mathfrak{i}'_{j})}}(x_{j} - y_{\mathrm{p}^{\mathbf{i}_{m}}(j,\mathfrak{i}'_{j})}) \\
        G_{a_{j}-b_{\mathrm{p}^{\mathbf{i}_{m}}(j,\mathfrak{i}_{j})}}(x_{j} - y_{\mathrm{p}^{\mathbf{i}_{m}}(j,\mathfrak{i}_{j})}) \\
    \end{bmatrix}_{\times}
    \begin{bmatrix}
        G_{(b_{j}-a_{j})/2}(y_{j} - x_{j}) \\
        \mathfrak{s}^{\theta}(b_{j}-a_{j}) \\
    \end{bmatrix}_{\times}.
\end{align}
Here the ancestor map $\mathrm{p}^{\mathbf{i}_{m}}$ is defined in \eqref{definition:map_p}, and we set
\begin{align*}
    b_{0} = 0, \quad y_{0} = (0,0) \in \mathbb{R}^{2}.
\end{align*}

\end{lemma}

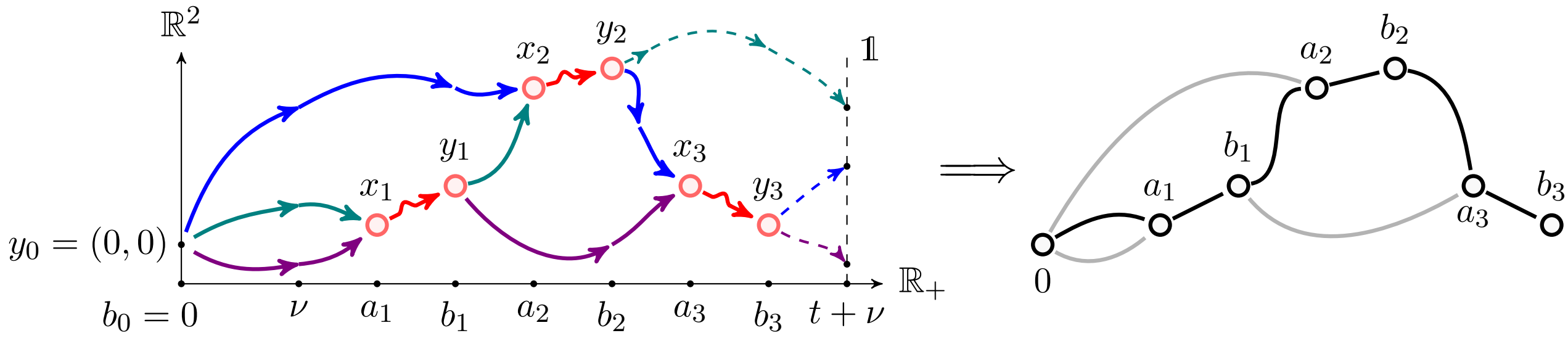
\begin{figure}
\centering
\begin{tikzpicture} [
        scale = 0.42,
        ->,    > = stealth',
       shorten > = 1pt,auto,
       auto,
   node distance = 3cm,
      decoration = {snake,   
                    pre length=3pt,post length=7pt,
                    }
                    ]

\path[black] (-5.8,-0.12) node[anchor=north]{$b_{0} = 0$};
\path[black] (-2,-0.12) node[anchor=north]{$\nu$};
\path[black] (0,-0.12) node[anchor=north]{$a_{1}$};
\path[black] (2,-0.12) node[anchor=north]{$b_{1}$};
\path[black] (4,-0.12) node[anchor=north]{$a_{2}$};
\path[black] (6,-0.12) node[anchor=north]{$b_{2}$};
\path[black] (8,-0.12) node[anchor=north]{$a_{3}$};
\path[black] (10,-0.12) node[anchor=north]{$b_{3}$};
\path[black] (12,-0.12) node[anchor=north]{$t+\nu$};

\filldraw[black] (-5,0) circle (2pt);
\filldraw[black] (0,0) circle (2pt);
\filldraw[black] (-2,0) circle (2pt);
\filldraw[black] (2,0) circle (2pt);
\filldraw[black] (4,0) circle (2pt);
\filldraw[black] (6,0) circle (2pt);
\filldraw[black] (8,0) circle (2pt);
\filldraw[black] (10,0) circle (2pt);
\filldraw[black] (12,0) circle (2pt);

\draw[->] (-5,0) -- (13,0) node[anchor=west]{$\mathbb{R}_{+}$};
\draw[->] (-5,0) -- (-5,6) node[anchor=south]{$\mathbb{R}^{2}$};

\draw[dashed,-] (12,0) -- (12,6) node[anchor=west]{$\mathds{1}$};

\node (a) at (-5,1) {};
\filldraw[black] (a) circle (2pt)
node[anchor=east]{ $y_{0} = (0,0)$};

\coordinate (b) at (-2,0.5);
\coordinate (c) at (-2,2);
\coordinate (d) at (-2,4.5);

\draw[->, violet, very thick, shorten >=0pt]
    (a) to[out=-30,in=-170] (b);

\draw[->, teal, very thick, shorten >=0pt]
    (a) to[out=30,in=-170] (c);

\draw[->, blue, very thick, shorten >=0pt]
    (a) to[out=70,in=-160] (d);

\node (A) at (0,1.5) {};
\filldraw[color=red!60, fill=red!5, very thick](0,1.5) circle (0.25);
\node (B) at (2,2.5) {};
\filldraw[color=red!60, fill=red!5, very thick](2,2.5) circle (0.25);
\node (C) at (4,5) {};
\filldraw[color=red!60, fill=red!5, very thick](4,5) circle (0.25);
\node (D) at (8,2.5) {};
\filldraw[color=red!60, fill=red!5, very thick](8,2.5) circle (0.25);
\node (E) at (10,1.5) {};
\filldraw[color=red!60, fill=red!5, very thick](10,1.5) circle (0.25);
\node (CC) at (6,5.5) {};
\filldraw[color=red!60, fill=red!5, very thick](6,5.5) circle (0.25);

\draw[->, very thick, violet]
    (-2,0.5) to[out=5,in=-140] (A);

\draw[->, very thick, teal]
    (-2,2) to[out=10,in=160] (A);

\draw[->, very thick, blue]
    (-2,4.5) to[out=30,in=160] (2,5);

\path[draw=red, very thick, decorate] (A) -- (B);
\path[draw=red, very thick, decorate] (C) -- (CC);
\path[draw=red, very thick, decorate] (D) -- (E);

\draw[->] [very thick,blue] (2,5) to[out=-30,in=-175] (C);
\draw[->] [very thick,teal] (B) to[out=10,in=-110] (C);
\draw[->] [very thick,violet] (B) to[out=-50,in=-150] (6,1);
\draw[->] [very thick,violet] (6,1) to[out=30,in=-150] (D);
\draw[->] [very thick,blue] (CC) to[out=-10,in=100] (6.7,4);
\draw[->] [very thick,blue] (6.7,4) to[out=-60,in=150] (D);
\draw[dashed,->] [thick,teal] (CC) to[out=25,in=-150] (7,6);
\draw[dashed,->] [thick,teal] (7,6) to[out=30,in=150] (10,6);
\draw[dashed,->] [thick,teal] (10,6) to[out=-30,in=130] (12,4.5);
\draw[dashed,->] [thick,blue] (E) to[out=40,in=-150] (12,3);
\draw[dashed,->] [thick,violet] (E) to[out=-30,in=130] (12,0.5);

\filldraw[black] (12,0.5) circle (2pt) ;
\filldraw[black] (12,3) circle (2pt);

\filldraw[black] (12,4.5) circle (2pt);

\node at (0,1.8)[anchor=south]{$x_{1}$};

\node at (2,2.8)[anchor=south]{$y_{1}$};

\node at (4,5.3)[anchor=south]{$x_{2}$};

\node at (8,2.8)[anchor=south]{$x_{3}$};

\node at (10,1.8)[anchor=south]{$y_{3}$};

\node at (6,5.8)[anchor=south]{$y_{2}$};

\node[anchor=west] at (14,3) {\Large $\xLongrightarrow{\quad}$};

\node (a') at (17,1) {};
\node at (17,0.7) [anchor=north]{$0$};
\filldraw[color=black, fill=gray!5, very thick](17,1) circle (0.25);

\coordinate (b') at (18,0.5);
\coordinate (c') at (18,2);
\coordinate (d') at (18,4.5);

\node (A') at (20,1.5) {};
\filldraw[color=black, fill=gray!5, very thick](20,1.5) circle (0.25);
\node (B') at (22,2.5) {};
\filldraw[color=black, fill=gray!5, very thick](22,2.5) circle (0.25);
\node (C') at (24,5) {};
\filldraw[color=black, fill=gray!5, very thick](24,5) circle (0.25);
\node (D') at (28,2.5) {};
\filldraw[color=black, fill=gray!5, very thick](28,2.5) circle (0.25);
\node (E') at (30,1.5) {};
\filldraw[color=black, fill=gray!5, very thick](30,1.5) circle (0.25);
\node (CC') at (26,5.5) {};
\filldraw[color=black, fill=gray!5, very thick](26,5.5) circle (0.25);

\draw[very thick, gray!60,-]
    (a') to[out=-30,in=-140] (A');

\draw[very thick, black,-]
    (a') to[out=30,in=160] (A');

\draw[very thick, gray!60,-]
    (a') to[out=60,in=160] (C');

\path[draw=black, very thick,-] (A') -- (B');
\path[draw=black, very thick,-] (C') -- (CC');
\path[draw=black, very thick,-] (D') -- (E');

\draw[-] [very thick,black] (B') to[out=10,in=-177] (C');
\draw[-] [very thick,gray!60] (B') to[out=-50,in=-150] (D');
\draw[-] [very thick,black] (CC') to[out=-10,in=100] (D');

\node at (20,1.8)[anchor=south]{$a_{1}$};

\node at (22,2.8)[anchor=south]{$b_{1}$};

\node at (24,5.3)[anchor=south]{$a_{2}$};

\node at (28,2.3)[anchor=north]{$a_{3}$};

\node at (30,1.8)[anchor=south]{$b_{3}$};

\node at (26,5.8)[anchor=south]{$b_{2}$};

\end{tikzpicture}
\begin{minipage}{0.9\textwidth}
\caption{ \justifying \footnotesize
The left figure is the Feynman diagram corresponds to $\int_{\mathbb{R}^{2\times 3}} dx'\mathbf{L}_{t}^{\theta;\textbf{i}_{3}}G_{\nu}^{\otimes 3}(x') $, where the setting is the same as in Figure \ref{Figure:Feynman_0}. Because of the Gaussian test function, all the particles start from $0 \in \mathbb{R}^{2}$ and the first collision time is at least $\nu$. In particular, since the initial condition is $\mathds{1}$, the dashed lines are deleted.  The right figure is the undirected weighted graph $G$ corresponding to the left Feynman diagram.  The black edges form a quasi-shortest spanning tree of $\mathrm{G}$.
}  \label{Figure:Feynman}
\end{minipage}
\end{figure}

In the diagrammatic interpretation described after \eqref{definition:time}, the integral in \eqref{equation:path_integral_1} corresponds to particles starting from the configuration $(0,\ldots,0)\in\mathbb{R}^{2N}$, with the first collision occurring after time $\nu$. The free motions following each particle’s last collision integrate to $1$ and are therefore omitted from the diagram. Thus, each particle path is drawn only up to the end of the last collision in which it participates. Particles that do not participate in any collision contribute a factor of $1$ and are omitted entirely. See Figure \ref{Figure:Feynman} for an example of the associated Feynman diagram.

To compute the spatial integral on the right-hand side of \eqref{equation:path_integral_1},  we now introduce the partition function formula of the Gaussian free field established in \cite[Lemma 3.2]{ganguly2025sharp}.
\begin{lemma}
[Ganguly and Nam \cite{ganguly2025sharp}] \label{lemma:GFF}
Let $H = (V,E,W)$ be a finite weighted connected graph with an arbitrary designated pinning vertex $o \in V$ and positive weight function $W$ defined on $E$. Set $X_{o} := (0,0) \in \mathbb{R}^{2}$.
Then it holds that
\begin{align*}
    \int_{(\mathbb{R}^{2})^{V\setminus\{o\}}} \exp\bkt{-\frac{1}{2} \sum_{(a,b) \in E} W((a,b)) |X_{a} - X_{b}|^{2} } \prod_{a \in V\setminus\{o\}} dX_{a} = (2\pi)^{|V|-1} \Bkt{\sum_{T \in \mathcal{T}(H)} W(T)}^{-1}.
\end{align*}
Here $\mathcal{T}(H)$ denotes the collection of all spanning trees of $H$ and 
\begin{align*}
    W(T) := \prod_{e \in E(T)} W(e).
\end{align*}
Multi-edges are counted separately  in the collection of spanning trees.
\end{lemma}

Applying Lemma \ref{lemma:GFF} to integrate out $dx_{j}dy_{j}$ on the right-hand side of \eqref{equation:path_integral_1} gives the spanning-tree representation of the path integral.
\begin{lemma} \label{lemma:representation_path_integral}
Assume that $\theta \in \mathbb{R}$ and $N\geq 3$. For each $t,\nu > 0$, $m\geq 1$, and $\mathbf{i}_{m} = ((\mathfrak{i}_{n},\mathfrak{i}'_{n}))_{1\leq n\leq m} \in \mathcal{E}_{N}^{\otimes' m}$, we have
\begin{align} \label{equation:representation_path_integral}
    \int_{\mathbb{R}^{2N}} dx'\mathbf{L}_{t}^{\theta;\textbf{i}_{m}}G_{\nu}^{\otimes N} (x') = &C^{m}\int_{\nu<a_{1}<b_{1}<\ldots<a_{m}<b_{m} < t+\nu} \prod_{j=1}^{m} da_{j} db_{j} \nonumber\\
    &\prod_{j=1}^{m} \begin{bmatrix}
        \displaystyle \frac{1}{(b_{j}-a_{j})} \\
        \displaystyle \mathfrak{s}^{\theta}(b_{j}-a_{j}) \\
    \end{bmatrix}_{\times} 
    \begin{bmatrix}
        \displaystyle \frac{1}{a_{j}-b_{\mathrm{p}^{\mathbf{i}_{m}}(j,\mathfrak{i}'_{j})}}\\
        \displaystyle \frac{1}{a_{j}-b_{\mathrm{p}^{\mathbf{i}_{m}}(j,\mathfrak{i}_{j})}}\\
    \end{bmatrix}_{\times} \Bkt{\sum_{T \in \mathcal{T}(\mathrm{G})} W(T) }^{-1},
\end{align}
where $C$ is an absolute positive constant and  $\mathrm{G} =(V,E,W)$ is the undirected weighted multigraph defined below.
\end{lemma}
For each choice of the temporal variables in \eqref{equation:representation_path_integral}, the vertex set consists of time $0$ and the endpoints of all collision intervals:
\begin{align} \label{equation:graph}
    V &:= \{0\}\cup \{a_{j},b_{j}: 1\leq j\leq m\}.
\end{align}
The terminal time $t+\nu$ is not included. The sequence $\mathbf{i}_{m}\in\mathcal{E}_{N}^{\otimes' m}$ determines the edges through the ancestor map:
\begin{align} \label{equation:graph_1}
    E &:= E_{S} \cup E_{L}, \nonumber\\
    E_{S} &:=  \{(a_{j},b_{j}): 1\leq j\leq m\}, \quad E_{L} := \{(b_{\mathrm{p}^{\mathbf{i}_{m}}(j,\mathfrak{i}'_{j})},a_{j}),(b_{\mathrm{p}^{\mathbf{i}_{m}}(j,\mathfrak{i}_{j})},a_{j}): 1\leq j\leq m\}. 
\end{align}
Here $E_{S}$ represents the collision intervals, while $E_{L}$ represents the free motion  leading to collisions. The edge collections are understood with multiplicity. Multi-edges occur when
\begin{align*}
    \mathrm{p}^{\mathbf{i}_{m}}(j,\mathfrak{i}'_{j}) = \mathrm{p}^{\mathbf{i}_{m}}(j,\mathfrak{i}_{j}),
\end{align*}
in which case the two incoming edges at $a_{j}$ have the same endpoints but are retained as distinct edges. See Figure \ref{Figure:Feynman} for an example.

Each edge weight is inversely proportional to the length of its corresponding time interval. More precisely, for $e=(t_{1},t_{2})$ with $t_{1}<t_{2}$, set
\begin{align} \label{equation:weight_function}
    W(e) := 
    \begin{cases}
        \displaystyle \frac{1}{t_{2}-t_{1}}, \quad &\text{if } e = (t_{1},t_{2}), \; e \in E_{L};\\
        \displaystyle \frac{2}{t_{2}-t_{1}}, \quad &\text{if } e = (t_{1},t_{2}), \; e \in E_{S}.
    \end{cases}
\end{align}

\begin{proof} [Proof of Lemma \ref{lemma:representation_path_integral}]
Expanding the heat kernels in \eqref{equation:path_integral_1}, we obtain
\begin{align*}
    &\int_{\mathbb{R}^{2N}} dx'\mathbf{L}_{t}^{\theta;\textbf{i}_{m}}G_{\nu}^{\otimes N}(x') \nonumber\\
    &= C_{0}^{m} \int_{\nu<a_{1}<b_{1}<\ldots<a_{m}<b_{m} < t+\nu} \prod_{j=1}^{m} da_{j} db_{j}  \prod_{j=1}^{m} \begin{bmatrix}
        \displaystyle \frac{1}{(b_{j}-a_{j})} \\
        \displaystyle \mathfrak{s}^{\theta}(b_{j}-a_{j}) \\
    \end{bmatrix}_{\times} 
    \begin{bmatrix}
        \displaystyle \frac{1}{a_{j}-b_{\mathrm{p}^{\mathbf{i}_{m}}(j,\mathfrak{i}'_{j})}}\\
        \displaystyle \frac{1}{a_{j}-b_{\mathrm{p}^{\mathbf{i}_{m}}(j,\mathfrak{i}_{j})}}\\
    \end{bmatrix}_{\times}\\
    &\int_{\mathbb{R}^{2\times 2\times m}} \prod_{j=1}^{m} dx_{j} dy_{j} 
    \exp \bkt{-\frac{1}{2}\sum_{j=1}^{m} \frac{|y_{j}-x_{j}|^{2}}{(b_{j}-a_{j})/2} -\frac{1}{2}\sum_{j=1}^{m} \frac{|x_{j} - y_{\mathrm{p}^{\mathbf{i}_{m}}(j,\mathfrak{i}'_{j})}|^{2}}{a_{j}-b_{\mathrm{p}^{\mathbf{i}_{m}}(j,\mathfrak{i}'_{j})}}-\frac{1}{2}\sum_{j=1}^{m}\frac{|x_{j} - y_{\mathrm{p}^{\mathbf{i}_{m}}(j,\mathfrak{i}_{j})}|^{2}}{a_{j}-b_{\mathrm{p}^{\mathbf{i}_{m}}(j,\mathfrak{i}_{j})}} },
\end{align*}
where $C_{0}>0$ is an absolute constant. Identify the spatial variables with the vertices of $\mathrm{G}$ by setting
\begin{align*}
    X_{0}:=(0,0), \quad X_{a_{j}}:=x_{j}, \quad X_{b_{j}}:=y_{j}, \quad 1\leq j\leq m.
\end{align*}
The exponent in the spatial integral then becomes
\begin{align*}
    -\frac{1}{2} \sum_{e=(a,b)\in E}W(e)|X_{a}-X_{b}|^{2}.
\end{align*}
The graph $\mathrm{G}$ is connected, since every vertex is joined to $0$ by a path through earlier collision times. Since $|V|=2m+1$, Lemma \ref{lemma:GFF} evaluates the spatial integral as
\begin{align*}
    (2\pi)^{2m} \Bkt{ \sum_{T\in\mathcal{T}(\mathrm{G})}W(T)}^{-1}
\end{align*}
Thus, \eqref{equation:representation_path_integral} follows with $C=(2\pi)^{2}C_{0}$
\end{proof}

\subsection{Spectral properties of self-adjoint operators}  \label{section:spectrum_properties}

In this subsection, we collect two general consequences of the spectral theorem for self-adjoint operators bounded from above. The first identifies the large-time exponential growth rate of the semigroup  in terms of the support of the associated
spectral measure.

\begin{lemma} \label{lemma:lemma_spectrum_theorem_support}
Let $\mathscr{L}$ be a self-adjoint operator on a Hilbert space $H$ that is bounded from above. For each non-zero vector $\psi \in H$, it holds that
\begin{align} \label{proof:lemma_spectrum_theorem_support}
    \lim_{t\to\infty} \frac{1}{t} \log \langle \psi, \exp(t \mathscr{L}) \psi \rangle_{H} = \sup \textup{supp}(\mu_{\psi}),
\end{align}
where $\mu_{\psi}$ is the spectral measure defined by
\begin{align*}
    \mu_{\psi}(A) :=  \langle \psi, \mathds{1}_{A}(\mathscr{L})\psi \rangle_{H}.
\end{align*}
Here $\mathds{1}_{A}(\mathscr{L})$ denotes the spectral projection of $\mathscr{L}$ associated with the Borel set $A$.
\end{lemma}

The next lemma shows that the top of the spectrum is determined by
the spectral measures associated with any family of vectors whose
linear span is dense.

\begin{lemma} \label{lemma:sup_representation}
Let $\mathscr{L}$ be a self-adjoint operator on a Hilbert space
$H$ that is bounded from above.
Suppose that $\mathcal{G}\subseteq H$ and that
$\textup{span}(\mathcal{G})$ is dense in $H$.
Then
\begin{align} \label{equation:sup_representation}
    \sup \sigma(\mathscr{L}) = \sup_{\psi \in \mathcal{G}} \sup \textup{supp}(\mu_{\psi}).
\end{align}
\end{lemma}
\begin{proof}
Let $\psi\in H$. For each $b>\sup\sigma(\mathscr{L})$, we have
$[b,\infty)\subseteq\rho(\mathscr{L})$, where $\rho(\mathscr{L})$
denotes the resolvent set of $\mathscr{L}$. Hence, \cite[Corollary 3.9]{teschl2009mathematical} gives
\begin{align} \label{equation:sup_representation_1}
    \mu_{\psi}([b,\infty)) = || \mathds{1}_{[b,\infty)}(\mathscr{L})\psi||_{H}^{2} = 0,
\end{align}
which implies $\sup \text{supp}(\mu_{\psi}) \leq b$. Letting $b \to \sup\sigma(\mathscr{L})$, we obtain  $\sup \textup{supp}(\mu_{\psi}) \leq \sup\sigma(\mathscr{L})$ for all $\psi \in \mathcal{G}$. Thus, we have
\begin{align} \label{opposite_inequality_-1}
    \sup_{\psi \in \mathcal{G}} \sup \textup{supp}(\mu_{\psi}) \leq \sup \sigma(\mathscr{L}).
\end{align}

Now let $b> \sup_{\psi \in \mathcal{G}} \sup \textup{supp}(\mu_{\psi})$. Then for each $\psi\in \mathcal{G}$, $\mu_{\psi}([b,\infty)) = 0$, and hence, the first equality in \eqref{equation:sup_representation_1} together with the density of $\textup{span}(\mathcal{G})$ then yields $\mathds{1}_{[b,\infty)}(\mathscr{L}) = 0$. This together with  \cite[Theorem 3.8]{teschl2009mathematical} implies $\sup \sigma(\mathscr{L}) \leq b$. By taking $b \to \sup_{\psi \in \mathcal{G}} \sup \textup{supp}(\mu_{\psi})$, we conclude
\begin{align} \label{opposite_inequality_1}
    \sup \sigma(\mathscr{L}) \leq \sup_{\psi \in \mathcal{G}} \sup \textup{supp}(\mu_{\psi}).
\end{align}
Combining \eqref{opposite_inequality_-1} and \eqref{opposite_inequality_1} implies Lemma \ref{lemma:sup_representation}. 
\end{proof}

\section{Proof ideas for the near-optimal moment bounds} \label{section:proof_idea_theorem:main_result_higher_moment}

To prove Theorem \ref{theorem:main_result_higher_moment}, let $(h,f)$ satisfy \eqref{eq:UBH}. By positivity and the bounds in \eqref{eq:UBH}, it suffices, after adjusting the constant, to consider
\begin{align*}
    h = \mathds{1},\quad f = G_{\nu}, \quad \nu>0.
\end{align*}
Indeed, the hypotheses $h\leq D_{h}$ and $f\leq C_{f}G_{\nu_{f}}$ allow the general moment to be bounded by $(D_{h}C_{f})^{N}$ times the corresponding moment with constant initial datum and a Gaussian test function. We establish the following estimate with an auxiliary parameter $\delta$.
\begin{theorem} \label{theorem:main_result_proof_idea}
Assume that $\theta \in \mathbb{R}$ and $N\geq 5$.  For each $t\geq 1$, $2/N \leq \delta < 2/3$, and $\nu>0$, there exists $C_{\eqref{equation:main_result_proof_idea}} = C_{\eqref{equation:main_result_proof_idea}}(\theta,\nu) > 0$ such that 
\begin{align} \label{equation:main_result_proof_idea}
    \mathbb{E}[\mathscr{Z}_{t}^{\theta}(G_{\nu})^{N}]
    \leq \exp\bkt{t\exp\bkt{\frac{C_{\eqref{theorem:main_result_proof_idea}}}{\delta}N^{1+\delta}} }. 
\end{align}
\end{theorem}

Taking $\delta=1/\log N$ gives $N^{1+\delta}/\delta = eN\log N$. Since this choice satisfies $\delta\geq2/N$ for $N\geq3$, Theorem \ref{theorem:main_result_higher_moment} follows from Theorem \ref{theorem:main_result_proof_idea} after adjusting the constant. 

To prove Theorem \ref{theorem:main_result_proof_idea}, it therefore suffices to bound the series on the right-hand side of \eqref{equation:path_integral_2}. Each path integral
\begin{align*}
    \int_{\mathbb{R}^{2N}} \mathbf{L}^{\theta;\mathbf{i}_{m}}_{t}G_{\nu}^{\otimes N}(x') dx'
\end{align*}
admits the spanning-tree representation \eqref{equation:representation_path_integral}.  In Sections \ref{section:quasi_tree}-\ref{section:collision_time_markov_chain}, we use this representation to derive an upper bound involving a path functional of the collision-time Markov chain; see Proposition \ref{proposition:FK_representation}. The key estimate is a Lyapunov-type inequality for the associated weighted transition operator, stated in Section \ref{section:lya}.

\subsection{A quasi-shortest spanning tree} \label{section:quasi_tree}

To bound \eqref{equation:representation_path_integral}, we retain a single suitably chosen term in the spanning-tree sum. This choice is expected to yield the optimal moment bound $\exp(t\exp(CN))$; see Section \ref{section:near_optimal_tree} for a heuristic justification. For any spanning tree $T_{\star}$ of $\mathrm{G}$, positivity of the weights gives
\begin{align*}
    \Bkt{\sum_{T\in\mathcal{T}(\mathrm{G})}W(T)}^{-1} \leq W(T_{\star})^{-1}.
\end{align*}
We use a particular spanning tree $T_{\star}$, defined below, to obtain the following bound.

\begin{lemma} \label{lemma:QS_spanning_tree}
Assume that $\theta \in \mathbb{R}$ and $N\geq 3$. For each $m \geq 1$ and $\mathbf{i}_{m} \in \mathcal{E}_{N}^{\otimes'm}$, we have
\begin{align} \label{equation:QS_spanning_tree}
    &\int_{\mathbb{R}^{2N}} \mathbf{L}^{\theta;\mathbf{i}_{m}}_{t}G_{\nu}^{\otimes N}(x') dx' \nonumber\\
    &\leq C^{m} \int_{\nu<a_{1}<b_{1}<\ldots<a_{m}<b_{m} < t+\nu} \prod_{j=1}^{m} da_{j} db_{j} \prod_{j=1}^{m} \begin{bmatrix}
        \displaystyle \frac{1}{(b_{j}-a_{j})} \\
        \displaystyle \mathfrak{s}^{\theta}(b_{j}-a_{j}) \\
    \end{bmatrix}_{\times} 
    \begin{bmatrix}
        \displaystyle \frac{1}{a_{j}-b_{\mathrm{p}^{\mathbf{i}_{m}}(j,\mathfrak{i}'_{j})}}\\
        \displaystyle \frac{1}{a_{j}-b_{\mathrm{p}^{\mathbf{i}_{m}}(j,\mathfrak{i}_{j})}}\\
    \end{bmatrix}_{\times} W(T_{\star})^{-1},
\end{align}
where $C>0$ is an absolute constant and $T_{\star}$ is given below.
\end{lemma}

We now describe the construction of $T_{\star}$. For an edge $e$, let $\mathrm{N}(e)$ denote the number of intervals of the form $(b_{k-1},a_{k})$, $1\leq k\leq m$, contained in the time interval corresponding to $e$. We refer to $\mathrm{N}(e)$ as the length of $e$.

\begin{definition} [Quasi-shortest spanning tree]
Given $\mathbf{i}_{m}\in\mathcal{E}_{N}^{\otimes' m}$, let $T_{\star}$ have vertex set $V$, and construct its edge set $E(T_{\star})$ as follows. First, include every edge in $E_{S}$. Next, for each $1\leq j\leq m$, consider the two incoming edges
\begin{align*}
    e_{j}:= \bigl(b_{\mathrm{p}^{\mathbf{i}_{m}}(j,\mathfrak{i}_{j})},a_{j}\bigr), \quad e'_{j}:= \bigl(b_{\mathrm{p}^{\mathbf{i}_{m}}(j,\mathfrak{i}'_{j})},a_{j}\bigr).
\end{align*}
If $\mathrm{N}(e_{j})<\mathrm{N}(e'_{j})$, include $e_{j}$ in $E(T_{\star})$; otherwise, include $e'_{j}$. Thus, we select the shorter incoming edge, choosing $e'_{j}$ when the lengths are equal. Throughout this paper, we call $T_{\star}$ a quasi-shortest spanning tree of $\mathrm{G}$.
\end{definition}
\begin{example}
Consider the configuration in Figure \ref{Figure:Feynman}. The lengths of the edges are given by 
\begin{align*}
    &N((a_{j},b_{j})) = 0, \quad j = 1,2,3, \\
    &N((0,a_{1})) = 1,\\
    &N((0,a_{2})) = 2, \quad N((b_{1},a_{2})) = 1,\\
    &N((b_{2},a_{3})) = 1, \quad N((b_{1},a_{3})) = 2.
\end{align*}
Hence, the spanning tree $T_{\star}$ contains the edges:
\begin{align*}
    (a_{j},b_{j}) \quad \forall j = 1,2,3, \quad (b_{1},a_{2}), \quad (b_{2},a_{3}), 
\end{align*}
and one copy of the multi-edge $(0,a_1)$.
\end{example}

Although the vertices are labelled by the temporal variables $a_{j}$ and $b_{j}$, the graph $(V,E)$, the lengths $\mathrm{N}(e)$, and the choice of $T_{\star}$ depend only on the collision sequence $\mathbf{i}_{m}$, not on the values of these variables. For convenience, we retain $a_{j}$ and $b_{j}$ as vertex labels. The weights $W(e)$, however, still depend on the temporal variables.

We now verify that the construction indeed yields a spanning tree.

\begin{lemma} \label{lemma:tree}
$T_{\star}$ is a spanning tree of $\mathrm{G}$.
\end{lemma} 
\begin{proof}
Every vertex is connected to $0$. Indeed, each $b_{j}$ is joined to $a_{j}$, and the selected incoming edge at $a_{j}$ joins it to either $0$ or some $b_{k}$ with $k<j$. Repeating this procedure produces a path to $0$. Moreover, $T_{\star}$ has $2m+1$ vertices and exactly $2m$ edges: $m$ edges from $E_{S}$ and one edge from $E_{L}$ for each collision. Hence, $T_{\star}$ is connected and has one fewer edge than vertices, so it is a spanning tree of $\mathrm{G}$.
\end{proof}

\subsection{Formulation via a simple random walk}

We now derive a probabilistic upper bound for the $N$-th moment. The requirement that consecutive collision labels be distinct naturally leads to a random walk $(\mathbf{X}_{n})_{n\geq0}$ on the collision-label set $\mathcal{E}_{N}\cup\{\mathbf{a}\}$, where $\mathbf{a}$ is an absorbing state. The initial state $\mathbf{X}_{0}$ is uniformly distributed on $\mathcal{E}_{N}$. From a state $\mathbf{x}\in\mathcal{E}_{N}$, the walk moves to a uniformly chosen state of $(\mathcal{E}_{N} \setminus \{\mathbf{x}\}) \cup \{\mathbf{a}\}$.
Once it reaches $\mathbf{a}$, it remains there. Denote the absorption time by
\begin{align*}
    T_{\mathbf{a}} := \inf\{n\geq1:\mathbf{X}_{n}=\mathbf{a}\}.
\end{align*}
Before absorption, write
\begin{align*}
    \mathbf{X}_{n} = \begin{pmatrix} \mathbf{X}_{n}(1)\\ \mathbf{X}_{n}(2) \end{pmatrix} \in\mathcal{E}_{N}, \quad 0\leq n<T_{\mathbf{a}}.
\end{align*}
Let $U_{0}$ be uniformly distributed on $(\mathfrak{u},1)$, and let $(U_{k})_{k\geq1}$ be independent and identically distributed uniform random variables on $(0,1)$, where
\begin{align} \label{equation:small_i}
    \mathfrak{u} := \frac{\nu}{t+\nu}.
\end{align}
Assume that $U_{0}$, the sequence $(U_{k})_{k\geq1}$, and the walk $(\mathbf{X}_{n})_{n\geq0}$ are independent. For $0\leq j\leq T_{\mathbf{a}}$ and $\mathrm{x} \in [N]$, define
\begin{align*}
    \mathbf{N}_{\mathrm{x}}(j) := \max\left( \{-1\}\cup \left\{ 0\leq k<j: \mathrm{x}\in\{\mathbf{X}_{k}(1),\mathbf{X}_{k}(2)\} \right\} \right).
\end{align*}
Thus, $\mathbf{N}_{\mathrm{x}}(j)$ is the most recent index $k<j$ for which the collision label $\mathbf{X}_{k}$ contains particle $\mathrm{x}$, with value $-1$ if no such index exists. For convenience, set
\begin{align*}
    \mathbf{n}_{j}(q) := \mathbf{N}_{\mathbf{X}_{j}(q)}(j), \quad 0\leq j< T_{\mathbf{a}}, \quad q\in\{1,2\}.
\end{align*}

Conditioned on $T_{\mathbf{a}}=m$, the sequence $(\mathbf{X}_{0},\ldots,\mathbf{X}_{m-1})$ is uniformly distributed on $\mathcal{E}_{N}^{\otimes' m}$. Note that the random walk is indexed from $0$, whereas the collision sequence $\mathbf i_{m}$ is indexed from $1$. More precisely, 
\begin{align} \label{index_adjustment}
    \mathbf{X}_{j-1} \quad \text{corresponds to} \quad (\mathfrak{i}_{j},\mathfrak{i}'_{j}) \quad \forall 1\leq j\leq m.
\end{align}
For each $m\geq1$, denote the conditional joint law of $\mathbf X$ and $(U_{n})_{n\geq0}$ by
\begin{align*}
    \mathds{P}(\cdot) := \mathbb P^{\mathbf X,U}(\cdot| T_{\mathbf{a}}=m),
\end{align*}
and write $\mathds{E}$ for the corresponding expectation. We suppress the dependence on $m$; in sums over $m$, the expectation in each summand is taken under the corresponding conditional law. With this notation, we obtain the following upper bound; see Section \ref{section:SRW_formulation_proof_ideas} for the proof. See also Section \ref{section:informal_explanation} for a heuristic derivation of the optimal bound from \eqref{equation:SRW_formulation_proof_ideas}.

\begin{lemma} \label{lemma:SRW_formulation_proof_ideas}
Assume that $\theta \in \mathbb{R}$ and $N\geq 3$. For each $t,\nu > 0$, $\delta\geq 0$, and $C_{\star} > \max(\theta,0)$, it holds that
\begin{align} \label{equation:SRW_formulation_proof_ideas}
    \mathbb{E}[\mathscr{Z}^{\theta}_{t}(G_{\nu})^{N}] &\leq 
    1+\exp(t \exp(C_{\star}N^{1+\delta})) \frac{t}{\nu}\nonumber \\
    &\sum_{m=1}^{\infty} \bkt{\frac{C_{\eqref{equation:SRW_formulation_proof_ideas}} N^{1-\delta}}{C_{\star} }}^{m} \mathds{E}\bigbkt{ \prod_{j = 1}^{m-1} \frac{1}{U_{j}+\ldots+U_{\mathbf{n}_{j}(1)\wedge \mathbf{n}_{j}(2)+1}}},
\end{align}
where $C_{\eqref{equation:SRW_formulation_proof_ideas}}>0$ is an absolute constant. For $m = 1$, the empty product is understood to equal $1$.
\end{lemma}

\subsection{Cutoff for the path-integral expansion}
We now bound the right-hand side of \eqref{equation:SRW_formulation_proof_ideas} by separating the contributions of short and long path integrals. Set
\begin{align} \label{equation:definition_of_m}
    \mathfrak{m} := \lfloor \exp(N) \rfloor
\end{align}
and write
\begin{align} \label{equation:cut_off}
    \text{(R.H.S.) of \eqref{equation:SRW_formulation_proof_ideas}}
    = \sum_{\textup{short}}+\sum_{\textup{long}}.
\end{align}
Here, the short- and long-path contributions are given as follows, respectively:
\begin{align}
    \sum_{\textup{short}} &:= 1 +  \exp(t \exp(C_{\star}N^{1+\delta})) \frac{t}{\nu} \sum_{m=1}^{\mathfrak{m}-1} \bkt{\frac{C_{\eqref{equation:SRW_formulation_proof_ideas}} N^{1-\delta}}{C_{\star} }}^{m} \mathbf{F}_{N}(m), \label{definition:short_path_integral_contribution}\\
    \sum_{\textup{long}}&:= \exp(t \exp(C_{\star}N^{1+\delta})) \frac{t}{\nu}\sum_{m=\mathfrak{m}}^{\infty} \bkt{\frac{C_{\eqref{equation:SRW_formulation_proof_ideas}} N^{1-\delta}}{C_{\star} }}^{m} \mathbf{F}_{N}(m),\label{definition:long_path_integral_contribution}
\end{align}
where the path functional $\mathbf{F}_{N}(m)$ is given by
\begin{align} \label{equation:F_m_functional}
    \mathbf{F}_{N}(m) &:= \mathds{E}\bigbkt{ \prod_{j = 1}^{m-1} \frac{1}{U_{j}+\ldots+U_{\mathbf{n}_{j}(1)\wedge \mathbf{n}_{j}(2)+1}}}.
\end{align}

The cutoff $\mathfrak{m}$ is chosen so that \cite[Lemma 5.4]{caravenna2019moments} bounds the finite sum in \eqref{definition:short_path_integral_contribution} by $\exp(\exp(CN))$. This yields the following estimate for the short-path contribution; see Section \ref{section:proof_short_path} for the proof.

\begin{lemma} \label{lemma:short_path_integral}
Assume that $\theta \in \mathbb{R}$ and $N \geq 3$. For each $t\geq 1$, $\nu,\delta>0$,  and $C_{\star} \geq \max( \theta,C_{0})$, where $C_{0}>0$ is a sufficiently large absolute constant, there is a constant $\widetilde{C}_{\star} =\widetilde{C}_{\star}(\theta,\nu,C_{\star})>C_{\star}$,  such that
\begin{align} \label{equation:short_path_integral}
    \sum_{\textup{short}}
    \leq 1+\exp(t \exp(\widetilde{C}_{\star}N^{1+\delta})).
\end{align}
\end{lemma}

The cutoff also isolates the range $m\geq\mathfrak{m}$, where a sharper estimate for $\mathbf{F}_{N}(m)$ is available. This estimate contains a factor that decays in $N$ and is raised to the power $m$, which allows us to control the long-path contribution; see Lemma \ref{lemma:estimate_for_path_functional}.

\begin{figure}
\centering
\begin{tikzpicture} [
        scale = 1,
        ->,    > = stealth',
       shorten > = 1pt,auto,
   node distance = 3cm,
      decoration = {snake,   
                    pre length=3pt,post length=7pt,
                    }
                    ]

\path[black] (-2,-0.12) node[anchor=north]{$-1$};
\node at (-2.2,0) [anchor=east]{\small \text{time origin}};
\path[black] (0,-0.12) node[anchor=north]{$0$};
\path[black] (0,0.12) node[anchor=south]{$\mathbf{X}_{0}$};
\path[black] (-1,0) node[anchor=north]{$U_{0}$};
\path[black] (2,-0.12) node[anchor=north]{$1$};
\path[black] (2,0.12) node[anchor=south]{$\mathbf{X}_{1}$};
\path[black] (1,0) node[anchor=north]{$U_{1}$};
\path[black] (4,-0.12) node[anchor=north]{$2$};
\path[black] (4,0.12) node[anchor=south]{$\mathbf{X}_{2}$};
\path[black] (3,0) node[anchor=north]{$U_{2}$};


\filldraw[black] (-2,0) circle (2pt);
\filldraw[black] (0,0) circle (2pt);
\filldraw[black] (2,0) circle (2pt);
\filldraw[black] (4,0) circle (2pt);


\draw[->] (-2,0) -- (7,0) node[anchor=west]{$\mathbb{R}_{+}$};

\draw[<->,dashed,thick] (-2,1) -- (0,1) node[anchor=west]{$\mathbf{Z}_{1}(0) = \mathbf{Z}_{2}(0) = \mathbf{Z}_{3}(0)$};

\node at (0.2,1.2) [anchor=south]{$\Big\Uparrow \mathbf{X}_{0}$ };

\draw[<->,dashed,thick] (0,2.2) -- (2,2.2) node[anchor=west]{$\mathbf{Z}_{1}(1) = \mathbf{Z}_{2}(1)$};
\draw[<->,dashed,thick] (-2,2.6) -- (2,2.6) node[anchor=west]{$\mathbf{Z}_{3}(1)$};

\node at (2.2,2.8) [anchor=south]{$\Big\Uparrow \mathbf{X}_{1}$ };

\draw[<->,dashed,thick] (2,3.8) -- (4,3.8) node[anchor=west]{$\mathbf{Z}_{2}(2) = \mathbf{Z}_{3}(2)$};
\draw[<->,dashed,thick] (0,4.2) -- (4,4.2) node[anchor=west]{$\mathbf{Z}_{1}(2)$};

\end{tikzpicture}
\begin{minipage}{0.9\textwidth}
\caption{
\justifying \footnotesize
This figure illustrates the last collision times in the example from
Figure~1. The label $-1$ marks the time origin, and each label
$j\geq0$ marks the $(j+1)$-st collision. The values of
$Z_{\cdot}(j)$, for $j=0,1,2$, are indicated by dashed lines.
The upward arrows indicate the updates at successive collisions.
}  \label{Figure:last_collision_time}
\end{minipage}
\end{figure}
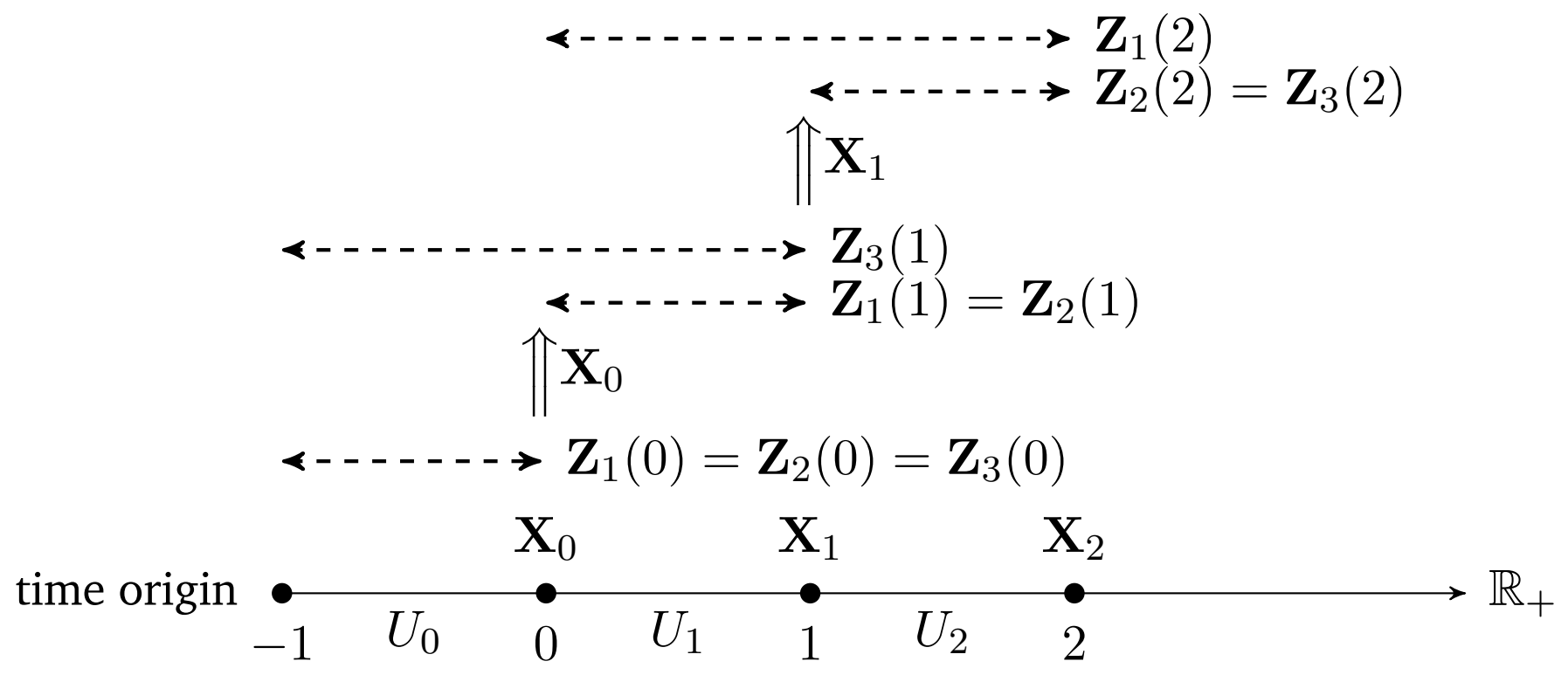

\subsection{The collision-time Markov chain} \label{section:collision_time_markov_chain} To bound the long-path contribution \eqref{definition:long_path_integral_contribution}, we estimate the path functional $\mathbf{F}_{N}(m)$ using an auxiliary Markov chain, defined in \eqref{equation:coulping_process}. Throughout this subsection, $\mathrm{x} \in [N]$ denotes a particle label, while $j$ denotes the collision index, starting from $0$; see \eqref{index_adjustment} for the index shift. Conditioned on $T_{\mathbf a}=m$, the two-body collisions indexed by $j=0,\ldots,m-1$ occur at times
\begin{align*}
    U_{0},\quad U_{0}+U_{1},\quad \ldots,\quad U_{0}+\cdots+U_{m-1}.
\end{align*}
respectively. In what follows, we focus on the new collision-time diagram shown in Figure \ref{Figure:last_collision_time}, rather than the earlier Feynman diagram.

\subsubsection{\textbf{Last collision times}}
Given $(U_{j})_{j \geq 0}$ and $(\mathbf{X}_{j})_{j\geq 0}$, for each particle $\mathrm{x} \in [N]$, define the last collision time
$\mathbf{Z}_{\mathrm{x}}(j)$ by the time between the last time particle $\mathrm{x}$ has a collision and the time of the $j$-th collision, namely, 
\begin{align} \label{definition:last_collition_time}
    \mathbf{Z}_{\mathrm{x}}(j) := \sum_{k = \mathbf{N}_{\mathrm{x}}(j)+1}^{j} U_{k}, \quad 0\leq j\leq T_{\mathbf{a}}, \quad \mathrm{x} \in [N].
\end{align}
See Figure \ref{Figure:last_collision_time} for an example. With these notations, we have the expression
\begin{align} \label{equation:U_and_Z}
    \frac{1}{U_{j}+\ldots+U_{\mathbf{n}_{j}(1)\wedge \mathbf{n}_{j}(2)+1}}
    = \frac{1}{\displaystyle \mathbf{Z}_{\mathbf{X}_{j}(1)}(j) \vee \mathbf{Z}_{\mathbf{X}_{j}(2)}(j)}.
\end{align}
In particular, the path functional \eqref{equation:F_m_functional} can be expressed as
\begin{align} \label{equation:functional_new_expression}
    \mathbf{F}_{N}(m)= \mathds{E}\bigbkt{\prod_{j = 1}^{m-1} \frac{1}{\displaystyle \mathbf{Z}_{\mathbf{X}_{j}(1)}(j) \vee \mathbf{Z}_{\mathbf{X}_{j}(2)}(j)}     }.
\end{align}
\subsubsection{\textbf{Evolution and transition operators}}
We estimate the path functional $\mathbf{F}_{N}(m)$ by iterated conditional expectations. The key is to exploit the underlying Markov chain defined in \eqref{equation:coulping_process}. To describe this chain, we first establish the following evolution equation for $\mathbf{Z}_{\cdot}(j)$.

\begin{lemma} \label{lemma:evolution_equation}
For each $m\geq 1$ and $1\leq j\leq m-1$, it holds that
\begin{align} \label{equation:Evolution_equation}
    \mathbf{Z}_{\cdot}(j+1) = \left( \mathbb{T}_{\mathbf{X}_{j},\mathbf{U}_{j+1}}\mathbf{Z}(j) \right)_{\cdot} \quad \mathds{P}(\cdot  )-\text{(a.s.)}.
\end{align}
\end{lemma}

For each $(X,u) \in \mathcal{E}_{N} \times \mathbb{R}_{+}$, the evolution map $\mathbb{T}_{X,u}: \mathcal{C}_{N} \mapsto \mathcal{C}_{N}$ is given by
\begin{align*}
    (\mathbb{T}_{X,u}Z)_{\mathrm{x}} := 
    \begin{cases}
        \displaystyle u, \quad &\text{if } \mathrm{x} \in \{X(1),X(2)\};\\
        \displaystyle Z_{\mathrm{x}}+u, \quad&\text{if } \mathrm{x} \not\in \{X(1),X(2)\}.
    \end{cases}
\end{align*}
Here $\mathcal{C}_{N} := (0,\infty)^{N}$ denotes the space of the last collision time.
\begin{proof} 
The equation \eqref{equation:Evolution_equation} directly follows from \eqref{definition:last_collition_time}.
\end{proof}

Due to \eqref{equation:Evolution_equation} and \eqref{equation:functional_new_expression}, we study the coupling process 
\begin{align} \label{equation:coulping_process}
    \mathscr{M}:= ((\mathbf{Z}_{\cdot}(j),\mathbf{X}_{j-1}))_{1\leq j\leq m}
\end{align}
rather than the last collision time $(\mathbf{Z}(j))_{j\geq 1}$ alone. We call $\mathscr{M}$ the
collision-time Markov chain. This process takes values almost surely in
\begin{align*}
    \mathcal{C}_{N} \otimes' \mathcal{E}_{N} := \set{(Z,X) \in \mathcal{C}_{N} \times \mathcal{E}_{N}: Z_{X(1)} = Z_{X(2)} < Z_{\mathrm{x}} \quad \forall \mathrm{x} \in [N]\setminus \{X(1),X(2)\} }.
\end{align*}
Indeed, the two particles involved in the previous collision have equal elapsed times, strictly smaller than those of all other particles. The following invariance property will play a key role in estimating the path functional.
\begin{lemma} \label{lemma:function_value_markov_chain}
For each $u \in (0,1)$, $Z \in \mathcal{C}_{N}$, and $X\in\mathcal{E}_{N}$, it holds that
\begin{align*}
    (\mathbb{T}_{X,u}Z,X) \in \mathcal{C}_{N} \otimes' \mathcal{E}_{N}.
\end{align*}
\end{lemma}
\begin{proof}
This property follows from the definition of the evolution map. 
\end{proof}

Let us now verify that $\mathscr{M}$ is a Markov chain, a property that is central to computing the path functional $\mathbf{F}_{N}(m)$.
For convenience, a function defined on $\mathcal{C}_{N} \otimes' \mathcal{E}_{N}$ is referred to as a terminal function.
\begin{lemma} \label{lemma:markov_process}
Assume that $m \geq 3$. For any  $1\leq j\leq m-1$  and every non-negative terminal function $f$, we have
\begin{align*}
    \mathds{E}[f(\mathbf{Z}_{\cdot}(j+1),\mathbf{X}_{j})  |  \mathscr{F}_{j}] = \mathscr{K}_{N}f(\mathbf{Z}_{\cdot}(j),\mathbf{X}_{j-1}) \quad \mathds{P}(\cdot )-\text{(a.s.)},
\end{align*}
where the filtration $\mathscr{F}_{j}$ is given by
\begin{align*}
    \mathscr{F}_{j} := \sigma\sbkt{U_{0},\ldots,U_{j},\mathbf{X}_{0},\ldots,\mathbf{X}_{j-1}}.
\end{align*}
Here the transition kernel $\mathscr{K}_{N}$, defined on the space of terminal functions, is given by 
\begin{align*}
    \mathscr{K}_{N}f(Z,X) := \frac{1}{M_{N}-1} \sum_{X'\in \mathcal{E}_{N} \setminus \{X\}  } \int_{0}^{1} f(\mathbb{T}_{X',u}Z,X') du, \quad (Z,X) \in \mathcal{C}_{N} \otimes' \mathcal{E}_{N},
\end{align*}
where we have set
\begin{align} \label{definition:M_N}
    M_{N} := \binom{N}{2}.
\end{align}
In particular,  $((\mathbf{Z}_{\cdot}(j),\mathbf{X}_{j-1}))_{1\leq j\leq m}$ is a $\mathcal{C}_{N} \otimes' \mathcal{E}_{N}$-valued Markov chain with respect to $(\mathscr{F}_{j})_{1\leq j\leq m}$. 
\end{lemma}
\begin{proof}
Let $1\leq j\leq m-1$. Applying \eqref{equation:Evolution_equation} gives
\begin{align*}
    \mathds{E}\bigbkt{f(\mathbf{Z}_{\cdot}(j+1),\mathbf{X}_{j}) \biggr |  \mathscr{F}_{j}}
    &= \mathds{E}\bigbkt{ f(\mathbb{T}_{\mathbf{X}_{j},U_{j+1}}\mathbf{Z}(j),\mathbf{X}_{j}) \biggr |\mathscr{F}_{j}}\\
    &= \frac{1}{M_{N}-1} \sum_{X' \in \mathcal{E}_{N}\setminus \{\mathbf{X}_{j-1}\} } \int_{0}^{1} du  f(\mathbb{T}_{X',u}\mathbf{Z}(j),X') = \mathscr{K}_{N}f(\mathbf{Z}(j),\mathbf{X}_{j-1}).
\end{align*}
Combining this with Lemma \ref{lemma:function_value_markov_chain}, we know that $\mathscr{M}$ is a $\mathcal{C}_{N} \otimes' \mathcal{E}_{N}$-valued Markov chain with respect to $(\mathscr{F}_{j})_{1\leq j\leq m}$. 
\end{proof}

\subsubsection{\textbf{Iterates of the weighted transition operator}}

We now generalize the path functional \eqref{equation:functional_new_expression} by replacing the constant terminal function $\mathds{1}$ with a general nonnegative terminal function. We express this generalized functional in terms of iterates of a weighted transition operator associated with the collision-time Markov chain $\mathscr{M}$.

\begin{proposition}
\label{proposition:FK_representation}
For every $m\geq2$ and every nonnegative Borel-measurable terminal function $f$,
\begin{align} \label{equation:FK_representation}
    \mathds{E}\bigbkt{\prod_{j = 1}^{m-1} \frac{1}{U_{j}+\ldots+U_{\mathbf{n}_{j}(1)\wedge \mathbf{n}_{j}(2)+1}}  f\sbkt{\mathbf{Z}_{\cdot}(m),\mathbf{X}_{m-1}}  \biggr | (\mathbf{Z}_{\cdot}(1),\mathbf{X}_{0}) = (Z,X) } = \mathscr{Q}_{N}^{m-1}f(Z,X).
\end{align}
Here the weighted transition operator $\mathscr{Q}_{N}$ acts on nonnegative terminal functions by
\begin{align} \label{equation:weight_transition_kernel}
    \mathscr{Q}_{N}f(Z,X) := \frac{1}{M_{N}-1} \sum_{X'\in \mathcal{E}_{N}\setminus \{X\} } \frac{1}{Z_{X'(1)} \vee Z_{X'(2)}} \int_{0}^{1} f(\mathbb{T}_{X',u}Z,X') du, \quad (Z,X) \in \mathcal{C}_{N} \otimes' \mathcal{E}_{N}.
\end{align}
\end{proposition}
\begin{proof}
We prove \eqref{equation:FK_representation} by induction on $m$. Rewrite the left-hand side of \eqref{equation:FK_representation} using \eqref{equation:U_and_Z}. Applying the Markov property stated in Lemma \ref{lemma:markov_process} yields
\begin{align*}
    &\mathds{E}\bigbkt{\prod_{j = 1}^{m-1} \frac{1}{\displaystyle \mathbf{Z}_{\mathbf{X}_{j}(1)}(j) \vee \mathbf{Z}_{\mathbf{X}_{j}(2)}(j)}   f\sbkt{\mathbf{Z}_{\cdot}(m),\mathbf{X}_{m-1}}  \biggr | (\mathbf{Z}_{\cdot}(1),\mathbf{X}_{0}) = (Z,X)  }\\
    &= \mathds{E}\bigbkt{ \frac{1}{\displaystyle \mathbf{Z}_{\mathbf{X}_{1}(1)}(1) \vee \mathbf{Z}_{\mathbf{X}_{1}(2)}(1)}  \\
    &\quad\quad\; \mathds{E}\bigbkt{\prod_{j = 2}^{m-1} \frac{1}{\displaystyle \mathbf{Z}_{\mathbf{X}_{j}(1)}(j) \vee \mathbf{Z}_{\mathbf{X}_{j}(2)}(j)}   f\sbkt{\mathbf{Z}_{\cdot}(m),\mathbf{X}_{m-1}}\biggr | (\mathbf{Z}_{\cdot}(2),\mathbf{X}_{1})     }  \biggr | (\mathbf{Z}_{\cdot}(1),\mathbf{X}_{0}) = (Z,X)  }\\
    &= \mathds{E}\bigbkt{ \frac{1}{\displaystyle \mathbf{Z}_{\mathbf{X}_{1}(1)}(1) \vee \mathbf{Z}_{\mathbf{X}_{1}(2)}(1)}  \mathscr{Q}_{N}^{m-2}f(\mathbf{Z}_{\cdot}(2),\mathbf{X}_{1})  \biggr | (\mathbf{Z}_{\cdot}(1),\mathbf{X}_{0}) = (Z,X)  } = \mathscr{Q}_{N}^{m-1}f(Z,X).
\end{align*}
The proof is complete. 
\end{proof}

\subsection{Contraction for the weighted transition operator} \label{section:lya}
The main advantage of \eqref{equation:FK_representation} is that it allows us to estimate the path functional $\mathbf{F}_{N}(m)$ through one-step bounds for $\mathscr{Q}_{N}$. Specifically, we construct a positive weight function $\mathcal{H}_{a}$ satisfying a Lyapunov-type inequality, with an explicit contraction rate depending on $N$ and $a$; see Proposition \ref{proposition:main_contractiv_inequality}. We then iterate this inequality using positivity of $\mathscr{Q}_{N}$ and apply the resulting bounds in \eqref{equation:main_contractiv_inequality}; see Section \ref{section:estimate_for_path_functional} for details. This approach avoids estimating the full product in \eqref{equation:F_m_functional} directly, which is difficult because its denominators share many of the same random variables. The proof of the following proposition is given in Section \ref{section:proof_contraction_inequality}.

\begin{proposition} \label{proposition:main_contractiv_inequality}
Suppose that $N\geq 3$ and $0<a<1/3$. Define the multiplicative and additive weight functions on $\mathcal{C}_{N}$ by
\begin{align*}
    \mathcal{M}_{a}(Z) := \prod_{x \in [N]} Z_{x}^{-a},  \quad \mathcal{A}_{a}(Z) := \sum_{x \in [N]} Z_{x}^{-\mathfrak{q}_{a}}, \quad Z \in \mathcal{C}_{N}, \quad \mathfrak{q}_{a} := 1-\frac{5}{2} a, 
\end{align*}
and set contractive weight function as
\begin{align*}
    \mathcal{H}_{a}(Z) := \mathcal{M}_{a}(Z) (1+\mathcal{A}_{a}(Z)).
\end{align*}
Then it holds that
\begin{align} \label{equation:main_contractiv_inequality}
    \mathscr{Q}_{N} \mathcal{H}_{a}(Z,X) \leq C_{\eqref{equation:main_contractiv_inequality}}\bkt{ \frac{1}{a^{2} N^{2}}+ \frac{1}{a N^{1-2a}}} \mathcal{H}_{a}(Z), \quad (Z,X) \in \mathcal{C}_{N} \otimes' \mathcal{E}_{N},
\end{align}
where $C_{\eqref{equation:main_contractiv_inequality}} >0$ is an absolute constant. 
\end{proposition}

The construction of $\mathcal{H}_{a}$ is guided by the following idea. We seek a finite collection $\mathbb{F} := \{\mathcal{F}_{\alpha}\}_{\alpha\in I}$ of nonnegative terminal functions, with $|I|$ independent of $N$, such that
\begin{align} \label{equation:Lya}
    \mathscr{Q}_{N}\mathcal{F}_{\alpha}(Z,X)
\leq
\sum_{\alpha'\in I}
C_{\alpha \rightarrow \alpha'}(N,a)\mathcal{F}_{\alpha'}(Z,X),
\quad \alpha\in I, \quad (Z,X) \in \mathcal{C}_{N} \otimes' \mathcal{E}_{N}.
\end{align}
In particular, the transition coefficients are non-negative and satisfy
\begin{align*}
    \max_{\alpha,\alpha'\in I}C_{\alpha \rightarrow \alpha'}(N,a)=o(1) \quad\text{as }N\to\infty.
\end{align*}
By linearity and nonnegativity, summing the inequalities \eqref{equation:Lya} gives
\begin{align*}
    \mathscr{Q}_{N}\left(\sum_{\alpha\in I}\mathcal{F}_{\alpha}\right)
\leq
|I|\max_{\alpha,\alpha'\in I}C_{\alpha \rightarrow \alpha'}(N,a)
\left(\sum_{\alpha\in I}\mathcal{F}_{\alpha}\right).
\end{align*}
Thus, the sum satisfies a contraction inequality for sufficiently large $N$. In our construction, $\mathbb{F} = \{\mathcal{M}_{a},\mathcal{M}_{a}\mathcal{A}_{a}\}$, whose sum is $\mathcal{H}_{a}$.

\subsection{Long path integrals and completion of the proof} \label{section:main_result_proof_idea}
The key to bounding the long-path integral \eqref{definition:long_path_integral_contribution} is to apply the Lyapunov-type inequality \eqref{equation:main_contractiv_inequality} to estimate the path functional $\mathbf{F}_{N}(m)$. To do this, we change the terminal function $\mathds{1}$ in the path functional $\mathbf{F}_{N}(m)$ to $\mathcal{H}_{a}$. The proof of Lemma \ref{lemma:estimate_for_path_functional} is given in Section \ref{section:estimate_for_path_functional}.
\begin{lemma} \label{lemma:estimate_for_path_functional}
Assume that $N\geq 3$.  For each $t\geq 1$, $\nu >0$, and $0<a<1/3$, we have
\begin{align} \label{equation:estimate_for_path_functional}
    \mathbf{F}_{N}(m) \leq  4\mathfrak{u}^{-aN} N^{2} (2C_{\eqref{equation:main_contractiv_inequality}})^{m} \bkt{\frac{1}{a^{2} N^{2}}+\frac{1}{aN^{1-2a}}}^{m} \quad \forall m\geq \mathfrak{m}.
\end{align}
\end{lemma}

Now we are ready to bound the long path integral. Recall that
\begin{align*}
    \sum_{\textup{long}} = \exp(t \exp(C_{\star}N^{1+\delta})) \frac{t}{\nu}\sum_{m=\mathfrak{m}}^{\infty} \bkt{\frac{C_{\eqref{equation:SRW_formulation_proof_ideas}} N^{1-\delta}}{C_{\star} }}^{m} \mathbf{F}_{N}(m).
\end{align*}
\begin{lemma} \label{lemma:long_path_integral}
Assume that $\theta \in \mathbb{R}$ and $N\geq 5$.  For each $t\geq 1$, $2/N \leq \delta <2/3$, and $\nu>0$, there exists $C_{\eqref{equation:long_path_integral}} = C_{\eqref{equation:long_path_integral}}(\theta,\nu) > 0$ such that 
\begin{align} \label{equation:long_path_integral}
    \sum_{\textup{long}}
    \leq \exp\bkt{t\exp\bkt{\frac{C_{\eqref{equation:long_path_integral}}}{\delta}N^{1+\delta}} }. 
\end{align}
\end{lemma}
\begin{proof} 
Now by choosing  
\begin{align*}
    a := \frac{\delta}{2}, \quad C_{\star} := \frac{4((2C_{\eqref{equation:SRW_formulation_proof_ideas}}C_{\eqref{equation:main_contractiv_inequality}})\vee \theta)}{a}, 
\end{align*}
applying Lemma \ref{lemma:estimate_for_path_functional} gives
\begin{align*}
    \sum_{\textup{long}} 
    &\leq \exp(t \exp(C_{\star}N^{1+\delta})) \frac{4t}{\nu} N^{2} \mathfrak{u}^{-a N} \sum_{m = \mathfrak{m}}^{\infty} \bkt{\bkt{\frac{2C_{\eqref{equation:SRW_formulation_proof_ideas}} C_{\eqref{equation:main_contractiv_inequality}}}{a C_{\star}} } \cdot a\bkt{\frac{1}{a^{2} N^{1+\delta}}+\frac{1}{aN^{\delta-2a}}}}^{m}\\
    &\leq \exp(t \exp(C_{\star}N^{1+\delta})) \frac{4t}{\nu} N^{2} \mathfrak{u}^{-a N} \sum_{m = \mathfrak{m}}^{\infty} \bkt{\frac{1}{2}}^{m},
\end{align*}
where we have used $a N^{1+\delta} \geq aN \geq 1$ due to $\delta \geq 2/N$. In addition, since
\begin{align*}
    \mathfrak{u}^{-aN} = \exp\bkt{aN \log\bkt{\frac{t+\nu}{\nu}}} \leq \exp\bkt{\frac{t}{\nu}\exp(N)} \leq \exp\bkt{t\exp\bkt{ (\nu^{-1} + 1) N }},
\end{align*}
choosing a large positive constant $C$, depending on both $\theta$ and $\nu$, gives
\begin{align*}
    \sum_{\textup{long}} \leq \exp\bkt{t \exp\bkt{\frac{C}{\delta} N^{1+\delta} } }.
\end{align*}
The proof is complete.  
\end{proof}

Theorem \ref{theorem:main_result_proof_idea} now follows by combining Lemma \ref{lemma:short_path_integral} and Lemma \ref{lemma:long_path_integral}. 

\subsection{Remarks on the proof} \label{section:proof_remarks}
In Sections~\ref{section:near_optimal_tree}-\ref{section:informal_explanation},
we explain heuristically why $T_{\star}$ should suffice to obtain the
optimal bound $\exp(t\exp(CN))$.
In Section~\ref{section:trade_off}, we illustrate the difficulties
that arise when applying the AM-GM inequality to establish
upper moment bounds.

\subsubsection{\textbf{Heuristic justification of the spanning-tree choice}} \label{section:near_optimal_tree}  Let $T_{\mathrm{max}}$ maximize $W(T)$ over $\mathcal{T}(G)$. Since $|\mathcal{T}(G)|\leq2^{|E|}=2^{3m}$, we have
\begin{align*}
    2^{-3m}W(T_{\mathrm{max}})^{-1} \leq \left(\sum_{T\in\mathcal{T}(G)}W(T)\right)^{-1} \leq W(T_{\mathrm{max}})^{-1}.
\end{align*}
Thus, retaining only the maximum-weight tree loses at most a factor $2^{3m}$. This loss merely changes the constant $C$ in \eqref{equation:representation_path_integral} and can be absorbed by increasing $C_{\star}$ in \eqref{equation:SRW_formulation_proof_ideas}, without affecting the target growth scale in $N$.

However, $T_{\mathrm{max}}$ depends on the temporal variables $a_{j},b_{j}$ as well as the collision sequence $\mathbf{i}_{m}$, making the resulting bound difficult to analyze. We therefore use $T_{\star}$, whose edge selection depends only on $\mathbf{i}_{m}$ and favors short edges with large weights. The heuristic comparison $W(T_{\star})\approx W(T_{\mathrm{max}})$ is meant at the level of the final moment bound: using $T_{\star}$ in \eqref{equation:QS_spanning_tree} is still expected to yield $\exp(t\exp(CN))$, as suggested by the following argument.

\subsubsection{\textbf{Heuristic for the optimal moment bound}} \label{section:informal_explanation}

To obtain the optimal bound $\exp(t\exp(CN))$ from \eqref{equation:SRW_formulation_proof_ideas}, it suffices to set
$\delta=0$ and bound the series over $m$ by $\exp(t\exp(C'N))$. The heuristic below suggests that this series should, in fact, be
bounded by $N^{O(N^{1+\varepsilon})}$ for any fixed $\varepsilon>0$.

To explain this, notice that the number of uniformly distributed random variables in a reciprocal sum, $j-\mathbf{n}_{j}(1)\wedge \mathbf{n}_{j}(2)$, has a geometric-type distribution. Once $j=j_{N}\geq N^{1+\varepsilon}$ and $m > j+\ell$ for some nonnegative integer $\ell$, this number is of order $N$:
\begin{align*}
        \bkt{\frac{j+k - \mathbf{n}_{j+k}(1)\wedge \mathbf{n}_{j+k}(2)}{N}}_{0\leq k\leq    \ell} \overset{d}{\longrightarrow} (e_{2k+1} \vee e_{2k+2})_{0\leq k\leq \ell},
\end{align*}
where $(e_{n})_{n\geq1}$ are independent exponential random variables with rate $2$. Combining this with the law of large numbers for the uniform increments gives 
\begin{align*}
        \bkt{\frac{U_{j+k}+\ldots+U_{\mathbf{n}_{j+k}(1)\wedge \mathbf{n}_{j+k}(2)+1}}{N}}_{0\leq k\leq \ell} \overset{d}{\longrightarrow} (\mathbb{E}[U_{1}] (e_{2k+1} \vee e_{2k+2}))_{0\leq k\leq \ell}.
\end{align*}
This suggests that each reciprocal sum contributes a factor of order $1/N$ after the first $N^{1+\varepsilon}$ collisions. If this decay could be accumulated, then choosing $C_{\star}$ sufficiently large would make the remaining series geometric. A coarser estimate for the initial part would then yield a bound of the form $N^{O(N^{1+\varepsilon})}$. However, making this heuristic rigorous requires balancing decay in $N$ against integrability.

\subsubsection{\textbf{Trade-off between decay in $N$ and integrability}} \label{section:trade_off}

Making the heuristic in Section \ref{section:informal_explanation}
rigorous requires quantitative control of the expectation in \eqref{equation:SRW_formulation_proof_ideas}. The heuristic suggests that each factor 
\begin{align*}
    I_{j}:= \frac{1}{U_{j}+\ldots+U_{\mathbf{n}_{j}(1)\wedge \mathbf{n}_{j}(2)+1}}.
\end{align*}
should contribute a decay of order $1/N$ for sufficiently large $j$. Turning this heuristic into an upper bound, however, requires a balance between preserving this decay and controlling the singularities of the full product. To explain this difficulty, set $\ell_{j}:= j-\mathbf{n}_{j}(1)\wedge \mathbf{n}_{j}(2)$. Since $\ell_{j}$ is typically of order $N$, retaining the dependence on $\ell_{j}$ is essential for recovering the decay $1/N$. A natural approach is to apply the AM-GM inequality to $I_{j}$:
\begin{align} \label{equation:am_gm_inequality}
    \prod_{j = 1}^{m-1} I_{j} \leq \prod_{j=1}^{m-1} \bkt{\frac{1}{\lfloor \ell_{j}/2 \rfloor} \prod_{k = 0}^{\lfloor\ell_{j}/2\rfloor-1} \frac{1}{(U_{j-2k}+U_{j-(2k+1)})^{1/\lfloor\ell_{j}/2\rfloor}}}.
\end{align}
The difficulty with \eqref{equation:am_gm_inequality} is that the denominators associated with different $I_{j}$ overlap. Consequently, $(U_{j}+U_{j-1})^{-1}$ may occur in several AM–GM upper bounds, and its accumulated exponent may become too large. 
More precisely, let $\sigma_{j} := \sum_{x \in K_{j}} 1/\lfloor \ell_{x}/2 \rfloor$ and $K_{j}$ be the subset of $[m-1]$ so that for every $x \in K_{j}$,  $U_{j}+U_{j-1}$ appears in $I_{x}$ and $x-j$ is even. For example, consider $\mathbf{i}_{9} = ((1,n))_{2\leq n\leq 10}$. We have $\sigma_{2} = \sum_{i = 1}^{4}1/i >2$. In this case, the expectation on the right-hand side of \eqref{equation:am_gm_inequality} is bounded below by
\begin{align*}
    \mathbb{E}\bigbkt{\frac{1}{(U_{2}+U_{1})^{\sigma_{2}}} \frac{1}{(U_{1}+U_{0})^{\sigma_{1}}} } = \infty.
\end{align*}
Thus, the AM-GM estimate can produce a nonintegrable upper bound. Even when no $\sigma_{j}$ reaches $2$, the overlap $\sigma_{j}>1$ requires additional control when seeking an estimate analogous to \cite[Lemma 5.4]{caravenna2019moments}.

The above illustrates the trade-off. Using many pairs in the AM-GM estimate preserves the small factor $1/\ell_{j}$, but repeated use of the same pairs across different $I_{j}$ may lead to a nonintegrable upper bound. Restricting the choice of pairs can limit this accumulation, but may force us to retain fewer terms in some denominators, thereby weakening the decay obtained from their lengths. 

Thus, the challenge is not merely to obtain an integrable upper bound, but to do so while retaining sufficient decay in $N$. This motivates the Markov-chain formulation introduced in Section \ref{section:collision_time_markov_chain}: rather than applying a pathwise bound to the product of $I_{j}$, we estimate its expectation through one-step inequalities for a
weighted transition operator.

\section{Derivation of the probabilistic upper bound}
\label{section:SRW_formulation_proof_ideas}

The goal of this section is to prove Lemma \ref{lemma:SRW_formulation_proof_ideas}. Throughout the proof, $C$ denotes a positive absolute constant
whose value may change from line to line.

By Lemma \ref{lemma:Gaussian_text_function}, Lemma \ref{lemma:QS_spanning_tree}, and the fact that $T_{\star}$ contains every edge in $E_{S}$, we obtain
\begin{align*}
    &\mathbb{E}[\mathscr{Z}^{\theta}_{t}(G_{\nu})^{N}] \leq 1+\sum_{m=1}^{\infty} C^{m} \sum_{ \mathbf{i}_{m} \in \mathcal{E}_{N}^{\otimes' m}} \\
    &\int_{\nu<a_{1}<b_{1}<\ldots<a_{m}<b_{m} < t+\nu} \prod_{j=1}^{m} da_{j} db_{j} \prod_{j=1}^{m} \begin{bmatrix}
        \displaystyle \frac{1}{(b_{j}-a_{j})} \\
        \displaystyle \mathfrak{s}^{\theta}(b_{j}-a_{j}) \\
    \end{bmatrix}_{\times} 
    \begin{bmatrix}
        \displaystyle \frac{1}{a_{j}-b_{\mathrm{p}^{\mathbf{i}_{m}}(j,\mathfrak{i}'_{j})}}\\
        \displaystyle \frac{1}{a_{j}-b_{\mathrm{p}^{\mathbf{i}_{m}}(j,\mathfrak{i}_{j})}}\\
    \end{bmatrix}_{\times} W(T_{\star})^{-1}\\
    &\leq 1+\sum_{m=1}^{\infty}  C^{m} \sum_{ \mathbf{i}_{m} \in \mathcal{E}_{N}^{\otimes' m}} \int_{\nu<a_{1}<b_{1}<\ldots<a_{m}<b_{m} < t+\nu} \prod_{j=1}^{m} da_{j} db_{j}  \prod_{j=1}^{m} \mathfrak{s}^{\theta}(b_{j}-a_{j}) \prod_{e\in E(T_{\star})^{c}\cap E_{L}} W(e).
\end{align*}
Here the weights of the edges in $T_{\star}$ cancel against the corresponding factors in the integrand, with the remaining factors absorbed into $C^{m}$.

Introduce the temporal variables
\begin{align*}
    u_{j-1} := a_{j} - b_{j-1}, \quad  v_{j} := b_{j} - a_{j} \quad \forall 1\leq j\leq m.
\end{align*}
For $j\geq2$, $u_{j-1}$ is the time between the end of the $(j-1)$-th collision  and the beginning of the $j$-th collision, while $u_{0}=a_{1}$ is the time before the first collision. The variable $v_{j}$ is the duration of the $j$-th collision; see Figure \ref{Figure:Feynman}.

We next separate the $u$- and $v$-integrals so that we can apply the Laplace transform of $\mathfrak s^{\theta}$ to the latter. For an edge $e=(b_{k},a_{j})\in E_{L}$, where $0\leq k<j\leq m$, we have
\begin{align*}
    W(e)^{-1} = a_{j}-b_{k} = \sum_{r=k}^{j-1}u_{r} + \sum_{r=k+1}^{j-1}v_{r} \geq \sum_{r=k}^{j-1}u_{r} =: W'(e)^{-1}.
\end{align*}
Then $W(e)\leq W'(e)$, and $W'(e)$ depends only on the $u$-variables. Replacing $W$ by $W'$ and enlarging the integration domain therefore gives
\begin{align*}
    &\mathbb{E}[\mathscr{Z}^{\theta}_{t}(G_{\nu})^{N}] \leq 1+\sum_{m=1}^{\infty} C^{m}\\
    &\sum_{ \mathbf{i}_{m} \in \mathcal{E}_{N}^{\otimes' m}}
    \int_{0<v_{1}+\ldots+v_{m} < t} \prod_{j=1}^{m} dv_{j} \prod_{j=1}^{m} \mathfrak{s}^{\theta}(v_{j})  \int_{\nu<u_{0}<t+\nu,0<u_{1}<t,\ldots,0<u_{m-1}<t} \prod_{j=0}^{m-1} du_{j} \prod_{e \in E(T_{\star})^{c}\cap E_{L}} W'(e).
\end{align*}

Now we estimate the integrals with respect to $dv_{j}$. To do this, set
\begin{align} \label{definition:lambda_N}
    \lambda_{N,\delta} := \exp(C_{\star} N^{1+\delta})
\end{align}
and apply the following trick to bound the integrals with respect to $dv_{j}$:
\begin{align*}
    &\int_{0<v_{1}+\ldots+v_{m} < t} \prod_{j=1}^{m} dv_{j} \prod_{j=1}^{m} \mathfrak{s}^{\theta}(v_{j})\\
    &\leq \exp(\lambda_{N,\delta} t) \int_{0<v_{1}+\ldots+v_{m} < t} \prod_{j=1}^{m} dv_{j} \prod_{j=1}^{m} \exp(-\lambda_{N,\delta} v_{j})  \mathfrak{s}^{\theta}(v_{j})\\
    &\leq \exp(\lambda_{N,\delta} t) \bkt{\int_{0}^{\infty} dv \exp(-\lambda_{N,\delta} v) \mathfrak{s}^{\theta}(v) }^{m} \leq \exp(\lambda_{N,\delta} t) \bkt{\frac{4\pi}{\log \lambda_{N,\delta} - \theta}}^{m}\\
    &\leq C^{m} \exp(\lambda_{N,\delta} t) \bkt{ \frac{1}{\log \lambda_{N,\delta}}}^{m} = \bkt{\frac{C}{C_{\star} N^{1+\delta}} }^{m} \exp(t \exp(C_{\star}N^{1+\delta})).
\end{align*}
Here we used the Laplace transform of $s^{\theta}(v)$:
\begin{align*}
    \int_{0}^{\infty} dv e^{-\lambda v}\mathfrak s^{\theta}(v)  = \frac{4\pi}{\log\lambda-\theta}, \quad \lambda>e^{\theta},
\end{align*}
which follows from \cite[Lemma 8.4]{Gu_2d} and \eqref{definition:beta_star}. Since $N\geq3$ and $\delta\geq0$,
\begin{align*}
    \log\lambda_{N,\delta} = C_{\star}N^{1+\delta} \geq3C_{\star}\geq2\theta,
\end{align*}
which implies $\log\lambda_{N,\delta}-\theta \geq (\log\lambda_{N,\delta})/2$. Consequently,
\begin{align*}
    \mathbb{E}[\mathscr{Z}^{\theta}_{t}(G_{\nu})^{N}] &\leq 1+\exp(t \exp(C_{\star}N^{1+\delta})) \\
    &\sum_{m=1}^{\infty} \bkt{\frac{C}{C_{\star} N^{1+\delta}}}^{m} \sum_{ \mathbf{i}_{m} \in \mathcal{E}_{N}^{\otimes' m}} \int_{\nu<u_{0}<t+\nu,0<u_{1}<t,\ldots,0<u_{m-1}<t} \prod_{j=0}^{m-1} du_{j} \prod_{e \in E(T_{\star})^{c}\cap E_{L}} W'(e).
\end{align*}

We now express the sum over collision sequences $\mathbf{i}_{m}\in\mathcal{E}_{N}^{\otimes' m}$ as an expectation with respect to $\mathbf{X}$. Conditioned on $T_{\mathbf a}=m$, the sequence $(\mathbf X_{0},\ldots,\mathbf X_{m-1})$ is uniformly distributed over this set. Identifying $\mathbf{X}_{j}=(\mathfrak{i}_{j+1},\mathfrak{i}'_{j+1})$, the choice of $T_{\star}$ gives
\begin{align*}
    \prod_{e\in E_{L}\setminus E(T_{\star})}W'(e) = \frac{1}{u_{0}} \prod_{j=1}^{m-1} \frac{1}{ u_{j}+\cdots+ u_{(\mathbf{n}_{j}(1)\wedge\mathbf{n}_{j}(2))+1} }.
\end{align*}
Next, rescale the $u$-variables by $t+\nu$ and enlarge the integration intervals for $u_{1},\ldots,u_{m-1}$ to $(0,1)$. The Jacobian cancels the scaling of the reciprocal factors. Using $u_{0}^{-1}\leq\mathfrak u^{-1}$, we bound the integral by the corresponding expectation over $U$, multiplied by $t/\nu$. Combining these observations with the preceding estimate yields
\begin{align} \label{equation:tight_upper_bound}
    \mathbb{E}[\mathscr{Z}^{\theta}_{t}(G_{\nu})^{N}] \leq 1+&\exp(t \exp(C_{\star}N^{1+\delta})) \frac{t}{\nu} \sum_{m=1}^{\infty} \bkt{\frac{C N^{1-\delta}}{C_{\star} }}^{m} \mathds{E}\bigbkt{ \prod_{j = 1}^{m-1} \frac{1}{U_{j}+\ldots+U_{\mathbf{n}_{j}(1)\wedge \mathbf{n}_{j}(2)+1}} }.
\end{align}
This proves Lemma \ref{lemma:SRW_formulation_proof_ideas}.

\section{Bounds for short path integrals} \label{section:proof_short_path}

In this section, we prove the short-path estimate in Lemma \ref{lemma:short_path_integral}. The key observation is that consecutive collision labels are distinct, which implies
\begin{align} \label{equation:length_at_least_two}
    j - \mathbf{n}_{j}(1) \wedge \mathbf{n}_{j}(2) \geq 2.
\end{align}
Consequently, $U_{j}+\cdots+ U_{(\mathbf{n}_{j}(1)\wedge\mathbf{n}_{j}(2))+1} \geq U_{j}+U_{j-1}$.
Using the definition \eqref{definition:short_path_integral_contribution}, we obtain
\begin{align*}
    \sum_{\textup{short}} &\leq 1 +  \exp(t \exp(C_{\star}N^{1+\delta})) \frac{t}{\nu} \sum_{m=1}^{\mathfrak{m}-1} \bkt{\frac{C N^{1-\delta}}{C_{\star} }}^{m} \mathds{E}\bigbkt{ \prod_{j = 1}^{m-1} \frac{1}{U_{j}+U_{j-1}}  }\\
    &\leq 1 +  \exp(t \exp(C_{\star}N^{1+\delta})) \frac{t}{\nu} \sum_{m=1}^{\mathfrak{m}-1} \bkt{\frac{C N^{1-\delta}}{C_{\star} }}^{m} \mathds{E}\bigbkt{ \prod_{j = 1}^{m-1} \frac{1}{U_{j}+U_{j-1}} },
\end{align*}
where $C>0$ is an absolute constant. By \cite[Lemma 5.4]{caravenna2019moments}, we know that 
\begin{align} \label{equation:lemma_5_4}
    \mathds{E}\bigbkt{ \prod_{j = 1}^{m-1} \frac{1}{U_{j}+U_{j-1}} }
    \leq \mathds{E}\bigbkt{ \prod_{j = 1}^{m-1} \frac{4}{\sqrt{U_{j-1} (U_{j-1}+U_{j})} } } \leq c^{m},
\end{align}
where $c>0$ is an absolute constant. In addition, since $m\leq \mathfrak{m}-1$, it follows that
\begin{align*}
    N^{m} \leq \exp(m \log N) \leq \exp(\exp(2N)),
\end{align*}
which implies
\begin{align*}
    \sum_{\textup{short}} &\leq 1+\exp(t \exp(C_{\star}N^{1+\delta})) \exp(\exp(2N)) \frac{t}{\nu} \sum_{m=1}^{\mathfrak{m}-1}\bkt{\frac{C \cdot c}{C_{\star} }}^{m} \\
    &\leq 1+\exp(t \exp(C_{\star}N^{1+\delta})) \exp(\exp(2N)) \frac{t}{\nu}\sum_{m=1}^{\infty} \frac{1}{2^{m}}.
\end{align*}
Here we have chosen $C_{\star} \geq 2C \cdot c$. This proves the short-path estimate \eqref{equation:short_path_integral} by choosing $\widetilde{C}_{\star}>C_{\star}$ sufficiently large.

\section{Contraction inequalities for the weighted transition operator} \label{section:proof_contraction_inequality}

In this section, we prove Proposition \ref{proposition:main_contractiv_inequality}, which follows by adding the two inequalities below and adjusting the absolute constant.

\begin{proposition}
Suppose that $N\geq 3$ and $0<a< 1/3$. Then, for each $(Z,X) \in \mathcal{C}_{N} \otimes' \mathcal{E}_{N}$,
\begin{align} 
    \mathscr{Q}_{N} \mathcal{M}_{a}(Z,X) &\leq C\bkt{\frac{1}{a^{2}N^{2}}+\frac{1}{a N^{1-2a}}} \mathcal{H}_{a}(Z), \label{lemma:contractiv_inequality_P}\\
    \mathscr{Q}_{N} (\mathcal{M}_{a} \mathcal{A}_{a}) (Z,X) &\leq \frac{C}{a N^{1-2a}} (\mathcal{M}_{a}\mathcal{A}_{a})(Z), \label{lemma:contractiv_inequality_P_S} 
\end{align}
where $C > 0$ is an absolute constant.
\end{proposition}

Throughout the proof, $C>0$ denotes an absolute constant whose value may change from line to line.

\subsection{The multiplicative weight function} \label{section:multiplicative_weight_function}
We first prove \eqref{lemma:contractiv_inequality_P}. Recall that
\begin{align} \label{equation:section_multiplicative_weight_function_0}
    \mathscr{Q}_{N} \mathcal{M}_{a}(Z,X) = \frac{1}{M_{N}-1} \sum_{X'\in \mathcal{E}_{N} \setminus \{X\} } \frac{1}{Z_{X'(1)} \vee Z_{X'(2)}} \int_{0}^{1} \mathcal{M}_{a}(\mathbb{T}_{X',u}Z,X') du.
\end{align}
Fix $0 < a < 1/3$ and $(Z,X) \in \mathcal{C}_{N} \otimes' \mathcal{E}_{N}$. Order the coordinates of $Z$ as
\begin{align} \label{order_z}
    z_{1} \leq z_{2} \leq \ldots \leq z_{N},
\end{align}
assigning ranks $1$ and $2$ to the particles $X(1)$ and $X(2)$, respectively. This is possible by the definition of $\mathcal{C}_{N}\otimes'\mathcal{E}_{N}$. Thus, we have 
\begin{align*}
    z_{1} =Z_{X(1)}=Z_{X(2)} = z_{2}.
\end{align*}
For each $X'\in\mathcal{E}_{N}\setminus\{X\}$, let $i<k$ be the ranks of its two particles: 
\begin{align} \label{particle_rank}
    z_{i} = \min(Z_{X'(1)},Z_{X'(2)}), \quad z_{k} = \max(Z_{X'(1)},Z_{X'(2)}).
\end{align}
Since $X'\neq X$, we necessarily have
\begin{align} \label{condition:k}
    k\geq 3.
\end{align}

\paragraph{Step 1.} Let us first bound $\mathscr{Q}_{N} \mathcal{M}_{a}(Z,X)$ by a sum over ranks. For a collision pair $X'$ with particle ranks $i<k$,
\begin{align*} 
    \frac{1}{z_{i} \vee z_{k}}\frac{\displaystyle \mathcal{M}_{a} (\mathbb{T}_{X',u}Z) }{\mathcal{M}_{a}(Z)} &= \frac{1}{z_{k}} u^{-2a} z_{i}^{a} z_{k}^{a} \prod_{n \in [N] \setminus \{i,k\} } \bkt{\frac{z_{n}}{z_{n}+u} }^{a} \\
    &\leq u^{-2a} z_{k}^{2a-1} \prod_{n \in [N] \setminus \{i,k\} } \bkt{\frac{z_{n}}{z_{n}+u} }^{a},
\end{align*}
where we have used the fact that $z_{i} \leq z_{k}$. For each $n > k$, apply
\begin{align*}
    \frac{z_{n}}{z_{n}+u} \leq 1,
\end{align*}
and for all $n\leq k$ so that $n \neq i$, the inequality $z_{n} \leq z_{k}$ gives
\begin{align*}
    \frac{z_{n}}{z_{n}+u} \leq \frac{z_{k}}{z_{k}+u}.
\end{align*}
We conclude
\begin{align} \label{equation:ratio_M_a}
    \frac{1}{z_{i} \vee z_{k}}\frac{\displaystyle \mathcal{M}_{a} (\mathbb{T}_{X',u}Z) }{\mathcal{M}_{a}(Z)}  
    \leq u^{-2a} z_{k}^{2a-1} \bkt{\frac{z_{k}}{z_{k}+u} }^{a(k-2)}.
\end{align}
For each fixed $k\geq3$, there are $k-1$ possible smaller ranks $i$. Summing over the collision pairs therefore yields
\begin{align} \label{equation:M_estimate}
    \frac{\displaystyle \mathscr{Q}_{N}\mathcal{M}_{a}(Z,X)}{\mathcal{M}_{a}(Z,X)} &\leq \frac{1}{M_{N}-1} \sum_{k=3}^{N} (k-1) \int_{0}^{1} du u^{-2a} z_{k}^{2a-1} \bkt{\frac{z_{k}}{z_{k}+u} }^{a(k-2)} \nonumber\\
    &= \frac{1}{M_{N}-1} \sum_{k=3}^{N} (k-1) \int_{0}^{\frac{1}{z_{k}}} du u^{-2a} (1+u)^{-a(k-2)},
\end{align}
where we have used a change of variable in the last equality. 

\paragraph{Step 2.}
To estimate the right-hand side of \eqref{equation:M_estimate}, let us concentrate on the integral with respect to $du$. To bound it, consider the following two cases: 
\begin{align*}
    3\leq k \leq k_{a} \wedge N \quad \text{(lower rank)} \quad \text{and} \quad k_{a}+1 \leq k \leq N \quad \text{(large rank)}, \qquad k_{a} := \biggl  \lfloor \frac{2}{a} \biggl \rfloor.
\end{align*}
For the large rank case, since the Beta function can be expressed as (see \cite[(6.2.1)]{abramowitz1965handbook}):
\begin{align} \label{equation:beta_function_formula}
    B(p,q) = \int_{0}^{\infty}   dx \frac{x^{p-1}}{(1+x)^{p+q}},
\end{align}
it follows that
\begin{align} \label{equation:M_estimate_large_rank}
    \int_{0}^{\frac{1}{z_{k}}} du u^{-2a} (1+u)^{-a(k-2)} &\leq \int_{0}^{\infty} du u^{-2a} (1+u)^{-a(k-2)} 
    = B(1-2a,ak-1) \nonumber\\
    &\leq  \frac{C\Gamma(1-2a)}{(ak)^{(1-2a)}} 
    \leq \frac{C}{a k^{(1-2a)}},
\end{align}
where we have used 
\begin{align*}
    ak-1 \geq \frac{ak}{2} \geq 1 \quad \forall k_{a} < k \leq N
\end{align*}
and Lemma \ref{lemma:beta_property_1} in the second-to-last inequality, and that both of $\Gamma(1-2a)$ and $a^{a}$ are uniformly bounded on $0<a<1/3$ in the last inequality. The proof of Lemma \ref{lemma:beta_property_1} is postponed to Appendix \ref{section:beta_property_1}.

\begin{lemma} \label{lemma:beta_property_1}
There exists an absolute constant $C>0$ such that
\begin{align*}
    B(p,k) \leq 
    \begin{cases}
        \displaystyle \frac{C}{k}, \quad &\text{if } 0<p<1, \quad 0<k\leq p;\\
        \displaystyle \frac{C\Gamma(p)}{k^{p}}, \quad &\text{if } 0<p<1, \quad 0<p\leq k.
    \end{cases}
\end{align*}
\end{lemma}

\paragraph{Step 3.}
Let us now turn to the lower rank case. We further decompose the lower rank case into
\begin{align*}
    z_{k} \geq 1 \quad \text{and} \quad 0<z_{k}<1.
\end{align*}
For the case of $z_{k} \geq 1$, notice that
\begin{align} \label{equation:M_estimate_1}
     \int_{0}^{\frac{1}{z_{k}}} du u^{-2a} (1+u)^{-a(k-2)} \leq \int_{0}^{1} du u^{-2a} \leq C,
\end{align}
where we have used $0<a<1/3$ in the last inequality.  We now focus on the case of $0<z_{k}<1$. In this case, we have
\begin{align} \label{equation:M_estimate_2_2}
    \int_{0}^{\frac{1}{z_{k}}} du u^{-2a} (1+u)^{-a(k-2)} \leq  C+\int_{1}^{\frac{1}{z_{k}}} du u^{-2a} (1+u)^{-a(k-2)}
    \leq  C+\int_{1}^{\frac{1}{z_{k}}} du u^{-ak}, 
\end{align}
where we have used $0<a<1/3$ in the second-to-last inequality. By \eqref{condition:k}, we know that
\begin{align*}
    u^{-ak} \leq u^{-5a/2} \quad \forall u>1,
\end{align*}
which implies
\begin{align} \label{equation:M_estimate_2_1}
    \int_{1}^{\frac{1}{z_{k}}} du u^{-ak} \leq \int_{1}^{\frac{1}{z_{k}}} du u^{-5a/2} \leq  C\frac{1}{z_{k}^{\mathfrak{q}_{a}}}.
\end{align}
Therefore, combining \eqref{equation:M_estimate_1},  \eqref{equation:M_estimate_2_1}, and \eqref{equation:M_estimate_2_2}, we know that, for all $3\leq k \leq k_{a}$ and $z_{k} > 0$, 
\begin{align} \label{equation:M_estimate_low_rank}
    \int_{0}^{\frac{1}{z_{k}}} du u^{-2a} (1+u)^{-a(k-2)}
    \leq C\bkt{ 1 + \frac{1}{z_{k}^{\mathfrak{q}_{a}}}}.
\end{align}

\paragraph{Step 4.}
Combining \eqref{equation:M_estimate_large_rank} and \eqref{equation:M_estimate_low_rank}, we conclude
\begin{align} \label{equation:transition_matrix}
    \frac{\displaystyle \mathscr{Q}_{N}\mathcal{M}_{a}(Z,X)}{\mathcal{M}_{a}(Z,X)} &\leq \frac{C}{M_{N}-1}\bkt{ \sum_{k=3}^{k_{a} \wedge N} k\bkt{ 1 + \frac{1}{z_{k}^{\mathfrak{q}_{a}}}}+\frac{1}{a}\sum_{k = k_{a}+1}^{N} k k^{-(1-2a)}} \nonumber\\
    &\leq \frac{C}{N^{2}}\bkt{ (k_{a}\wedge N)^{2} + (k_{a} \wedge N)\mathcal{A}_{a}(Z)+\frac{1}{a} N^{1+2a}},
\end{align}
which implies
\begin{align*}
    \mathscr{Q}_{N}\mathcal{M}_{a}(Z,X) \leq C\bkt{\frac{1}{a^{2}N^{2}}+\frac{1}{a N^{1-2a}}} \mathcal{H}_{a}(Z).
\end{align*}
The proof of \eqref{lemma:contractiv_inequality_P} is complete. 

\subsection{The mixed weight function}
We now prove the mixed-weight estimate \eqref{lemma:contractiv_inequality_P_S}. Recall that 
\begin{align} \label{equation:section_mixed_weight_function_6}
    \mathscr{Q}_{N} (\mathcal{A}_{a}\mathcal{M}_{a})(Z,X) = \frac{1}{M_{N}-1} \sum_{X'\in \mathcal{E}_{N} \setminus \{X\} } \frac{1}{Z_{X'(1)} \vee Z_{X'(2)}} \int_{0}^{1} (\mathcal{A}_{a}\mathcal{M}_{a})(\mathbb{T}_{X',u}Z,X') du.
\end{align}
Fix $N\geq3$, $0<a<1/3$, and $(Z,X)\in\mathcal{C}_{N}\otimes'\mathcal{E}_{N}$. We use the coordinate ordering \eqref{order_z} and particle ranks \eqref{particle_rank} introduced in the preceding subsection, and recall that $\mathfrak{q}_{a} = 1-5a/2$.
\paragraph{Step 1.} For a collision pair $X'$ whose particles have ranks \eqref{particle_rank}, we have
\begin{align*}
    \mathcal{A}_{a}(\mathbb{T}_{X',u} Z) = 2u^{-\mathfrak{q}_{a}} + \sum_{n \in [N] \setminus \{i,k\}} \bkt{\frac{1}{z_{n}+u}}^{\mathfrak{q}_{a}}.
\end{align*}
Combining this with the ratio estimate \eqref{equation:ratio_M_a} from the preceding subsection and summing over the $k-1$ possible smaller ranks gives
\begin{align*}
    &\frac{\displaystyle \mathscr{Q}_{N} (\mathcal{M}_{a} \mathcal{A}_{a})(Z,X)}{\mathcal{M}_{a}(Z,X)} \\
    &\leq \frac{C}{M_{N}-1} \sum_{k=3}^{N} (k-1) \int_{0}^{1} du u^{-2a} z_{k}^{2a-1} \bkt{\frac{z_{k}}{z_{k}+u} }^{a(k-2)} \bkt{ 2u^{-\mathfrak{q}_{a} } + \sum_{n \in [N] \setminus \{k\}} \frac{1}{(z_{n}+u)^{\mathfrak{q}_{a}}} }\\
    &\leq I_{1}+I_{2},\\
    &I_{1} := \frac{C}{M_{N}-1} \sum_{k=3}^{N} (k-1) \int_{0}^{1} du u^{-2a-\mathfrak{q}_{a}} z_{k}^{2a-1} \bkt{\frac{z_{k}}{z_{k}+u} }^{a(k-2)}, \\
    &I_{2} := \frac{C}{M_{N}-1} \sum_{k=3}^{N} (k-1) \sum_{n \in [N] \setminus \{k\}}  \int_{0}^{1} du u^{-2a} z_{k}^{2a-1} \bkt{\frac{z_{k}}{z_{k}+u} }^{a(k-2)} \frac{1}{(z_{n}+u)^{\mathfrak{q}_{a}}}.
\end{align*}

\paragraph{Step 2.}
Let us first analyze the first term $I_{1}$.  The integral corresponding to $du$ can be bounded by 
\begin{align*}
    \int_{0}^{1} du u^{-2a-\mathfrak{q}_{a}} z_{k}^{2a-1} \bkt{\frac{z_{k}}{z_{k}+u} }^{a(k-2)}
    &\leq z_{k}^{-\mathfrak{q}_{a}} \int_{0}^{\frac{1}{z_{k}}} du u^{-2a-\mathfrak{q}_{a}} \bkt{\frac{1}{1+u}}^{a(k-2)}\\
    &\leq z_{k}^{-\mathfrak{q}_{a}} B(1-2a-\mathfrak{q}_{a},a(k-2)-(1-2a-\mathfrak{q}_{a}))\\
    &= z_{k}^{-\mathfrak{q}_{a}} B(1-2a-\mathfrak{q}_{a},ak-1+\mathfrak{q}_{a}).
\end{align*}
Both parameters are positive, since
\begin{align*}
    1-2a-\mathfrak{q}_{a} = 1-2a-\bkt{1-\frac{5}{2}a} = \frac{1}{2} a > 0, \quad ak-1+\mathfrak{q}_{a} \geq 3a-1+\bkt{1-\frac{5}{2}a} = \frac{1}{2}a > 0.
\end{align*}
By \eqref{condition:k}, we know that 
\begin{align} \label{parameter_lower_bound}
    ak-1+\mathfrak{q}_{a} = a \bkt{k - \frac{5}{2}} \geq \frac{ak}{6}.
\end{align}
This together with Lemma \ref{lemma:beta_property_1} gives
\begin{align*}
    \int_{0}^{1} du u^{-2a-\mathfrak{q}_{a}} z_{k}^{2a-1} \bkt{\frac{z_{k}}{z_{k}+u} }^{a(k-2)}
    \leq z_{k}^{-\mathfrak{q}_{a}} \frac{C\Gamma(a/2)}{(a k)^{a/2}} \leq C \frac{z_{k}^{-\mathfrak{q}_{a}}}{ak^{a/2}}.
\end{align*}
In the last inequality, we have used $\Gamma(a/2) \leq c/a$ and the  uniform boundedness of $a^{-a}$  on $0<a<1/3$.
Consequently, 
\begin{align} \label{inequality:I_1_path_functional}
    I_{1} \leq \frac{C}{M_{N}-1} \sum_{k=3}^{N} k z_{k}^{-\mathfrak{q}_{a}} \frac{1}{ak^{a/2}} \leq \frac{C}{a N^{1+a/2}} \mathcal{A}_{a}(Z).
\end{align}

\paragraph{Step 3}
We next estimate $I_{2}$, distinguishing between $n<k$ and $n>k$. 
\paragraph{Step 3-1} We begin by bounding the integral with respect to $du$. For the case of $n<k$, we have $z_{n} \leq z_{k}$, which implies
\begin{align*}
    \frac{z_{n}}{z_{n}+u} \leq \frac{z_{k}}{z_{k}+u}.
\end{align*}
Similarly, by \eqref{equation:beta_function_formula} and Lemma \ref{lemma:beta_property_1}, the integral with respect to $du$ can be bounded by
\begin{align} \label{equation:section_mixed_weight_function_1}
    \int_{0}^{1} du u^{-2a} z_{k}^{2a-1} \bkt{\frac{z_{k}}{z_{k}+u} }^{a(k-2)} \frac{1}{(z_{n}+u)^{\mathfrak{q}_{a}}} 
    &=  z_{k}^{2a-1} z_{n}^{-\mathfrak{q}_{a}} \int_{0}^{1} du u^{-2a}  \bkt{\frac{z_{k}}{z_{k}+u} }^{a(k-2)} \bkt{\frac{z_{n}}{z_{n}+u}}^{\mathfrak{q}_{a}} \nonumber\\
    &\leq C z_{n}^{-\mathfrak{q}_{a}} B(1-2a,ak+\mathfrak{q}_{a} -1) \leq  \frac{C z_{n}^{-\mathfrak{q}_{a}}}{a k^{1-2a}}.
\end{align}
Here we have used \eqref{parameter_lower_bound} and that $\Gamma(1-2a)$ is uniformly bounded on $0<a<1/3$.

For the case of $n > k$, we know that
\begin{align*}
    z_{k} \leq z_{n}.
\end{align*}
Similarly, by \eqref{equation:beta_function_formula} and Lemma \ref{lemma:beta_property_1},  the integral with respect to $du$ can be bounded by
\begin{align} \label{equation:section_mixed_weight_function_2}
    &\int_{0}^{1} du u^{-2a} z_{k}^{2a-1} \bkt{\frac{z_{k}}{z_{k}+u} }^{a(k-2)} \frac{1}{(z_{n}+u)^{\mathfrak{q}_{a}}} \leq \int_{0}^{1} du u^{-2a} z_{k}^{2a-1} \bkt{\frac{z_{k}}{z_{k}+u} }^{a(k-2)} \frac{1}{(z_{k}+u)^{\mathfrak{q}_{a}}} \nonumber\\
    &= z_{k}^{2a-1-\mathfrak{q}_{a}} \int_{0}^{1} du u^{-2a}  \bkt{\frac{z_{k}}{z_{k}+u} }^{a(k-2)+\mathfrak{q}_{a}} 
    \leq z_{k}^{-\mathfrak{q}_{a}} B(1-2a,ak+\mathfrak{q}_{a} -1) 
    \leq \frac{C z_{k}^{-\mathfrak{q}_{a}} }{a k^{1-2a}}.
\end{align}

\paragraph{Step 3-2} As a result, combining \eqref{equation:section_mixed_weight_function_1} and \eqref{equation:section_mixed_weight_function_2} gives
\begin{align} \label{equation:section_mixed_weight_function_3}
    I_{2} &\leq \frac{C}{M_{N}-1} \sum_{k=3}^{N} (k-1) \sum_{n: 1\leq n<k } z_{n}^{-\mathfrak{q}_{a}} \frac{1}{a k^{1-2a}}+\frac{C}{M_{N}-1} \sum_{k=3}^{N} (k-1) \sum_{n: k < n\leq N}  z_{k}^{-\mathfrak{q}_{a}} \frac{1}{a k^{1-2a}}.
\end{align}
For the first term on the right-hand side of \eqref{equation:section_mixed_weight_function_3}, interchanging the sums yields
\begin{align*}
    &\frac{C }{M_{N}-1} \sum_{k=3}^{N} (k-1) \sum_{n: 1\leq n<k } z_{n}^{-\mathfrak{q}_{a}} \frac{1}{a k^{1-2a}}\leq \frac{C}{a N^{2}} \sum_{n=1}^{N} z_{n}^{-\mathfrak{q}_{a}} \sum_{k: N\geq k\geq \max\{3,n\}} \frac{k}{k^{1-2a}}
    \leq \frac{C}{a N^{1-2a}} \mathcal{A}_{a}(Z).
\end{align*}
Similarly, for the second term on the right-hand side of \eqref{equation:section_mixed_weight_function_3},
\begin{align*}
    &\frac{C }{M_{N}-1} \sum_{k=3}^{N} (k-1) \sum_{n: k < n\leq N}  z_{k}^{-\mathfrak{q}_{a}} \frac{1}{a k^{1-2a}}
    \leq \frac{C }{a N^{2}} \sum_{k=3}^{N} k^{2a} (N-k) z_{k}^{-\mathfrak{q}_{a}} 
    \leq \frac{C}{a N^{1-2a}} \mathcal{A}_{a}(Z).
\end{align*}
Consequently,
\begin{align} \label{inequality:I_2_path_functional}
    I_{2} \leq \frac{C}{a N^{1-2a}} \mathcal{A}_{a}(Z).
\end{align}

\paragraph{Step 4.}
Therefore, combining \eqref{inequality:I_1_path_functional} and \eqref{inequality:I_2_path_functional} yields
\begin{align*}
    &\frac{\displaystyle \mathscr{Q}_{N} (\mathcal{M}_{a} \mathcal{A}_{a})(Z,X)}{\mathcal{M}_{a}(Z,X)}  \leq \frac{C}{a N^{1-2a}} \mathcal{A}_{a}(Z).
\end{align*}
The proof is complete.

\section{Estimates for the path functional} \label{section:estimate_for_path_functional}

In this section, we prove Lemma \ref{lemma:estimate_for_path_functional}.  To begin with, we change the terminal function $\mathds{1}$ to $\mathcal{H}_{a}$. Since
\begin{align*}
    \mathbf{Z}_{x}(m) \leq  m+1 \quad \forall x \in [N] \quad \text{(a.s.)},
\end{align*}
we know that
\begin{align*}
    \mathcal{M}_{a}(\mathbf{Z}_{\cdot}(m)) &= \prod_{x \in [N]} \frac{1}{\mathbf{Z}_{x}(m)^{a}} \geq \frac{1}{(m+1)^{a N}},\quad \mathcal{A}_{a}(\mathbf{Z}_{\cdot}(m)) &= \sum_{x \in [N]} \frac{1}{\mathbf{Z}_{x}(m)^{\mathfrak{q}_{a}}} \geq \frac{N}{(m+1)^{\mathfrak{q}_{a}}} \quad \text{(a.s.)}.
\end{align*}
It follows that
\begin{align*}
    \frac{(m+1)^{aN+\mathfrak{q}_{a}}}{N} \mathcal{H}_{a}(\mathbf{Z}_{\cdot}(m)) = \frac{(m+1)^{aN+\mathfrak{q}_{a}}}{N} \mathcal{M}_{a}(\mathbf{Z}_{\cdot}(m))(1+\mathcal{A}_{a}(\mathbf{Z}_{\cdot}(m))) \geq 1 \quad \text{(a.s.)},
\end{align*}
which implies
\begin{align} \label{equation:change_terminal_function}
    \mathds{1} \leq \frac{(m+1)^{aN+\mathfrak{q}_{a}}}{N} \mathcal{H}_{a} (\mathbf{Z}_{\cdot}(m)) \leq 2^{m} \mathcal{H}_{a} (\mathbf{Z}_{\cdot}(m)) \quad \text{(a.s.)}. 
\end{align}
Here we have used that for all $m \geq \mathfrak{m}$, $N\geq 3$, and $0<a<1/3$,
\begin{align*}
    (m+1)^{aN+\mathfrak{q}_{a}} = \exp((aN+\mathfrak{q}_{a}) \log (m+1)) \leq \exp(m\log 2) = 2^{m}.
\end{align*}
Consequently, \eqref{equation:change_terminal_function} together with the positivity of the weighted transition operator gives
\begin{align} \label{eqaution:section_estimate_for_path_functional_0}
    \mathbf{F}_{N}(m) 
    \leq 2^{m} \mathds{E}\bigbkt{ \mathscr{Q}_{N}^{m-1}\mathcal{H}_{a}(\mathbf{Z}_{\cdot}(1),\mathbf{X}_{0})    }.
\end{align}

Now we are ready to bound $\mathbf{F}_{N}(m)$. By using \eqref{eqaution:section_estimate_for_path_functional_0}, we obtain
\begin{align} \label{equation:section_main_result_proof_idea_0}
    \mathbf{F}_{N}(m) 
    &\leq (2C_{\eqref{equation:main_contractiv_inequality}})^{m} \bkt{\frac{1}{a^{2} N^{2}}+ \frac{1}{a N^{1-2a}}}^{m-1} \mathds{E}\bigbkt{ \mathcal{H}_{a}(\mathbf{Z}_{\cdot}(1),\mathbf{X}_{0})    } \nonumber\\
    &\leq a N (2C_{\eqref{equation:main_contractiv_inequality}})^{m} \bkt{\frac{1}{a^{2} N^{2}}+\frac{1}{aN^{1-2a}}}^{m} \mathds{E}\bigbkt{ \mathcal{H}_{a}(\mathbf{Z}_{\cdot}(1),\mathbf{X}_{0})    },
\end{align}
where we have used Proposition \ref{proposition:main_contractiv_inequality}. Let us now estimate the expectation on the right-hand side of \eqref{equation:section_main_result_proof_idea_0}. To do this, note that 
\begin{align*}
    \mathbf{Z}_{\mathbf{X}_{0}(1)}(1) = \mathbf{Z}_{\mathbf{X}_{0}(2)}(1) = U_{1} \quad \text{(a.s.)}
\end{align*}
and that 
\begin{align} \label{equation:lower_bound_U_0}
    \mathbf{Z}_{x}(1) = U_{1}+U_{0} \geq \mathfrak{u}>0 \quad \forall x\not \in \{\mathbf{X}_{0}(1),\mathbf{X}_{0}(2)\} \quad  \text{(a.s.)}.
\end{align}
This implies 
\begin{align*}
    \mathcal{M}_{a}(\mathbf{Z}_{\cdot}(1)) \leq \mathfrak{u}^{-a(N-2)} U_{1}^{-2 a}, \quad \mathcal{A}_{a}(\mathbf{Z}_{\cdot}(1)) \leq N U_{1}^{-\mathfrak{q}_{a}} \quad \ \text{(a.s.)}.
\end{align*}
It follows that 
\begin{align} \label{equation:section_main_result_proof_idea_3}
    \mathds{E}[\mathcal{H}_{a}(\mathbf{Z}_{\cdot}(1),\mathbf{X}_{0}) ]
    \leq N \mathfrak{u}^{-a(N-2)} \mathds{E}[U_{1}^{-2a-\mathfrak{q}_{a}}]+\mathfrak{u}^{-a(N-2)} \mathds{E}[U_{1}^{-2a}] \leq \frac{4}{a} \mathfrak{u}^{-aN} N,
\end{align}
where we have used $2a+\mathfrak{q}_{a} = 1-a/2$ in the last inequality. 
Combining \eqref{equation:section_main_result_proof_idea_0} and \eqref{equation:section_main_result_proof_idea_3}
gives \eqref{equation:estimate_for_path_functional}. The proof is complete.

\section{Tail probability of the SHF} \label{section:tail_probability_SHF}

To prove the tail bound in Theorem \ref{corollary:two_side_bound_tail_probability},  we first extend Ganguly and Nam's lower moment bound as follows, where the proof is given in Section~\ref{section:upper_bounds}.

\begin{theorem}\label{theorem:localization}
Assume that $\theta \in \mathbb{R}$. For each $\varphi,g \in C_{c}(\mathbb{R}^{2};\mathbb{R}_{\geq 0})$, not identically zero, and $t\geq 1$, there exist $c_{\eqref{equation:localization}} > 0$ and $N_{0}  \geq 3$ such that
\begin{align} \label{equation:localization}
    \mathbb{E}[\mathscr{Z}_{t}^{\theta}(\varphi,g)^{N}] \geq \exp(t\exp(c_{\eqref{equation:localization}}N)) \quad \forall N \geq N_{0}.
\end{align}
Here $c_{\eqref{equation:localization}}$ is independent of $t$ and $N$.
\end{theorem}

\begin{proof} [Proof of Theorem \ref{corollary:two_side_bound_tail_probability}]

Let $(h,f) \in \eqref{eq:UBH}$, and take $z$ sufficiently large so that  
\begin{align*}
    N = \biggl \lfloor \frac{1}{C_{\eqref{equation:main_result_higher_moment}}} \frac{\log((\log z) / t)}{\log \log((\log z) / t)} \biggr \rfloor \geq 3.
\end{align*}
Using Markov’s inequality, \begin{align*}
    \mathbb{P}(\mathscr{Z}_{t}^{\theta}(h,f) \geq z) \leq \frac{\mathbb{E}[\mathscr{Z}_{t}^{\theta}(h,f)^{N}]}{z^{N}}.
\end{align*}
Since $N\ge 3$,  the required upper bound \eqref{equation:upper_bound_tail_probability} follows from the upper moment bound \eqref{equation:main_result_higher_moment}.

Now let $h,f\in C_{c}(\mathbb{R}^{2};\mathbb{R}_{\geq0})$
be nonzero. Then $(h,f)$ satisfies \eqref{eq:UBH}, so the upper
moment bound \eqref{equation:main_result_higher_moment} applies. This together with Theorem~\ref{theorem:localization} yields that there exist constants $C,c>0$, which are independent of $t$, such that, for all sufficiently
large real number $p$,
\begin{align} \label{equation:tail_probability_lower_bound_0}
    \exp(t\exp(c p)) \leq \mathbb{E}[\mathscr{Z}_{t}^{\theta}(h,f)^{p}] \leq \exp(t\exp(Cp \log p)).
\end{align}
For sufficiently large $z$, take a sufficiently large real number $p$ such that
\begin{align*}
    z= \exp\bkt{ \frac{t}{p}\exp(cp) }.
\end{align*}
It follows that
\begin{align*}
    &\mathbb{P}(\mathscr{Z}_{t}^{\theta}(h,f) \geq z/2) 
    \geq \mathbb{P}\bkt{\mathscr{Z}_{t}^{\theta}(h,f) \geq \frac{1}{2} \mathbb{E}[\mathscr{Z}_{t}^{\theta}(h,f)^{p}]^{\frac{1}{p}} } \geq \bkt{1-\frac{1}{2^{p}} }^{2} \frac{\mathbb{E}[\mathscr{Z}_{t}^{\theta}(h,f)^{p}]^{2}}{\mathbb{E}[\mathscr{Z}_{t}^{\theta}(h,f)^{2p}]}\\
    &\geq \bkt{1-\frac{1}{2^{p}} }^{2} \frac{\exp(2t\exp(cp))}{\exp(t\exp(C 2p \log (2p)))}
    \geq \bkt{1-\frac{1}{2^{p}} }^{2} \exp(-t\exp(C 2p \log (2p))).
\end{align*}
Here we used \eqref{equation:tail_probability_lower_bound_0} and applied the Paley–Zygmund inequality to $\mathscr{Z}_{t}^{\theta}(h,f)^{p}$.  From $\log (p (\log z)/t )=cp$, we see that
$p\asymp \log (\frac{\log z} t)$ and hence, 
\begin{align} \label{approximation_p}
    p \log p \asymp (\log \log ((\log z)/t))(\log ((\log z)/t)),
\end{align}
which completes the proof of \eqref{equation:lower_bound_tail_probability}.

\end{proof}

\section{Proof of the ground-state energy bounds} \label{section:2d_DBG}

The goal of this section is to prove Theorem~\ref{theorem:2d_DBG}.
For simplicity, we focus on the unrestricted Hamiltonians.
All the results in this section also hold for their bosonic restrictions,
and the same proof ideas apply.

\subsection{Large-time growth of the semigroup} \label{section:lower_bounds}

The lower bound on the ground-state energy follows from the upper moment
bound in Theorem \ref{theorem:main_result_higher_moment}. Since that theorem  imposes restrictions on the test functions, we cannot estimate the operator norm of
$\exp(-t\mathscr{H}_{N}^{\theta})$ directly. Instead, we use the spectral theorem to recover $\inf \mathscr{H}_{N}^{\theta}$ from the large-time growth of the semigroup.  We begin with the following density property.
\begin{lemma} \label{lemma:lemma_spectrum_theorem_below_G}
Suppose that $N \geq 2$. Then $\textup{span}(\mathcal{G})$ is dense in $L^{2}(\mathbb{R}^{2N})$, where
\begin{align*}
    \mathcal{G} := \{\psi \in L^{2}(\mathbb{R}^{2N}): 0 \leq \psi \leq G_{1}^{\otimes N} \}.
\end{align*}
\end{lemma}
\begin{proof}
Fix $u \in L^{2}(\mathbb{R}^{2N})$. To prove $u \in \overline{\text{span}(\mathcal{G})}$, without loss of generality, we may assume that $u \geq 0$. Observe the following cutoff of $u$:
\begin{align} \label{u_cut_off}
    u_{n} := u \wedge (n G_{1}^{\otimes N}) = n \bkt{ \frac{u}{n} \wedge G_{1}^{\otimes N} } \in \text{span}(\mathcal{G}).
\end{align}
The strict positivity of $G_{1}^{\otimes N}$ implies
$u_{n}\uparrow u$ almost everywhere.
Since $0\leq u_{n}\leq u$, the dominated convergence theorem gives
$u_{n}\to u$ in $L^{2}(\mathbb{R}^{2N})$.
This proves the density assertion.
\end{proof}

Combining this density property with the spectral results in
Section \ref{section:spectrum_properties}  yields the following large-time limit.

\begin{lemma} \label{lemma:spectrum_theorem_1}
Assume that $\theta \in \mathbb{R}$. Then for each $N \geq 2$, 
\begin{align} \label{equation:large_time_semigroup}
    \lim_{t\to\infty} \frac{1}{t} \log \langle G^{\otimes N}_{1}, \exp(t \mathscr{L}_{N}^{\theta}) G_{1}^{\otimes N} \rangle_{L^{2}(\mathbb{R}^{2N})} = \sup \sigma(\mathscr{L}_{N}^{\theta}),
\end{align}
where we have set
\begin{align} \label{change_sign}
    \mathscr{L}_{N}^{\theta}  := -\mathscr{H}_{N}^{\theta}.
\end{align}
\end{lemma}
\begin{proof} 
For $\psi\in L^{2}(\mathbb{R}^{2N})$, let $\mu_{\psi}$ denote the
spectral measure
\begin{align*}
    \mu_{\psi}(A)
    :=\langle\psi,
    \mathds{1}_{A}(\mathscr{L}_{N}^{\theta})\psi
    \rangle_{L^{2}(\mathbb{R}^{2N})},
    \quad A\in\mathcal{B}(\mathbb{R}).
\end{align*}
By Lemma \ref{lemma:lemma_spectrum_theorem_support},
\begin{align} \label{limsup_G_0}
    \lim_{t\to\infty}\frac{1}{t}
    \log\langle G_{1}^{\otimes N},
    \exp(t\mathscr{L}_{N}^{\theta})G_{1}^{\otimes N}
    \rangle_{L^{2}(\mathbb{R}^{2N})}
    =\sup\textup{supp}(\mu_{G_{1}^{\otimes N}}).
\end{align}
Since $G_{1}^{\otimes N} \in \mathcal{G}$, applying Lemma
\ref{lemma:sup_representation} then gives
\begin{align} \label{limsup_G}
    \lim_{t\to\infty}\frac{1}{t}
    \log\langle G_{1}^{\otimes N},
    \exp(t\mathscr{L}_{N}^{\theta})G_{1}^{\otimes N}
    \rangle_{L^{2}(\mathbb{R}^{2N})}
    \leq\sup\sigma(\mathscr{L}_{N}^{\theta}).
\end{align}
For the reverse, fix a non-zero $\psi\in\mathcal{G}$.
Since $\exp(t\mathscr{L}_{N}^{\theta})$ is positivity preserving and
$0\leq\psi\leq G_{1}^{\otimes N}$,
\begin{align*}
    \langle G_{1}^{\otimes N},
    \exp(t\mathscr{L}_{N}^{\theta})G_{1}^{\otimes N}
    \rangle_{L^{2}(\mathbb{R}^{2N})}
    \geq
    \langle\psi,\exp(t\mathscr{L}_{N}^{\theta})\psi
    \rangle_{L^{2}(\mathbb{R}^{2N})}
    \quad\forall t>0.
\end{align*}
Taking logarithms, dividing by $t$, and applying
Lemma \ref{lemma:lemma_spectrum_theorem_support} gives
\begin{align*}
    \lim_{t\to\infty}\frac{1}{t}
    \log\langle G_{1}^{\otimes N},
    \exp(t\mathscr{L}_{N}^{\theta})G_{1}^{\otimes N}
    \rangle_{L^{2}(\mathbb{R}^{2N})}
    \geq\sup\textup{supp}(\mu_{\psi}).
\end{align*}
Taking the supremum over $\psi\in\mathcal{G}$ and using
Lemma \ref{lemma:lemma_spectrum_theorem_below_G} and Lemma
\ref{lemma:sup_representation}, we conclude that
\begin{align} \label{liminf_G}
    \lim_{t\to\infty}\frac{1}{t}
    \log\langle G_{1}^{\otimes N},
    \exp(t\mathscr{L}_{N}^{\theta})G_{1}^{\otimes N}
    \rangle_{L^{2}(\mathbb{R}^{2N})}
    \geq\sup_{\psi\in\mathcal{G}}
    \sup\textup{supp}(\mu_{\psi})
    =\sup\sigma(\mathscr{L}_{N}^{\theta}).
\end{align}
Combining \eqref{limsup_G} and \eqref{liminf_G} completes the proof.
\end{proof}

\subsection{Proof of Theorem \ref{theorem:2d_DBG}}

The upper bound on the ground-state energy uses the localized lower moment bound in Theorem~\ref{theorem:localization}, whose proof is given in the next subsection.

Firstly, by Theorem \ref{theorem:main_result_higher_moment}, for each $t \geq 1$ and $N\geq 3$,
\begin{align} \label{equation:top_spectrum_less_moment_exponent}
    \langle G_{1}^{\otimes N}, \exp(t \mathscr{L}_{N}^{\theta}) G_{1}^{\otimes N} \rangle_{L^{2}(\mathbb{R}^{2N})} = \mathbb{E}[\mathscr{Z}_{t}^{\theta}(G_{1},G_{1})^{N}] \leq \exp(t\exp(CN \log N)),
\end{align}
where $C>0$ only depends on $\theta$. Thus, by using Lemma \ref{lemma:spectrum_theorem_1}, taking logarithms, dividing by $t$, and letting $t\to\infty$ then gives
\begin{align*} 
    \sup \sigma(\mathscr{L}_{N}^{\theta}) \leq \exp(C N (\log N)) \quad \forall N \geq 3,
\end{align*}
which implies the desired inequality.

We now turn to the upper bound. Let $\varphi,g \in C_{c}(\mathbb{R}^{2};\mathbb{R}_{\geq 0})$ be non-zero. Then $\varphi \vee g \lesssim G_{1}$. Thus, applying the spectral theorem together with Lemma \ref{lemma:sup_representation} implies that there exists $C = C(\varphi,g)>0$ such that
\begin{align} \label{upper_bound_sp}
    \mathbb{E}[\mathscr{Z}_{1}^{\theta}(\varphi,g)^{N}] = \langle \varphi^{\otimes N}, \exp( \mathscr{L}_{N}^{\theta}) g^{\otimes N} \rangle_{L^{2}(\mathbb{R}^{2N})} 
    \leq \exp( \sup \sigma(\mathscr{L}_{N}^{\theta})) C^{N} \quad \forall N\geq 3.
\end{align}
Therefore, by using the lower moment bounds in Theorem \ref{theorem:localization} with $t = 1$, there exists $c = c(\theta,\varphi,g)>0$ such that for all sufficiently  large $N$, $\exp(\exp(cN))\leq \mathbb{E}[\mathscr{Z}_{1}^{\theta}(\varphi,g)^{N}]$. This together with \eqref{upper_bound_sp} then gives $\exp(\exp(\widetilde{c}N))\leq \exp(\sup \sigma(\mathscr{L}_{N}^{\theta}))$, where $\widetilde{c} = \widetilde{c}(\theta,\varphi,g) > 0$. Taking logarithms then yields the upper bound in Theorem \ref{theorem:2d_DBG}. 

\subsection{Proof of the localized lower moment bound} \label{section:upper_bounds}

The proof of Theorem \ref{theorem:localization} relies on a comparison
between translated terminal test functions.
For $a\in\mathbb{R}^{2}$, define
\begin{align*}
    g_{a}(z):=g(z-a),\quad z\in\mathbb{R}^{2}.
\end{align*}
The following proposition compares the moment associated with $g_{a}$
with the one associated with $g$.
All $L^{N}$-norms of random variables below are taken with respect to
$\mathbb{P}$.

\begin{proposition} [Translation estimate]\label{proposition:change_of_ceneter_a}
Assume that $\theta \in \mathbb{R}$ and $\varphi,g \in C_{c}(\mathbb{R}^{2};\mathbb{R}_{\geq 0})$ are non-zero. Then there exists $D = D(\varphi,g)>0$  such that for each $N\geq 3$, $a \in \mathbb{R}^{2}$, and $t\geq 1$,
\begin{align} \label{estimate:change_of_ceneter_a}
    ||\mathscr{Z}_{t}^{\theta}(\varphi,g_{a}) ||_{L^{N}} \leq \mathbb{D}(a,t) ||\mathscr{Z}_{t}^{\theta}(\varphi,g)||_{L^{N}},
\end{align}
where 
\begin{align*}
     \mathbb{D}(a,t) := \exp\bkt{\frac{-(|a|-D)^{2}}{2t} } \exp\bkt{\frac{D^{2}}{2t} }.
\end{align*}
\end{proposition}

\begin{proof} [Proof of Theorem \ref{theorem:localization}]
The starting point is the lower moment bound of Ganguly and Nam \cite[Theorem 1.2]{ganguly2025sharp}, in which one endpoint is a constant function. Set
\begin{align*}
    \mathscr{Z}_{t}^{\theta}(\varphi) := \mathscr{Z}_{t}^{\theta}(\varphi,\mathds{1}).
\end{align*}
Assume that $\theta \in \mathbb{R}$. For each $\varphi \in C_{c}(\mathbb{R}^{2};\mathbb{R}_{\geq 0})$ such that $\varphi(0) > 0$,  their result gives constants $c_{\eqref{equation:localization_0}}>0$ and $N_{0}\geq 3$ such that
\begin{align} \label{equation:localization_0}
    \exp(\exp(c_{\eqref{equation:localization_0}}N))  \leq \mathbb{E}[\mathscr{Z}_{1}^{\theta}(\varphi)^{N}] \quad \forall N \geq N_{0}.
\end{align}

We now relax these assumptions. Spatial translation invariance of the SHF
\cite[Corollary 1.6]{tsai2024stochastic} allows the condition
$\varphi(0)>0$ to be replaced by $\varphi\not\equiv 0$.
Moreover, \cite[Remark 1.3]{ganguly2025sharp} gives the same lower bound
at every fixed time $t>0$, with the threshold $N_{0}$ allowed to
depend on $t$.
For any $0<c'<c_{\eqref{equation:localization_0}}$ and fixed $t\geq 1$,
\begin{align*}
    \exp\bigl(\exp(c_{\eqref{equation:localization_0}}N)\bigr)
    \geq\exp\bigl(t\exp(c'N)\bigr)
\end{align*}
for all sufficiently large $N$.
Consequently, there exists
$c_{\eqref{equation:localization_1}}
=c_{\eqref{equation:localization_1}}(\theta,\varphi)>0$
such that, for every $t\geq 1$, there is an integer
$N_{0}=N_{0}(\theta,\varphi,t)\geq 3$ for which
\begin{align} \label{equation:localization_1}
    \mathbb{E}\bigl[\mathscr{Z}_{t}^{\theta}(\varphi)^{N}\bigr]
    \geq\exp\bigl(t\exp(c_{\eqref{equation:localization_1}}N)\bigr)
    \qquad\forall N\geq N_{0}.
\end{align}

To localize the remaining endpoint, fix non-zero functions $\varphi,g \in C_{c}(\mathbb{R}^{2},\mathbb{R}_{\geq 0})$. Tonelli’s theorem gives
\begin{align*}
     \mathscr{Z}_{t}^{\theta}(\varphi,\mathds{1}) = ||g||_{L^{1}(\mathbb{R}^{2})}^{-1} \cdot \int_{\mathbb{R}^{2}} da \mathscr{Z}_{t}^{\theta}(\varphi,g_{a}).
\end{align*}
This together with Minkowski’s inequality yields
\begin{align} \label{equation:Minkowski}
    ||\mathscr{Z}_{t}^{\theta}(\varphi)||_{L^{N}} \leq ||g||_{L^{1}(\mathbb{R}^{2})}^{-1} \cdot\int_{\mathbb{R}^{2}} da ||\mathscr{Z}_{t}^{\theta}(\varphi,g_{a}) ||_{L^{N}}.
\end{align}
One easily shows that there exists a constant $C = C(\varphi,g)>0$ such that 
\begin{align*}
    \int_{\mathbb{R}^{2}} da \mathbb{D}(a,t) \leq Ct \quad \forall t\geq 1.
\end{align*}
This together with \eqref{equation:Minkowski}, Proposition \ref{proposition:change_of_ceneter_a}, and \eqref{equation:localization_1}, implies Theorem \ref{theorem:localization}. Here we have used that there exists $c = c(\theta,\varphi,g) > 0$ such that for all $t\geq 1$ and $N\geq N_{0}$, $\exp(t\exp(c N)) \leq \exp(t\exp(c_{\eqref{equation:localization_1}}N))/(Ct)^{N}$. 
\end{proof}

We now turn to the translation estimate in Proposition~\ref{proposition:change_of_ceneter_a}. Fix $a \in \mathbb{R}^{2}$. By \eqref{equation:connection_SHF_DBG}, we have
\begin{align} \label{equation:connection_moment_on_a}
    \mathbb{E}[\mathscr{Z}_{t}^{\theta}(\varphi,g_{a})^{N}] = \langle \varphi^{\otimes N}, \exp(t\mathscr{L}_{N}^{\theta}) g_{a}^{\otimes N} \rangle_{L^{2}(\mathbb{R}^{2N})}. 
\end{align}
To isolate the dependence on $a$, we separate the center-of-mass motion from the relative motion using the Jacobi coordinates
\begin{align}  \label{equation:jacobi}
    q_{0} = \frac{1}{\sqrt{N}} \sum_{k= 1}^{N} x_{k}, \quad q_{k} = \frac{1}{\sqrt{k(k+1)}}\bkt{\sum_{i=1}^{k} x_{i} - k x_{k+1} } \quad \forall k = 1,\ldots,N-1.
\end{align}
Write
\begin{align*}
    q = (q_{0},q^{\perp}), \quad q^{\perp} := (q_{1},\ldots,q_{N-1}) \in \mathbb{R}^{2(N-1)}.
\end{align*}
We denote the coordinate transformation by $x\mapsto q(x)$ and its
inverse by $q\mapsto x(q)$.
\begin{proposition}[Semigroup decomposition in Jacobi coordinates] \label{proposition:decompistion}
For each $\theta \in \mathbb{R}$, $N\geq 3$, and $t>0$, there exists a non-negative measurable function $\exp(t \mathscr{R}_{N}^{\theta})(q^{\perp},q'^{\perp})$ such that for each $f,f' \in C_{c}(\mathbb{R}^{2N})$,
\begin{align} \label{equation:decomposition}
    &\int_{\mathbb{R}^{4\times N}} dxdx' f(x) \exp(t \mathscr{L}_{N}^{\theta})(x,x') f'(x') \nonumber\\
    &= \int_{\mathbb{R}^{4\times N}} dxdx' f(x) G_{t}(q_{0}(x) - q'_{0}(x'))\exp(t \mathscr{R}_{N}^{\theta})(q^{\perp}(x),q^{\perp}(x')) f'(x').
\end{align}
\end{proposition}

The proof of this decomposition is deferred to
Section \ref{section:decomposition}.
The translation estimate follows by comparing the center-of-mass
heat kernels.

\begin{proof} [Proof of Proposition \ref{proposition:change_of_ceneter_a}]
Fix $a\in\mathbb{R}^{2}$ and write $\mathbf{a} = (a,\ldots,a) \in \mathbb{R}^{2N}$. Using Proposition \ref{proposition:decompistion} gives
\begin{align} \label{equation:Jacobi_representation_0}
    &\langle \varphi^{\otimes N}, \exp(t\mathscr{L}_{N}^{\theta}) g_{a}^{\otimes N} \rangle_{L^{2}(\mathbb{R}^{2N})} \nonumber\\
    &= \int_{\mathbb{R}^{4N}} dxdx' \varphi^{\otimes N}(x) G_{t}(q_{0}(x)-q_{0}(x')) \exp(t \mathscr{R}_{N}^{\theta})(q^{\perp}(x),q^{\perp}(x')) g_{0}^{\otimes N}(x'-\mathbf{a}) \nonumber\\
    &= \int_{\mathbb{R}^{4N}} dxdx' \varphi^{\otimes N}(x) G_{t}(q_{0}(x)-q_{0}(x'+\mathbf{a})) \exp(t \mathscr{R}_{N}^{\theta})(q^{\perp}(x),q^{\perp}(x'+\mathbf{a})) g_{0}^{\otimes N}(x') \nonumber\\
    &= \int_{\mathbb{R}^{4N}} dxdx' \varphi^{\otimes N}(x) G_{t}(q_{0}(x)-q_{0}(x'+\mathbf{a})) \exp(t \mathscr{R}_{N}^{\theta})(q^{\perp}(x),q^{\perp}(x')) g_{0}^{\otimes N}(x').
\end{align}
In the second equality, we changed variables; in the third, we used
$q^{\perp}(x'+\mathbf{a})=q^{\perp}(x')$. Moreover, \eqref{equation:jacobi} gives
\begin{align} \label{application_q0_0}
    G_{t}(q_{0}(x)-q_{0}(x'+\mathbf{a})) = G_{t}(q_{0}(x)-q_{0}(x')-\sqrt{N}a).
\end{align}
Since $\varphi$ and $g$ have compact support,  there exists $D = D(\varphi,g_{0})>0$ such that
\begin{align} \label{application_q0_1}
    \sup_{x_{i} \in \text{supp}(\varphi),\; x'_{i} \in \text{supp}(g_{0})} |q_{0}(x) - q_{0}(x')| \leq \sqrt{N} D, \quad D = D(\varphi,g_{0}),
\end{align}
which implies
\begin{align} \label{application_q0_2}
    &G_{t}(q_{0}(x)-q_{0}(x')-\sqrt{N}a) \nonumber\\
    &\leq G_{t}(q_{0}(x)-q_{0}(x')) \exp\bkt{\frac{- |a|^{2} N}{2t} } \exp\bkt{\frac{|a|ND}{t}} = G_{t}(q_{0}(x)-q_{0}(x'))  \mathbb{D}(a,t)^{N}.
\end{align}
Therefore, by the nonnegativity of the relative-coordinate kernel,
\eqref{equation:connection_moment_on_a},
\eqref{equation:Jacobi_representation_0}, and
\eqref{application_q0_0}--\eqref{application_q0_2} imply 
\begin{align*}
    \mathbb{E}[\mathscr{Z}_{t}^{\theta}(\varphi,g_{a})^{N}]  \leq  \mathbb{D}(a,t)^{N} \langle \varphi^{\otimes N}, \exp(t\mathscr{L}_{N}^{\theta}) g_{0}^{\otimes N} \rangle_{L^{2}(\mathbb{R}^{2N})} = \mathbb{D}(a,t)^{N} \mathbb{E}[\mathscr{Z}_{t}^{\theta}(\varphi,g_{0})^{N}].
\end{align*}
The proof of Proposition \ref{proposition:change_of_ceneter_a} is complete. 
\end{proof}

\subsection{Proof of the semigroup decomposition} \label{section:decomposition}

We first establish the semigroup decomposition for the regularized operators
\begin{align} \label{equation:L_varepsilon}
    \mathscr{L}_{N;\varepsilon}^{\theta} := \frac{1}{2}\sum_{i = 1}^{N} \Delta_{x_{i}} + \beta_{\varepsilon}^{2} \sum_{(\mathfrak{i},\mathfrak{i}')\in \mathcal{E}_{N}} R_{\varepsilon}(x_{\mathfrak{i}'} - x_{\mathfrak{i}}), \quad x = (x_{1},\ldots,x_{N}) \in \mathbb{R}^{2N},
\end{align}
where $\beta_{\varepsilon}$ and $R_{\varepsilon}(x)$ are given in \eqref{equation:critical_scaling} and \eqref{eq:R}, respectively. We then pass to the limit in the space of locally finite positive
Borel measures equipped with the vague topology.
The following properties of the Jacobi transformation will be used.
\begin{lemma} \label{lemma:jacobi}
The following properties of the Jacobi transformation hold. 
\begin{enumerate} [label = (\roman*).]
    \item The change-of-coordinates operator
    \begin{align*}
        (\mathcal{U}f)(q):=f(x(q)),
        \quad f\in L^{2}(\mathbb{R}^{2N}),
    \end{align*}
    is unitary, with adjoint
    \begin{align*}
        (\mathcal{U}^{*}\widetilde{f})(x)
        :=\widetilde{f}(q(x)).
    \end{align*}
    \item For each $(\mathfrak{i},\mathfrak{i}')\in \mathcal{E}_{N}$, we have
    \begin{align*}
        x_{\mathfrak{i}'} - x_{\mathfrak{i}} &= A_{\mathfrak{i},\mathfrak{i}'} (q^{\perp}(x)), 
    \end{align*}
    where
    \begin{align*}
        A_{\mathfrak{i},\mathfrak{i}'} (q^{\perp}) := \mathds{1}_{\{\mathfrak{i}\geq 2\}} \sqrt{\frac{\mathfrak{i}-1}{\mathfrak{i}}} q_{\mathfrak{i}-1}- \sum_{k=\mathfrak{i}}^{\mathfrak{i}'-2} \frac{1}{\sqrt{k(k+1)}} q_{k} - \sqrt{\frac{\mathfrak{i}'}{\mathfrak{i}'-1}} q_{\mathfrak{i}'-1}.
    \end{align*}
    \item The Jacobian determinant is $1$, i.e. $ dx_{1}\ldots dx_{N} = dq_{0}dq^{\perp}$.
    \item Let $B = (B_{1},\ldots,B_{N})$ be a $2N$-dimensional Brownian motion starting from $x$. Then $Q := q(B)$ is a $2N$-dimensional Brownian motion starting from $q(x)$.
\end{enumerate}
\end{lemma}

\begin{proof} [Proof of Proposition \ref{proposition:decompistion}]

By the Feynman–Kac formula and parts (ii) and (iv) of Lemma \ref{lemma:jacobi}, a change of measure shows that for each $f,f' \in C_{c}(\mathbb{R}^{2N})$,
\begin{align} \label{FK_0}
    &\int_{\mathbb{R}^{4N}} dxdx' f(x)\exp(t \mathscr{L}_{N;\varepsilon}^{\theta})(x,x') f'(x') \nonumber\\
    &= \int_{\mathbb{R}^{2N}} dx f(x) \mathbf{E}^{B}_{x}\bigbkt{\exp\bkt{\sum_{(\mathfrak{i},\mathfrak{i}')\in \mathcal{E}_{N}} \int_{0}^{t} \beta_{\varepsilon}^{2}R_{\varepsilon}(B_{\mathfrak{i}'}(s)-B_{\mathfrak{i}}(s)) ds } f'(B(t)) } \nonumber\\
    &= \int_{\mathbb{R}^{2N}} dx f(x) \mathbf{E}^{Q}_{q(x)}\bigbkt{\exp\bkt{\sum_{(\mathfrak{i},\mathfrak{i}')\in \mathcal{E}_{N}}\int_{0}^{t} \beta_{\varepsilon}^{2}R_{\varepsilon}(A_{\mathfrak{i},\mathfrak{i}'}Q^{\perp}(s)) ds } f'(x(Q(t))) }.
\end{align}
Define the relative-coordinate operator
\begin{align*}
    \mathscr{R}_{N;\varepsilon}^{\theta} := \frac{1}{2} \Delta_{q^{\perp}} + \beta_{\varepsilon}^{2} \sum_{(\mathfrak{i},\mathfrak{i}')\in \mathcal{E}_{N}} R_{\varepsilon}(A_{\mathfrak{i},\mathfrak{i}'}(q^{\perp})).
\end{align*}
Applying the change of variable $q(x) \mapsto q$ together with (iii) in Lemma \ref{lemma:jacobi} yields
\begin{align} \label{FK_1}
    \text{(R.H.S.) of \eqref{FK_0}}
    &= \int_{\mathbb{R}^{2N}} dq \widetilde{f}(q) \mathbf{E}^{Q}_{q}\bigbkt{\exp\bkt{\sum_{(\mathfrak{i},\mathfrak{i}')\in \mathcal{E}_{N}} \int_{0}^{t} \beta_{\varepsilon}^{2}R_{\varepsilon}(A_{\mathfrak{i},\mathfrak{i}'}Q^{\perp}(s)) ds } \widetilde{f'}(Q(t)) } \nonumber\\
    &=\int_{\mathbb{R}^{4N}} dq dq'\widetilde{f}(q) G_{t}(q_{0}-q'_{0}) \exp(t \mathscr{R}_{N;\varepsilon}^{\theta})(q^{\perp},q'^{\perp})\widetilde{f'}(q'), \quad \widetilde{f}(q) := \mathcal{U}f(q).
\end{align}
To compute the limit of $\exp(t \mathscr{R}_{N;\varepsilon}^{\theta})$, we  move $G_{t}(q_{0}-q'_{0})$ to the other side. To this end, we combine \eqref{FK_0} and \eqref{FK_1} and extend
the resulting identity from product test functions to general
test functions. Thus, for each $F\in C_{c}(\mathbb{R}^{4N})$, we have
\begin{align} \label{equation:Jacobi_representation_varepsilon}
    &\int_{\mathbb{R}^{4N}} dxdx' F(x,x') \exp(t \mathscr{L}_{N;\varepsilon}^{\theta})(x,x') \nonumber\\
    &= \int_{\mathbb{R}^{4N}} dq dq'\widetilde{F}(q,q') G_{t}(q_{0}-q'_{0}) \exp(t \mathscr{R}_{N;\varepsilon}^{\theta})(q^{\perp},q'^{\perp}), \quad \widetilde{F}(q,q') := F(x(q),x(q')).
\end{align}
This together with the pointwise convergence and the pointwise bound in \cite[Theorem 2.3]{chen}  then yields that for each $\widetilde{F} \in C_{c}(\mathbb{R}^{4N})$, 
\begin{align} \label{FK_2}
    \lim_{\varepsilon \to 0}\int_{\mathbb{R}^{4N}} dqdq'\widetilde{F}(q,q') \exp(t \mathscr{R}_{N;\varepsilon}^{\theta})(q^{\perp},q'^{\perp}) 
    = \int_{\mathbb{R}^{4N}} dqdq' \widetilde{F}(q,q') \exp(t \mathscr{L}_{N}^{\theta})(x(q),x'(q')) G_{t}(q_{0} - q'_{0})^{-1},
\end{align}
which defines a positive distribution. Thus, the Riesz-Markov representation theorem \cite[Theorem 6.22]{lieb2001analysis} yields that \eqref{FK_2} defines a Borel measure in $\mathcal{M}_{+}(\mathbb{R}^{4N})$. The right-hand side of \eqref{FK_2} shows that this measure has a density function with respect to $dqdq'$. The left-hand side of \eqref{FK_2} implies that the density function only depends on $(q^{\perp},q'^{\perp})$. We denote this density function by $\exp(t \mathscr{R}_{N}^{\theta})(q^{\perp},q'^{\perp})$. Thus, \eqref{FK_2} implies that
for each $f,f' \in C_{c}(\mathbb{R}^{2N})$,
\begin{align*}
    \int_{\mathbb{R}^{4\times N}} dxdx' f(x) \exp(t \mathscr{L}_{N}^{\theta})(x,x') f'(x')
    = \int_{\mathbb{R}^{4N}} dq dq' \widetilde{f}(q) G_{t}(q_{0}-q'_{0}) \exp(t \mathscr{R}_{N}^{\theta})(q^{\perp},q'^{\perp}) \widetilde{f'}(q'),
\end{align*}
which implies \eqref{equation:decomposition}. The proof is complete.     
\end{proof}

\section{Ground-state energy formula and mean-field approximation} \label{section:Ground_state_energy_and_mean_field_approximation}

In this section, we develop in detail the ground-state energy formula for the bosonic Hamiltonian stated in Theorem \ref{theorem:informal_formula_ground_state_energy}.

\subsection{Resolvent formula and the principal operator} \label{section:Tsai_setting}
Fix $N\geq 3$. For each $(\mathfrak{i},\mathfrak{i}') \in \mathcal{E}_{N}$, define the linear transformation
\begin{align*}
    \mathcal{S}_{\mathfrak{i},\mathfrak{i}'}: \mathbb{R}^{2N-2} \mapsto \mathbb{R}^{2N}
\end{align*}
by
\begin{align*}
    \mathcal{S}_{\mathfrak{i},\mathfrak{i}'}(x_{2},\ldots,x_{N})
    := (x_{3},\ldots, \underbrace{x_{2}}_{\mathfrak{i}\text{-th}},\ldots,\underbrace{x_{2}}_{\mathfrak{i}'\text{-th}},\ldots,x_{N}) \in \mathbb{R}^{2N}, \quad x_{k} \in \mathbb{R}^{2},
\end{align*}
where $x_{2}$ is placed in the $\mathfrak{i}$-th and $\mathfrak{i}'$-th components of $\mathbb{R}^{2N}$. This transformation defines  the  unbounded operator 
\begin{align*}
    \mathcal{S}_{\mathfrak{i},\mathfrak{i}'}:L^{2}(\mathbb{R}^{2N}) \mapsto L^{2}(\mathbb{R}^{2N-2})
\end{align*}
by
\begin{align} \label{definition:lowering_operator}
    \mathcal{S}_{\mathfrak{i},\mathfrak{i}'}\psi
    := \psi \circ \mathcal{S}_{\mathfrak{i},\mathfrak{i}'}.
\end{align}
Intuitively, this operator restricts a function to the collision configuration, where the variable $x_{2}$ above represents the collision position. For every $q \in \mathbb{R}^{2N}$ and $2\leq k\leq N$, denote
\begin{align*}
    \quad q_{k-N} := (q_{k},\ldots,q_{N}) \in \mathbb{R}^{2(N-k+1)}.
\end{align*}
Using the following convention for the Fourier transform:
\begin{align*}
    \widehat{f}(p) := \frac{1}{(2\pi)^{d/2}} \int_{\mathbb{R}^{d}} dx \exp(-i x\cdot p) f(x), \quad p\in \mathbb{R}^{d},
\end{align*}
we have
\begin{align} \label{definition:adjoint_lowering_operator_-1}
    \widehat{\mathcal{S}_{\mathfrak{i},\mathfrak{i}'}\psi}(q_{2-N}) = \frac{1}{2\pi}
    \int_{\mathbb{R}^{2}} \widehat{\psi}(\mathcal{M}_{\mathfrak{i},\mathfrak{i}'}(q)) dq_{1}, 
\end{align}
where the linear transform $\mathcal{M}_{\mathfrak{i},\mathfrak{i}'}$ is given by
\begin{align*}
    \mathcal{M}_{\mathfrak{i},\mathfrak{i}'}(q):= \bkt{q_{3},\ldots,\underbrace{\frac{1}{2}q_{2}+q_{1}}_{\mathfrak{i}\text{-th}},\ldots,\underbrace{\frac{1}{2}q_{2}-q_{1}}_{\mathfrak{i}'\text{-th}},\ldots,q_{N}} \in \mathbb{R}^{2N}, \quad q:= (q_{1},\ldots,q_{N}) \in \mathbb{R}^{2N}.
\end{align*}
Similarly, the above definition means that $q_{2}/2+q_{1}$ and $q_{2}/2-q_{1}$ are placed in the $\mathfrak{i}$-th and $\mathfrak{i}'$-th components of $\mathbb{R}^{2N}$, respectively. The Fourier representation
\eqref{definition:adjoint_lowering_operator_-1} provides a natural
domain for $\mathcal{S}_{\mathfrak{i},\mathfrak{i}'}$, given by
\begin{align*}
    D(\mathcal{S}_{\mathfrak{i},\mathfrak{i}'})
    :=
    \set{
    \psi\in L^{2}(\mathbb{R}^{2N}):
    \widehat{\mathcal{S}_{\mathfrak{i},\mathfrak{i}'}\psi}
    \in L^{2}(\mathbb{R}^{2N-2})
    },
\end{align*}
where the Fourier transform on the right-hand side is defined by
the integral in \eqref{definition:adjoint_lowering_operator_-1}.

Define the operator 
\begin{align*}
    \mathcal{S}^{*}_{\mathfrak{i},\mathfrak{i}'}&: L^{2}(\mathbb{R}^{2N-2}) \mapsto \bigcap_{a>1} H^{-a}(\mathbb{R}^{2N}), 
\end{align*}
by
\begin{align} \label{definition:adjoint_lowering_operator}
    \widehat{\mathcal{S}^{*}_{\mathfrak{i},\mathfrak{i}'}\varphi}(q) &:= \frac{1}{2\pi} \widehat{\varphi}(q_{\mathfrak{i}}+q_{\mathfrak{i}'},q_{\overline{\mathfrak{i}\mathfrak{i}'}}), \quad q = (q_{1},\ldots,q_{N}) \in \mathbb{R}^{2N},  
\end{align}
where $q_{\overline{\mathfrak{i}\mathfrak{i}'}} \in \mathbb{R}^{2N-4}$ is obtained by removing $q_{\mathfrak{i}}$ and $q_{\mathfrak{i}'}$ from $q$. Here $H^{a}(\mathbb{R}^{2N})$ denotes the Sobolev space of degree $a \in \mathbb{R}$. The operator $\mathcal{S}^{*}_{\mathfrak{i},\mathfrak{i}'}$ is the
dual operator of $\mathcal{S}_{\mathfrak{i},\mathfrak{i}'}$ with respect to Sobolev duality, in the sense that, for every $a>1$, $\varphi\in L^{2}(\mathbb{R}^{2N-2})$, and $\psi\in H^{a}(\mathbb{R}^{2N})$,
\begin{align} \label{equation:adjoint}
    \langle \varphi,\mathcal{S}_{\mathfrak{i},\mathfrak{i}'}\psi \rangle_{L^{2}(\mathbb{R}^{2N-2})} = \langle \mathcal{S}^{*}_{\mathfrak{i},\mathfrak{i}'}\varphi,\psi \rangle_{L^{2}(\mathbb{R}^{2N})}.
\end{align}
Denote the free resolvent operator by
\begin{align*}
    \mathcal{G}_{z} := \bkt{-\frac{1}{2}\Delta_{x} - z I_{\mathbb{R}^{2N}}}^{-1}, \quad x = (x_{1},\ldots,x_{N}) \in \mathbb{R}^{2N}, \quad z \in \mathbb{C}\setminus \mathbb{R}_{\geq 0}.
\end{align*}
For each $z \in \mathbb{C}_{<0}$,
define the operators
\begin{align} \label{definition:domain_adjusting}
    \mathcal{S}_{z} : L^{2}_{\textup{sym}}(\mathbb{R}^{2N}) \mapsto L^{2}_{\textup{sym}'}(\mathbb{R}^{2N-2}), \quad \widetilde{\mathcal{S}}_{z} : L^{2}_{\textup{sym}'}(\mathbb{R}^{2N-2}) \mapsto L^{2}_{\textup{sym}}(\mathbb{R}^{2N}),
\end{align}
by
\begin{align*}
    \mathcal{S}_{z} &:= \sum_{(\mathfrak{i},\mathfrak{i}')\in \mathcal{E}_{N}} \mathcal{S}_{\mathfrak{i},\mathfrak{i}'}\mathcal{G}_{z}, \quad \widetilde{\mathcal{S}}_{z} := \sum_{(\mathfrak{i},\mathfrak{i}')\in \mathcal{E}_{N}} \mathcal{G}_{z} \mathcal{S}^{*}_{\mathfrak{i},\mathfrak{i}'}. 
\end{align*}
By Lemma \ref{lemma:estimate_from_GQT}, both operators in \eqref{definition:domain_adjusting} are bounded.

Consider the subspace $L^{2}_{\textup{sym}'}(\mathbb{R}^{2N-2})$ of $L^{2}(\mathbb{R}^{2N-2})$ given by
\begin{align*}
    L^{2}_{\textup{sym}'}(\mathbb{R}^{2N-2})
    := \set{\psi\in L^{2}(\mathbb{R}^{2N-2}): \psi(y_{2},y_{3},\ldots,y_{N}) = \psi(y_{2},y_{\sigma_{3}},\ldots,y_{\sigma_{N}}) \quad \text{(a.e.)} \quad \forall \sigma},
\end{align*}
where $\sigma$ denotes a permutation of the set $\{3,\ldots,N\}$. For each $z \in \mathbb{C}_{<0}$,  define the  Fourier multiplier $\mathcal{J}_{z}$ on $L^{2}_{\textup{sym}'}(\mathbb{R}^{2N-2})$ by
\begin{align*}
    \widehat{\mathcal{J}_{z}\varphi}(q_{2-N}) := \log\bkt{\frac{1}{2}|q|_{2-N}^{2}-z }\widehat{\varphi}(q_{2-N}), \quad 
    |q|_{2-N}^{2} := \frac{1}{2}|q_{2}|^{2} +|q_{3}|^{2} + \ldots+|q_{N}|^{2},
\end{align*}
with $z$-independent domain 
\begin{align} \label{definition:domain_J}
    D(\mathcal{J}_{z}) := D(\mathcal{J}_{-1})
\end{align}
where
\begin{align} \label{definition:domain_J_0}
    D(\mathcal{J}_{-1})  := \set{\varphi \in L_{\textup{sym}'}^{2}(\mathbb{R}^{2N-2}): \int_{\mathbb{R}^{2N-2}} |\widehat{\varphi}(q_{2-N}) \log (|q|_{2-N}^{2}+1)|^{2} dq_{2-N} < \infty}.
\end{align}
Here $\log(z)$ denotes the principal branch of the complex logarithm. One easily verifies that for each $z \in \mathbb{C}_{<0}$, the natural domain of $\mathcal{J}_{z}$  coincides with the set \eqref{definition:domain_J_0}. With these notations, define the unbounded operator $\mathcal{A}^{N;\theta}_{z}$ on $L^{2}_{\textup{sym}'}(\mathbb{R}^{2N-2})$ by 
\begin{align} \label{equation:resolvent_representation_1}
    \mathcal{A}^{N;\theta}_{z} &:= \frac{1}{4\pi} (\mathcal{J}_{z} - \beta_{\star}(\theta)I') - \frac{2}{N(N-1)} \sum_{(\mathfrak{i},\mathfrak{i}')\neq (\mathfrak{j},\mathfrak{j}')} \mathcal{S}_{\mathfrak{i},\mathfrak{i}'}\mathcal{G}_{z} \mathcal{S}_{\mathfrak{j},\mathfrak{j}'}^{*}
\end{align}
with domain
\begin{align} \label{definition:domain_A}
    D(\mathcal{A}_{z}^{N;\theta}) = D(\mathcal{J}_{-1}),
\end{align}
where $\beta_{\star}(\theta)$ is given in \eqref{definition:beta_star} and $I'$ denotes the identity operator on $L^{2}_{\textup{sym}'}(\mathbb{R}^{2N-2})$. Following Rajeev \cite{rajeev1999condensation}, we refer to $\mathcal{A}_{z}^{N;\theta}$  as the \textit{principal operator} at energy  $z$.  By Lemma \ref{lemma:estimate_from_GQT},
\begin{align} \label{definition:domain_adjusting_1}
    \sum_{(\mathfrak{i},\mathfrak{i}')\neq (\mathfrak{j},\mathfrak{j}')} \mathcal{S}_{\mathfrak{i},\mathfrak{i}'}\mathcal{G}_{z} \mathcal{S}_{\mathfrak{j},\mathfrak{j}'}^{*}: L^{2}_{\textup{sym}'}(\mathbb{R}^{2N-2}) \mapsto L^{2}_{\textup{sym}'}(\mathbb{R}^{2N-2})
\end{align}
is a bounded operator. Consequently, the domain of $\mathcal{A}_{z}^{N;\theta}$ coincides with the domain of $\mathcal{J}_{z}$.

For an operator $\mathcal{T}$, $\rho(\mathcal{T})$ denotes its resolvent set. For a Hilbert space $H$, $B(H)$ denotes the space of bounded linear operators on $H$.  With these notations, combining  \cite[Proposition 1.4 \& Theorem 1.6 \& Remark 1.8]{Gu_2d} gives the following results.
\begin{theorem} [Gu, Quastel, and Tsai \cite{Gu_2d}]
\label{theorem:resolvent_formula}
Assume that $\theta \in \mathbb{R}$. Then there exists $C_{\ref{theorem:resolvent_formula}} > 0$ such that for every $N\geq 3$, one has 
\begin{align*}
    \widetilde{\mathscr{O}}_{N}^{\theta} \subseteq \rho(\mathscr{H}_{N}^{\theta;\textup{sym}}), 
\end{align*}
where 
\begin{align} \label{definition:half_plane}
    \widetilde{\mathscr{O}}_{N}^{\theta} := \{z \in \mathbb{C}:\Re(z) < -\exp(C_{\ref{theorem:resolvent_formula}}N^{2}+\beta_{\star}(\theta))\}.
\end{align}
In particular, it holds that
\begin{align} \label{equation:resolvent_representation}
    (\mathscr{H}_{N}^{\theta;\textup{sym}} - z)^{-1} = \mathcal{G}_{z} + \frac{2}{N(N-1)} \widetilde{\mathcal{S}}_{z} (\mathcal{A}_{z}^{N;\theta})^{-1}\mathcal{S}_{z} \quad \forall z \in \widetilde{\mathscr{O}}_{N}^{\theta}.
\end{align}
\end{theorem}

\subsection{Analytic continuation and the ground-state energy formula} \label{section:AC_GSE_formula}

\begin{figure}
\centering
\begin{tikzpicture}[
    scale=0.44,
    shorten > = 1pt,
    auto,
    node distance=3cm,
    decoration={
        snake,
        pre length=3pt,
        post length=7pt
    }
]


\filldraw[lightgray,opacity=0.7] (-12,-4) rectangle (-15,4);

\draw[-] (-15,0) -- (-4,0);
\draw (-3,-0.25) node[anchor=south]{\footnotesize $\ldots$};
\draw[->] (-2,0) -- (0,0);

\draw (0.35,0) node[anchor=south]{\tiny $\Im(z)=0$};

\draw[->] (-1,-4) -- (-1,4) node[anchor=south]{\tiny $\Re(z)=0$};

\draw[blue,very thick] (-5,-4) -- (-5,4);
\draw (-4.5,4) node[anchor=south]{\tiny $-\exp(C'N)$};

\draw[-] (-7,-4) node[anchor=north]{\tiny $\inf \mathscr{H}_{N}^{\theta;\mathrm{sym}}$} -- (-7,4);

\draw[blue,very thick] (-9,-4) -- (-9,4);
\draw (-9.5,4) node[anchor=south]{\tiny $-\exp(CN\log N)$};

\draw[dashed] (-12,-4) node[anchor=north]{\tiny $-\exp(C''N^{2})$} -- (-12,4);


\draw (4.5,0) node[anchor=east]{\large $\xLongrightarrow{\quad}$};


\draw[-] (5,0) -- (16,0);
\draw (17,-0.25) node[anchor=south]{\footnotesize $\ldots$};
\draw[->] (18,0) -- (20,0);

\draw (20.35,0) node[anchor=south]{\tiny $\Im(z)=0$};

\draw[->] (19,-4) -- (19,4) node[anchor=south]{\tiny $\Re(z)=0$};

\draw[blue,very thick] (15,-4) -- (15,4);

\draw[-] (13,-4) node[anchor=north]{\tiny $\inf \mathscr{H}_{N}^{\theta;\mathrm{sym}}$} -- (13,4);

\draw[blue,very thick] (11,-4) -- (11,4);

\draw[red,very thick] (5,0) -- (13,0);

\filldraw[lightgray,opacity=0.7] (8,-4) rectangle (5,4);

\draw[dashed] (8,-4) node[anchor=north]{\tiny $-\exp(C''N^{2})$} -- (8,4);

\filldraw[color=black, fill=black] (10,0) circle (2.5pt);
\draw (10,0) node[anchor=south]{\tiny $z$};

\fill[lightgray,opacity=0.7]
    (13,0)
    to[out=95,in=-10] (8,3.5)
    -- (8,-3.5)
    to[out=10,in=-95] (13,0)
    -- cycle;

\draw[dashed] (13,0) to[out=95,in=-10] (8,3.5);
\draw[dashed] (13,0) to[out=-95,in=10] (8,-3.5);

\filldraw[color=black, fill=white] (13,0) circle (2.5pt);

\draw (9.5,3.5) node[anchor=south]{\tiny $\mathscr{O}_{N}^{\theta}$};


\end{tikzpicture}
\begin{minipage}{0.9\textwidth}
\caption{ \justifying \footnotesize
The left and right figures denote the domain $\mathscr{O}_{N}^{\theta}$ of $\Psi_{N}^{\theta}(z)$ before and after analytic continuation, respectively. The blue vertical lines indicate the two-sided bounds on the ground-state energy obtained in Theorem \ref{theorem:2d_DBG}. The red horizontal line indicates the portion of the real axis on which the extended resolvent formula \eqref{equation:resolvent_representation_2} is valid, allowing us to take $z \uparrow  \inf \mathscr{H}_{N}^{\theta;\textup{sym}}$.
}\label{Figure:1}
\end{minipage}
\end{figure}
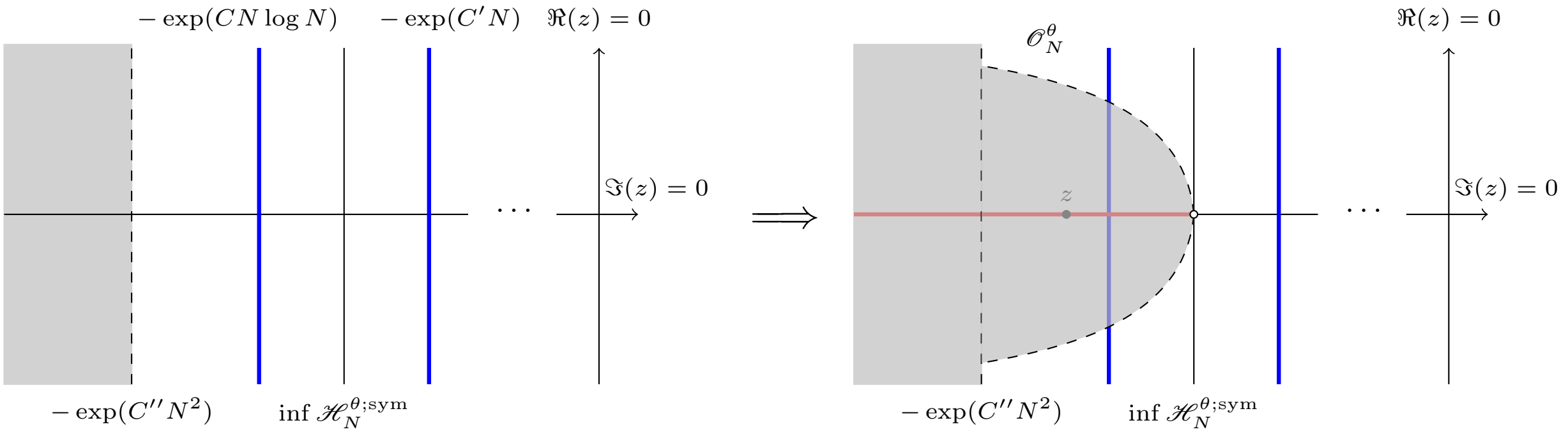

We now introduce the main ideas behind the proof of the formula \eqref{equation:informal_formula_ground_state_energy}. The \textit{key step} is to extend Theorem \ref{theorem:resolvent_formula} by establishing an analytic continuation of 
\begin{align} \label{definition:Psi_0}
    \Psi_{N}^{\theta}(z) := \mathcal{G}_{z} + \frac{2}{N(N-1)} \widetilde{\mathcal{S}}_{z} (\mathcal{A}_{z}^{N;\theta})^{-1}\mathcal{S}_{z}.
\end{align}
Throughout this paper, for each Hilbert space $H$ and open subset $\mathscr{O} \subseteq \mathbb{C}$, we call $\Psi: \mathscr{O}  \mapsto B(H)$ a $B(H)$-valued holomorphic function if for all $\psi,\varphi \in H$, $z \in \mathscr{O} \mapsto \langle \psi,\Psi(z)\varphi \rangle_{H}$ is  holomorphic. Theorem \ref{theorem:resolvent_extension} extends the resolvent representation from the half-plane \eqref{definition:half_plane} to a domain containing the real line below the ground-state energy (see Figure \ref{Figure:1}).

\begin{theorem} \label{theorem:resolvent_extension}
Suppose that $\theta\in \mathbb{R}$  and $N\geq 3$. There exists a domain $\mathscr{O}_{N}^{\theta} \subseteq \mathbb{C}_{<0}$ such that 
\begin{align*}
    (-\infty,\inf \mathscr{H}_{N}^{\theta;\textup{sym}}) \subseteq \mathscr{O}_{N}^{\theta}.
\end{align*}
In addition, 
\begin{align*}
    \Psi_{N}^{\theta}:\mathscr{O}_{N}^{\theta} \mapsto B(L^{2}_{\textup{sym}}(\mathbb{R}^{2N}))
\end{align*}    
defines a holomorphic function, and 
\begin{align} \label{equation:resolvent_representation_2}
    (\mathscr{H}_{N}^{\theta;\textup{sym}} - z)^{-1} = \Psi_{N}^{\theta}(z)  \quad \forall z \in \mathscr{O}_{N}^{\theta}.
\end{align}
Here the domain $\mathscr{O}_{N}^{\theta}$ is given in \eqref{definition:set_O}.
\end{theorem}

We prove Theorem \ref{theorem:resolvent_extension} in Section \ref{section:proof_analytic_continuation}. 
To do this, since the free resolvent and the operators $\widetilde{\mathcal{S}}_{z}$ and $\mathcal{S}_{z}$ are holomorphic on $\mathbb{C}_{<0}$ (see Lemma \ref{lemma:S_is_holomorphic}), the main task is to analytically continue $(\mathcal{A}_{z}^{N;\theta})^{-1}$. The key observation is that, although $\mathcal{A}_{z}^{N;\theta}$
is unbounded because of the Fourier multiplier  $\mathcal{J}_{z}$ on the right-hand side of \eqref{recall_principal_operator}, the difference
$\mathcal{A}_{z}^{N;\theta}-\mathcal{A}_{w}^{N;\theta}$ extends to a bounded operator.
The identity \eqref{equation:resolvent_identity_principal_operator} for this difference is a key ingredient in establishing
the boundedness of $(\mathcal{A}_{z}^{N;\theta})^{-1}$ beyond the half-plane
in \eqref{definition:half_plane}. It also yields a differential equation for $(\mathcal{A}_{z}^{N;\theta})^{-1}$, from which the holomorphicity of $(\mathcal{A}_{z}^{N;\theta})^{-1}$ follows; see Lemma \ref{lemma:Analytic_continuation}.

\begin{lemma}
\label{proposition:Hilbert_identity}
Assume that $\theta \in \mathbb{R}$ and $N\geq 3$.  For each $z,w \in \mathbb{C}_{<0}$, the difference
$\mathcal{A}_{z}^{N;\theta} - \mathcal{A}_{w}^{N;\theta}$ extends to a bounded operator on $L^{2}_{\textup{sym}'}(\mathbb{R}^{2N-2})$, and its extension satisfies
\begin{align} \label{equation:resolvent_identity_principal_operator}
    \mathcal{A}_{z}^{N;\theta} - \mathcal{A}_{w}^{N;\theta} = - \frac{2}{N(N-1)}(z-w) \mathcal{S}_{z} \widetilde{\mathcal{S}}_{w} \quad \text{in} \quad B(L^{2}_{\textup{sym}'}(\mathbb{R}^{2N-2})).
\end{align}
In particular, for any smooth curve $\gamma:[0,1]\to\mathbb{C}_{<0}$ with $\gamma(0)=w$ and $\gamma(1)=z$,
we have
\begin{align} \label{equation:integral_principal_operator}
    \mathcal{A}_{z}^{N;\theta} - \mathcal{A}_{w}^{N;\theta} = -\frac{2}{N(N-1)} \int_{0}^{1} \mathcal{S}_{\gamma(s)} \widetilde{\mathcal{S}}_{\gamma(s)} \gamma'(s)  ds \quad \text{in} \quad B(L^{2}_{\textup{sym}'}(\mathbb{R}^{2N-2})),
\end{align}
where the integral is understood in the Bochner sense.
\end{lemma}

The extended resolvent formula \eqref{equation:resolvent_representation_2} implies the ground-state energy formula below. We prove Proposition \ref{theorem:formula_ground_state_energy} in Section \ref{section:proof_formula_ground_state_energy}.
\begin{proposition} \label{theorem:formula_ground_state_energy}
Assume that $\theta \in \mathbb{R}$. For each $N\geq 3$, it holds that
\begin{align} \label{equation:formula_ground_state_energy}
    \inf \mathscr{H}_{N}^{\theta;\textup{sym}} = -\exp(-4\pi E_{N}^{\theta}(-1)).
\end{align}
Here for each $z < 0$, $\mathcal{A}_{z}^{N;\theta}$ is self-adjoint, and we set
\begin{align} \label{defnition:bottom_spectrum_A}
    E_{N}^{\theta}(z) := \inf_{\psi \in D(\mathcal{A}_{z}^{N;\theta}), \; ||\psi||_{L^{2}_{\textup{sym}'}(\mathbb{R}^{2N-2})} = 1} \langle \psi,\mathcal{A}_{z}^{N;\theta}\psi \rangle_{L^{2}_{\textup{sym}'}(\mathbb{R}^{2N-2})}.
\end{align}
\end{proposition}

A key step in the proof is to show that 
\begin{align*}
    E_{N}^{\theta}(z) > 0 \quad \forall z \in (-\infty,\inf \mathscr{H}_{N}^{\theta;\textup{sym}}).
\end{align*}
Thus, by the continuity of $E_{N}^{\theta}(z)$, taking $z \uparrow \inf \mathscr{H}_{N}^{\theta;\textup{sym}}$ yields that $E_{N}^{\theta}(\inf \mathscr{H}_{N}^{\theta;\textup{sym}}) \geq 0$. If this quantity were strictly positive, $\Psi_{N}^{\theta}(z)$
would remain uniformly bounded in operator norm as $z$ approaches
$\inf\mathscr{H}_{N}^{\theta;\textup{sym}}$ from below.
However, by the extended resolvent formula
\eqref{equation:resolvent_representation_2},
$\Psi_{N}^{\theta}(z)$ agrees with
$(\mathscr{H}_{N}^{\theta;\textup{sym}}-z)^{-1}$
for every $z\in(-\infty,\inf\mathscr{H}_{N}^{\theta;\textup{sym}})$.
This would contradict the divergence of the resolvent norm as
$z\uparrow\inf\mathscr{H}_{N}^{\theta;\textup{sym}}$.
Therefore,
\begin{align} \label{equation:equation_for_ground_state_energy}
    E_{N}^{\theta}
    (\inf\mathscr{H}_{N}^{\theta;\textup{sym}})=0.
\end{align}
Applying a scaling argument to
\eqref{equation:equation_for_ground_state_energy} then yields
the ground-state energy formula \eqref{equation:formula_ground_state_energy}.

\subsection{Rajeev's mean-field approximation} \label{section:mean_field_approximation}

As an application, let us apply Rajeev's idea to the ground-state energy formula \eqref{equation:formula_ground_state_energy} and explain how this approximation relates to the coefficient $C_{4-\textup{particles}}$ in \eqref{conjecture:Rajeev1999}.

Rajeev conjectures that $E_{N}^{\theta}(-1)$ can be approximated by minimizing over the class  $\mathscr{A}_{\textup{M-F}}(\mathbb{R}^{2N-2})$, a variational ansatz corresponding to mean-field theory, given by 
\begin{align} \label{definition:mean_field_pair}
    \mathscr{A}_{\textup{M-F}}(\mathbb{R}^{2N-2}):= &\set{ v\star u :=  v(x_{2}) u^{\otimes (N-2)}(x_{3-N}): ||v||_{L^{2}(\mathbb{R}^{2})} = ||u||_{L^{2}(\mathbb{R}^{2})} = 1, \nonumber \\ 
    &\int_{\mathbb{R}^{2}} dp |\widehat{u}(p) \log (|p|^{2}+1) |^{2} \bigvee \int_{\mathbb{R}^{2}} dp |\widehat{v}(p) \log (|p|^{2}+1) |^{2} < \infty}.
\end{align}
See \cite[(27)]{rajeev1999condensation} for details. It is expected that \eqref{equation:Raj_idea} holds. Here $E_{N,\theta}^{\textup{M-F}}$ is referred to as the \textit{mean-field exponent}, given by
\begin{align} \label{equation:Raj_energy}
    E_{N,\theta}^{\textup{M-F}} := \inf_{v\star u \in \mathscr{A}_{\textup{M-F}}(\mathbb{R}^{2N-2})}  \langle v\star u,\mathcal{A}_{-1}^{N;\theta} v\star u \rangle_{L^{2}_{\textup{sym}'}(\mathbb{R}^{2N-2})}.
\end{align}
One easily verifies that
$\mathscr{A}_{\textup{M-F}}(\mathbb{R}^{2N-2})$
consists precisely of the functions in
$D(\mathcal{A}_{-1}^{N;\theta})$ of the form $v\star u$, with
$\|v\|_{L^{2}(\mathbb{R}^{2})}
=\|u\|_{L^{2}(\mathbb{R}^{2})}=1$. Our main result identifies the leading-order term of the mean-field exponent. We prove Theorem \ref{proposition:mean_field_two_side} in Section \ref{section:mean_field_exponent}.
\begin{theorem} \label{proposition:mean_field_two_side}
Assume that $\theta \in \mathbb{R}$. For all sufficiently large $N$, 
\begin{align} \label{equation:mean_field_two_side}
    -C_{\textup{G-N}}N-O(\sqrt{N}) \leq E_{N,\theta}^{\textup{M-F}} \leq -C_{\textup{G-N}}N+o(N).
\end{align}
\end{theorem}

The two possible configurations of distinct collision labels in the principal operator $\mathcal{A}_{-1}^{N;\theta}$ give rise to the three- and four-particle terms:
\begin{align} \label{three_and_four_contribution}
    \frac{2}{N(N-1)}\sum_{|\{\mathfrak{i},\mathfrak{i}'\}\cap \{\mathfrak{j},\mathfrak{j}'\}| = 1} \mathcal{S}_{\mathfrak{i},\mathfrak{i}'}\mathcal{G}_{-1} \mathcal{S}_{\mathfrak{j},\mathfrak{j}'}^{*}, \quad \frac{2}{N(N-1)}\sum_{|\{\mathfrak{i},\mathfrak{i}'\}\cap \{\mathfrak{j},\mathfrak{j}'\}| = 0} \mathcal{S}_{\mathfrak{i},\mathfrak{i}'}\mathcal{G}_{-1} \mathcal{S}_{\mathfrak{j},\mathfrak{j}'}^{*}, \quad \text{respectively}.
\end{align}
In the estimate \eqref{equation:mean_field_two_side}, these terms contribute $o(N)$ and $-C_{\textup{G-N}}N$ to
$E_{N,\theta}^{\textup{M-F}}$, respectively.

The heuristic behind these expected orders is as follows. A formal application of the law of large
numbers suggests that the resolvent $\mathcal{G}_{-1}$ contributes a factor of order $1/N$. Since there are $O(N^{3})$ three-particle configurations and the normalization is of order $1/N^{2}$, the resulting contribution is expected to be $O(N^{3}/N^{2} \cdot 1/N) = O(1)$. The same heuristic also explains why the four-particle contribution in \eqref{three_and_four_contribution} is expected to be of linear order in $N$.

A notable feature of \eqref{equation:mean_field_two_side} is that the leading-order coefficient of $E_{N,\theta}^{\textup{M-F}}$ is universal, in the sense that it is independent of $\theta$. This is because the principal operator $\mathcal{A}_{-1}^{N;\theta}$ depends on $\theta$ only through the term $-\beta_{\star}(\theta)/(4\pi) I'$. For fixed $\theta$, this term contributes only $O(1)$ to $E_{N,\theta}^{\textup{M-F}}$ as $N\to\infty$, and therefore does not affect the leading-order term.

\section{Analytic continuation of the resolvent formula}
\label{section:proof_analytic_continuation}

The goal of this section is to prove Theorem \ref{theorem:resolvent_extension}. To do this, the main step is to show that $\Psi_{N}^{\theta}$ admits a $B(L_{\textup{sym}}^{2}(\mathbb{R}^{2N}))$-valued holomorphic extension on $\mathscr{O}_{N}^{\theta}$, where
\begin{align} \label{definition:set_O}
    \mathscr{O}^{\theta}_{N} := \mathscr{G}^{\theta}_{N}\cup \mathscr{O}'^{\theta}_{N}, \quad 
    \mathscr{G}^{\theta}_{N} := \set{z \in \mathbb{C}: \Re(z) < \inf \mathscr{H}_{N}^{\theta;\textup{sym}}, \quad |\Im(z)| < \frac{(-\inf \mathscr{H}_{N}^{\theta;\textup{sym}})E_{N}^{\theta}(\Re(z)) }{2 C_{\eqref{definition:constant_set_O}} N^{2}} }.
\end{align}
Here $\mathscr{O}'^{\theta}_{N}$ is the open half-plane defined in \eqref{definition:Psi}, which is contained in the half-plane in \eqref{definition:half_plane}.
Since the resolvent
$z \mapsto (\mathscr{H}_{N}^{\theta;\textup{sym}}-z)^{-1}$
is a $B(L_{\textup{sym}}^{2}(\mathbb{R}^{2N}))$-valued holomorphic
function on $\mathscr{O}_{N}^{\theta}$, Theorem \ref{theorem:resolvent_formula} shows that 
$\Psi_{N}^{\theta}$ is consistent with the resolvent on
$\mathscr{O}'^{\theta}_{N}\subseteq\mathscr{O}_{N}^{\theta}$. Thus,
\eqref{equation:resolvent_representation_2} then follows from the identity
theorem for holomorphic functions.

\subsection{Preliminary operator estimates}

To begin with, we introduce several auxiliary results. Combining \cite[Lemma 4.2 and Lemma 5.1]{Gu_2d} gives the following estimates. 
\begin{lemma} \label{lemma:estimate_from_GQT} 
Assume that $\theta \in \mathbb{R}$ and $N\geq 3$. There exists an absolute constant $C>0$ such that for each $z \in \mathbb{C}_{<0}$ and $(\mathfrak{i},\mathfrak{i}'),(\mathfrak{j},\mathfrak{j}') \in \mathcal{E}_{N}$ such that $(\mathfrak{i},\mathfrak{i}') \neq (\mathfrak{j},\mathfrak{j}')$, 
\begin{align*}
    ||\mathcal{S}_{\mathfrak{i},\mathfrak{i}'}\mathcal{G}_{z}||_{B({L^{2}(\mathbb{R}^{2N})},{L^{2}(\mathbb{R}^{2N-2})})} &= ||\mathcal{G}_{\overline{z}}\mathcal{S}_{\mathfrak{i},\mathfrak{i}'}^{*}||_{B({L^{2}(\mathbb{R}^{2N-2})},{L^{2}(\mathbb{R}^{2N})})} \leq C \Re(-z)^{-1/2},\\
    ||\mathcal{S}_{\mathfrak{i},\mathfrak{i}'}\mathcal{G}_{z}\mathcal{S}^{*}_{\mathfrak{j},\mathfrak{j}'} ||_{B({L^{2}(\mathbb{R}^{2N-2})})} &\leq C.
\end{align*}
\end{lemma}

To prove that $\Psi_{N}^{\theta}(\cdot)$ is holomorphic, we now prove the holomorphicity of $\mathcal{S}_{\cdot}$ and $\widetilde{\mathcal{S}}_{\cdot}$.
\begin{lemma} \label{lemma:S_is_holomorphic}
Suppose that $N\geq 3$. Then the following maps are holomorphic:
\begin{align*}
    &\mathcal{S}_{\cdot}:z\in\mathbb{C}_{<0} \mapsto \mathcal{S}_{z} \in B(L^{2}_{\textup{sym}}(\mathbb{R}^{2N}),L^{2}_{\textup{sym}'}(\mathbb{R}^{2N-2})), \\
    &\widetilde{\mathcal{S}}_{\cdot}:z\in\mathbb{C}_{<0} \mapsto \widetilde{\mathcal{S}}_{z} \in B(L^{2}_{\textup{sym}'}(\mathbb{R}^{2N-2}),L^{2}_{\textup{sym}}(\mathbb{R}^{2N})).
\end{align*}
\end{lemma}
\begin{proof}
Fix $\varphi \in L^{2}_{\textup{sym}'}(\mathbb{R}^{2N-2})$ and $\psi \in L^{2}_{\textup{sym}}(\mathbb{R}^{2N})$. Applying Plancherel's theorem  gives
\begin{align*}
    &\langle \varphi,\mathcal{S}_{z} \psi \rangle_{L^{2}_{\textup{sym}'}(\mathbb{R}^{2N-2})} \\
    &= \sum_{(\mathfrak{i},\mathfrak{i}')\in \mathcal{E}_{N}} \frac{1}{2\pi} \int_{\mathbb{R}^{2N}} dq_{1}dq_{2-N} \overline{\widehat{\varphi}(q_{2-N}) } \widehat{\psi}(\mathcal{M}_{\mathfrak{i},\mathfrak{i}'}(q_{1},q_{2-N})) \frac{1}{|\mathcal{M}_{\mathfrak{i},\mathfrak{i}'}(q_{1},q_{2-N})|^{2}/2-z}.
\end{align*}
By the Cauchy–Schwarz inequality, we have
\begin{align*}
    &\int_{\mathbb{R}^{2N}} dq_{1}dq_{2-N} \abs{\overline{\widehat{\varphi}(q_{2-N}) } \widehat{\psi}(\mathcal{M}_{\mathfrak{i},\mathfrak{i}'}(q_{1},q_{2-N})) \frac{1}{|\mathcal{M}_{\mathfrak{i},\mathfrak{i}'}(q_{1},q_{2-N})|^{2}/2-z}}\\
    &\leq  \Norm{\widehat{\varphi} \frac{1}{|\mathcal{M}_{\mathfrak{i},\mathfrak{i}'}(q_{1},q_{2-N})|^{2}/2-z} }_{L^{2}(\mathbb{R}^{2N})}  ||\psi||_{L^{2}_{\textup{sym}}(\mathbb{R}^{2N})}\\
    &\leq C(z) ||\varphi||_{L^{2}_{\textup{sym'}}(\mathbb{R}^{2N-2})} ||\psi||_{L^{2}_{\textup{sym}}(\mathbb{R}^{2N})}, \quad C(z) :=  \norm{(|q_{1}|^{2}-\Re(z))^{-2}}_{L^{1}(\mathbb{R}^{2})}^{1/2} < \infty.
\end{align*}
Indeed, since
\begin{align*}
    ||\mathcal{M}_{\mathfrak{i},\mathfrak{i}'}(q_{1},q_{2-N})|^{2}/2-z| \geq \frac{1}{2}\abs{\frac{1}{2}q_{2}-q_{1}}^{2}+\frac{1}{2}\abs{\frac{1}{2}q_{2}+q_{1}}^{2}-\Re(z)\geq |q_{1}|^{2}-\Re(z),
\end{align*}
it follows that $(|q_{1}|^{2}-\Re(z))^{-2} \in L^{1}(\mathbb{R}^{2})$.
In particular, for any smooth closed curve $\Gamma$ in $\mathbb{C}_{<0}$, 
\begin{align*}
    \sup_{z \in \Gamma} C(z) = C\left ( \sup_{z \in \Gamma} \Re(z) \right) <\infty.
\end{align*}
Thus,  applying Fubini's theorem gives
\begin{align*}
    \int_{\Gamma} dz \langle \varphi,\mathcal{S}_{z} \psi \rangle_{L^{2}_{\textup{sym}'}(\mathbb{R}^{2N-2})} = 0.
\end{align*}
Therefore, by Morera's theorem, $\mathcal{S}_{\cdot}$ is holomorphic. The holomorphicity  of $\widetilde{\mathcal{S}_{\cdot}}$ is proved in a similar way. The proof is complete. 

\end{proof}

Now consider 
\begin{align} \label{definition:Psi_3}
    \Psi_{N}^{\theta}:\mathscr{O}'^{\theta}_{N} \mapsto \text{(R.H.S.) of } \eqref{definition:Psi_0}
\end{align}
restricted to the half-plane
\begin{align} \label{definition:Psi}
    \mathscr{O}'^{\theta}_{N}:=\{z \in \mathbb{C}: \Re(z) < -\exp(C_{\eqref{definition:Psi}}N^{2}+\beta_{\star}(\theta))\},
\end{align}
and recall that 
\begin{align} \label{recall_principal_operator}
    \mathcal{A}_{z}^{N;\theta}= \frac{1}{4\pi} (\mathcal{J}_{z} - \beta_{\star}(\theta)I') - \frac{2}{N(N-1)} \sum_{(\mathfrak{i},\mathfrak{i}')\neq (\mathfrak{j},\mathfrak{j}')} \mathcal{S}_{\mathfrak{i},\mathfrak{i}'}\mathcal{G}_{z} \mathcal{S}_{\mathfrak{j},\mathfrak{j}'}^{*}.
\end{align}
Here $C_{\eqref{definition:Psi}} := C_{\eqref{theorem:resolvent_formula}}\vee C_{\eqref{equation:A_operator_norm}}$, where $C_{\eqref{equation:A_operator_norm}}$ is given as follows.
Throughout this section, $I$ and $I'$ denote the identity operators on $L^{2}_{\textup{sym}}(\mathbb{R}^{2N})$ and $L^{2}_{\textup{sym}'}(\mathbb{R}^{2N-2})$, respectively.  
\begin{lemma} \label{lemma:A_operator_norm}
Assume that $\theta \in \mathbb{R}$. There exists $C_{\eqref{equation:A_operator_norm}} = C_{\eqref{equation:A_operator_norm}}(\theta)$ such that for all $N\geq 3$,
\begin{align*}
    (\mathcal{A}_{z}^{N;\theta})^{-1} \in B(L^{2}_{\textup{sym}'}(\mathbb{R}^{2N-2})) \quad \forall z \in \mathscr{O}''^{\theta}_{N}
\end{align*}
and
\begin{align} \label{equation:A_operator_norm}
    ||(\mathcal{A}_{z}^{N;\theta})^{-1}||_{B(L^{2}_{\textup{sym}'}(\mathbb{R}^{2N-2}))} \leq \frac{C_{\eqref{equation:A_operator_norm}}}{\log (-\Re(z)) - \beta_{\star}(\theta)},
\end{align}
where
\begin{align*}
    \mathscr{O}''^{\theta}_{N}:=\{z \in \mathbb{C}: \Re(z) < -\exp(C_{\eqref{equation:A_operator_norm}}N^{2}+\beta_{\star}(\theta))  \}.
\end{align*}
In particular, the map $z\in \mathscr{O}''^{\theta}_{N} \mapsto (\mathcal{A}_{z}^{N;\theta})^{-1}$ is holomorphic.
\end{lemma}
\begin{proof}
Since
\begin{align} \label{estimate_J}
    \Norm{\bkt{\frac{1}{4\pi} (\mathcal{J}_{z} - \beta_{\star}(\theta)I')}^{-1}}_{B(L^{2}_{\textup{sym}'}(\mathbb{R}^{2N-2}))} \leq \frac{C}{\log (-\Re(z)) - \beta_{\star}(\theta)} \quad \forall\Re(z) < -\exp(\beta_{\star}(\theta)),
\end{align}
applying Lemma \ref{lemma:estimate_from_GQT} gives
\begin{align} \label{equattion:G_series_ratio}
    \Norm{\bkt{\frac{1}{4\pi} (\mathcal{J}_{z} - \beta_{\star}(\theta)I')}^{-1}\frac{2}{N(N-1)} \sum_{(\mathfrak{i},\mathfrak{i}')\neq (\mathfrak{j},\mathfrak{j}')} \mathcal{S}_{\mathfrak{i},\mathfrak{i}'}\mathcal{G}_{z} \mathcal{S}_{\mathfrak{j},\mathfrak{j}'}^{*}}_{B(L^{2}_{\textup{sym}'}(\mathbb{R}^{2N-2}))}\leq \frac{C N^{2}}{\log (-\Re(z)) - \beta_{\star}(\theta)},
\end{align}
where $C > 0$ is an absolute constant. By choosing $\Re(z)$ sufficiently negative such that the right-hand side of \eqref{equattion:G_series_ratio} is less than $1/2$, $\mathcal{A}_{z}^{N;\theta}$ is invertible and its inverse is bounded, which can be viewed as
\begin{align} \label{equation:A_inverse_formula}
    (\mathcal{A}_{z}^{N;\theta})^{-1} = \bkt{I' - \bkt{\frac{1}{4\pi} (\mathcal{J}_{z} - \beta_{\star}(\theta)I')}^{-1}\frac{2}{N(N-1)} \sum_{(\mathfrak{i},\mathfrak{i}')\neq (\mathfrak{j},\mathfrak{j}')} \mathcal{S}_{\mathfrak{i},\mathfrak{i}'}\mathcal{G}_{z} \mathcal{S}_{\mathfrak{j},\mathfrak{j}'}^{*} }^{-1} \bkt{\frac{1}{4\pi} (\mathcal{J}_{z} - \beta_{\star}(\theta)I')}^{-1},
\end{align}
The first factor on the right-hand side of \eqref{equation:A_inverse_formula} is understood via the Neumann series $(I-A)^{-1} = \sum_{n=0}^{\infty} A^{n}$, which converges absolutely in the operator norm. The estimate \eqref{equation:A_operator_norm} follows by using \eqref{estimate_J} and \eqref{equattion:G_series_ratio}. In particular, the holomorphicity of $(\mathcal{A}_{z}^{N;\theta})^{-1}$ follows from the expression \eqref{equation:A_inverse_formula}, together with Lemma \ref{lemma:S_is_holomorphic}. 
\end{proof}

\subsection{Holomorphic extension of the resolvent formula} \label{section:properties_principal_operator}

In view of \eqref{definition:Psi_0}, the key step in establishing an
analytic continuation of $\Psi_{N}^{\theta}$ to $\mathscr{O}_{N}^{\theta}$
is the following holomorphicity for the inverse of the principal
operator.

\begin{proposition} \label{proposition:A_inverse_is_bounded}
Suppose that $\theta \in \mathbb{R}$ and $N\geq 3$. Then 
\begin{align}\label{equation:A_inverse_is_bounded}
    (\mathcal{A}_{z}^{N;\theta})^{-1}: \mathscr{O}_{N}^{\theta} \mapsto  B(L^{2}_{\textup{sym}'}(\mathbb{R}^{2N-2}))
\end{align}
defines a holomorphic function. 
\end{proposition}

\subsubsection{\textbf{Identities for the renormalized principal operator}} \label{section:renormalized_principal_operator}

Let us now prove Lemma \ref{proposition:Hilbert_identity}. The key step is the following bounded-extension identity for $\mathcal{J}_{z} - \mathcal{J}_{w}$.

\begin{lemma} \label{lemma:resolvent_identity}
Assume that  $N\geq 3$.  For each $z,w \in \mathbb{C}_{<0}$, the difference
$\mathcal{J}_{z} - \mathcal{J}_{w}$ and $\mathcal{S}_{\mathfrak{i},\mathfrak{i}'}(\mathcal{G}_{z} - \mathcal{G}_{w})\mathcal{S}^{*}_{\mathfrak{i},\mathfrak{i}'}$ extend to bounded operators on $L^{2}_{\textup{sym}'}(\mathbb{R}^{2N-2})$, and their extensions satisfy, for each $(\mathfrak{i},\mathfrak{i}') \in \mathcal{E}_{N}$,
\begin{align} \label{equation:resolvent_identity}
    -\frac{1}{4\pi} (\mathcal{J}_{z} - \mathcal{J}_{w}) = \mathcal{S}_{\mathfrak{i},\mathfrak{i}'}(\mathcal{G}_{z} - \mathcal{G}_{w})\mathcal{S}^{*}_{\mathfrak{i},\mathfrak{i}'}  \quad \text{in} \quad B(L^{2}_{\textup{sym}'}(\mathbb{R}^{2N-2}))  \quad \forall (\mathfrak{i},\mathfrak{i}') \in \mathcal{E}_{N}.
\end{align}
\end{lemma}

Roughly speaking, the main reason why the individual $\mathcal{S}_{\mathfrak{i},\mathfrak{i}'}\mathcal{G}_{z}\mathcal{S}_{\mathfrak{i},\mathfrak{i}'}^{*}$ does not define a bounded operator on $L^{2}_{\textup{sym}'}(\mathbb{R}^{2N-2})$, whereas the difference extends to a bounded operator, is that taking the resolvent difference provides additional decay in the Fourier domain:
\begin{align*}
    \frac{1}{|p|^{2}/2-z} - \frac{1}{|p|^{2}/2-w} = \frac{z-w}{(|p|^{2}/2-z) (|p|^{2}/2-w)}.
\end{align*}
This additional decay compensates for the loss of integrability caused by $\mathcal{S}_{\mathfrak{i},\mathfrak{i}'}$ and $\mathcal{S}_{\mathfrak{i},\mathfrak{i}'}^{*}$. Similarly, the unbounded growth of $\mathcal{J}_{z}$ and $\mathcal{J}_{w}$ at infinity cancels in their difference.

Let us first prove Lemma \ref{proposition:Hilbert_identity} by using Lemma \ref{lemma:resolvent_identity}. 
\begin{proof} [Proof of Lemma \ref{proposition:Hilbert_identity}]
Using \eqref{equation:resolvent_identity} gives
\begin{align*}
    \frac{1}{4\pi} (\mathcal{J}_{z} - \mathcal{J}_{w}) = \frac{2}{N(N-1)}\sum_{(\mathfrak{i,\mathfrak{i}'}) \in \mathcal{E}_{N}} \mathcal{S}_{\mathfrak{i},\mathfrak{i}'}(\mathcal{G}_{w} - \mathcal{G}_{z})\mathcal{S}^{*}_{\mathfrak{i},\mathfrak{i}'}   \quad \text{in} \quad B(L^{2}_{\textup{sym}'}(\mathbb{R}^{2N-2})) .
\end{align*}
Thus, on the $z$-independent domain \eqref{definition:domain_A},
\begin{align} \label{equation:differen_equal_SS_on_dense_set}
    \mathcal{A}_{z}^{N;\theta} - \mathcal{A}_{w}^{N;\theta} 
    = \frac{2}{N(N-1)}\sum_{(\mathfrak{i,\mathfrak{i}'}), (\mathfrak{j},\mathfrak{j}') \in \mathcal{E}_{N}} \mathcal{S}_{\mathfrak{i},\mathfrak{i}'}(\mathcal{G}_{w} - \mathcal{G}_{z})\mathcal{S}^{*}_{\mathfrak{j},\mathfrak{j}'} 
    = -\frac{2}{N(N-1)}(z-w) \mathcal{S}_{z} \widetilde{\mathcal{S}}_{w}. 
\end{align}
Here we have used the resolvent identity for Laplace operator
\begin{align} \label{free_resolvent_identity}
    \mathcal{G}_{\omega} - \mathcal{G}_{\omega'} = (\omega-\omega') \mathcal{G}_{\omega}\mathcal{G}_{\omega'} \quad \forall \omega,\omega'\in \mathbb{C}_{<0}
\end{align}
in the second equality. By Lemma \ref{lemma:estimate_from_GQT}, the right-hand side of \eqref{equation:differen_equal_SS_on_dense_set} is bounded. One easily shows that the domain \eqref{definition:domain_A} is dense in $L^{2}_{\textup{sym}'}(\mathbb{R}^{2N-2})$. Thus, the difference operator $\mathcal{A}_{z}^{N;\theta} - \mathcal{A}_{w}^{N;\theta}$ has a bounded extension \eqref{equation:differen_equal_SS_on_dense_set}. This implies \eqref{equation:resolvent_identity_principal_operator}.  Finally, \eqref{equation:integral_principal_operator} follows from \eqref{equation:resolvent_identity_principal_operator} by Riemann-sum approximation.

\end{proof}

Let us now turn to the proof of Lemma \ref{lemma:resolvent_identity}.  We first prove that both sides of \eqref{equation:resolvent_identity} extend to bounded operators on $L^{2}_{\textup{sym}'}(\mathbb{R}^{2N-2})$. 
\begin{lemma} \label{lemma:both_bounded}
For each $z,w \in \mathbb{C}_{<0}$, there exists a positive constant $C = C(z,w)$ such that
\begin{align} \label{equation:both_bounded}
    ||\mathcal{S}_{\mathfrak{i},\mathfrak{i}'}(\mathcal{G}_{z} - \mathcal{G}_{w})\mathcal{S}^{*}_{\mathfrak{i},\mathfrak{i}'}||_{B(L^{2}(\mathbb{R}^{2N-2}))} \bigvee ||\mathcal{J}_{z} - \mathcal{J}_{w}||_{B(L^{2}(\mathbb{R}^{2N-2}))}\leq C   \quad \forall (\mathfrak{i},\mathfrak{i}') \in \mathcal{E}_{N}.
\end{align}
\end{lemma}
\begin{proof}
Let us first focus on the operator $\mathcal{S}_{\mathfrak{i},\mathfrak{i}'}(\mathcal{G}_{z} - \mathcal{G}_{w}) \mathcal{S}^{*}_{\mathfrak{i},\mathfrak{i}'}$. Without loss of generality, we may assume that $(\mathfrak{i},\mathfrak{i}') = (1,2)$. Fix $\psi \in L^{2}(\mathbb{R}^{2N-2})$ and write $f := \mathcal{S}_{1,2}h$ and $h:= (\mathcal{G}_{z} - \mathcal{G}_{w}) \mathcal{S}^{*}_{1,2} \psi$. From \eqref{definition:adjoint_lowering_operator_-1} and \eqref{definition:adjoint_lowering_operator}, we know that
\begin{align} \label{equation:representation_fourier_f}
    \widehat{f}(q_{2-N}) &= \frac{1}{2\pi} \int_{\mathbb{R}^{2}} dq_{1} \widehat{h}(\mathcal{M}_{1,2}(q_{1},q_{2},\ldots,q_{N})) \nonumber\\
    &= \frac{1}{2\pi} \int_{\mathbb{R}^{2}}  dq_{1} \bkt{\frac{1}{|p|^{2}/2-z} - \frac{1}{|p|^{2}/2-w} } \frac{1}{2\pi} \widehat{\psi}(p_{1}+p_{2},p_{\overline{12}}) \nonumber\\
    &= \frac{1}{(2\pi)^{2}} \widehat{\psi}(q_{2},q_{\overline{12}})
    \int_{\mathbb{R}^{2}}  dq_{1} \bkt{\frac{1}{|p|^{2}/2-z} - \frac{1}{|p|^{2}/2-w} }. 
\end{align}
Here we have set
\begin{align*}
    p = \mathcal{M}_{1,2}(q)
    = \bkt{\frac{1}{2}q_{2}+q_{1},\frac{1}{2}q_{2}-q_{1},q_{3},\ldots,q_{N}},
\end{align*}
which implies
\begin{align*}
    \widehat{\psi}(p_{1}+p_{2},p_{\overline{12}}) = \widehat{\psi}(q_{2},q_{\overline{12}}).
\end{align*}

Now we compute the integral on the right-hand side of \eqref{equation:representation_fourier_f}. Fix $R \geq 1$. Since
\begin{align*}
    |p|^{2} = 2|q_{1}|^{2} + \frac{1}{2}|q_{2}|^{2}+|q_{\overline{12}}|^{2}
    = 2|q_{1}|^{2} + |q|_{2-N}^{2},
\end{align*}
we decompose this integral as follows:
\begin{align} \label{renormalized_identity_R}
    &\int_{\mathbb{R}^{2}}  dq_{1} \bkt{\frac{1}{|p|^{2}/2-z} - \frac{1}{|p|^{2}/2-w} }
    = (2\pi) \int_{0}^{\infty} rdr \bkt{\frac{1}{r^{2}+|q|_{2-N}^{2}/2 - z} - \frac{1}{r^{2}+|q|_{2-N}^{2}/2 - w}} \nonumber\\
    &= \pi \int_{0}^{\infty} dr \bkt{\frac{1}{r+|q|_{2-N}^{2}/2 - z} - \frac{1}{r+|q|_{2-N}^{2}/2 - w}} 
    = \int_{0}^{R}\ldots + \int_{R}^{\infty}\ldots = \pi\sum_{j=0}^{2} I_{j}(q_{2},q_{\overline{12}}),
\end{align}
where we have set
\begin{align*}
    I_{j}(q_{2},q_{\overline{12}}) := 
    \begin{cases}
        \displaystyle \log \bkt{\frac{|q|_{2-N}^{2}/2 - w}{|q|_{2-N}^{2}/2 - z} }, \quad &\text{if } j = 0;\\
        \displaystyle \log \bkt{\frac{R+|q|_{2-N}^{2}/2 - z}{R+|q|_{2-N}^{2}/2 - w} }, \quad &\text{if } j = 1;\\
        \displaystyle \int_{R}^{\infty} dr \bkt{\frac{z-w}{(r+|q|_{2-N}^{2}/2 - z)(r+|q|_{2-N}^{2}/2 - w)}}, \quad &\text{if } j = 2.
    \end{cases}
\end{align*}
To bound the case of $j = 0,1$, write
\begin{align*}
    \frac{|q|_{2-N}^{2}/2 - w}{|q|_{2-N}^{2}/2 - z} = 1+\frac{z-w}{|q|_{2-N}^{2}/2 - z}, \quad \frac{R+|q|_{2-N}^{2}/2 - z}{R+|q|_{2-N}^{2}/2 - w} = 1+\frac{w-z}{R+|q|_{2-N}^{2}/2 - w}.
\end{align*}
It follows that
\begin{align} \label{bound:z_and_w}
    \sup_{q_{2-N} \in \mathbb{R}^{2N-2}} \abs{\frac{|q|_{2-N}^{2}/2 - w}{|q|_{2-N}^{2}/2 - z}}
    \leq 1+\frac{|z-w|}{|\Re(z)|}, \quad \sup_{q_{2-N} \in \mathbb{R}^{2N-2}} \abs{\frac{R+|q|_{2-N}^{2}/2 - z}{R+|q|_{2-N}^{2}/2 - w}} \leq 1+|w-z|.
\end{align}
In particular, interchanging $z$ and $w$ in \eqref{bound:z_and_w} also yields similar uniform bounds. Thus,
\begin{align} \label{equation:I_j_q_2_q_12}
    \sup_{j \in \{0,1,2\}, \; R>1,\; q_{2-N} \in \mathbb{R}^{2N-2}} |I_{j}(q_{2},q_{\overline{12}})| \leq C_{\eqref{equation:I_j_q_2_q_12}},
\end{align}
where $C_{\eqref{equation:I_j_q_2_q_12}} = C_{\eqref{equation:I_j_q_2_q_12}}(z,w)$. As a consequence,
\begin{align*}
    || \mathcal{S}_{1,2}(\mathcal{G}_{z} - \mathcal{G}_{w}) \mathcal{S}^{*}_{1,2} \psi||_{L^{2}(\mathbb{R}^{2N-2})}
    \leq \frac{3 C_{\eqref{equation:I_j_q_2_q_12}}}{4\pi}  || \psi||_{L^{2}(\mathbb{R}^{2N-2})}.
\end{align*}

To bound $\mathcal{J}_{z} - \mathcal{J}_{w}$, notice that its Fourier multiplier is
\begin{align*}
    \log\bkt{\frac{|q|_{2-N}^{2}}{2} -z} - \log\bkt{\frac{|q|_{2-N}^{2}}{2} -w} = -I_{0}(q_{2},q_{\overline{12}}).
\end{align*}
Similarly, we have
\begin{align*}
    ||\mathcal{J}_{z} - \mathcal{J}_{w}||_{B(L^{2}_{\textup{sym}'}(\mathbb{R}^{2N-2}))} \leq C_{\eqref{equation:I_j_q_2_q_12}}.
\end{align*}
The proof is complete.
\end{proof}

Now we are ready to prove \eqref{equation:resolvent_identity}. 
\begin{proof} [Proof of Lemma \ref{lemma:resolvent_identity}] 
Without loss of generality, we may
assume that $(\mathfrak{i},\mathfrak{i}') = (1,2)$. Let $\psi,\varphi \in L^{2}_{\textup{sym}'}(\mathbb{R}^{2N-2})$. Fix $R>1$. By using \eqref{equation:representation_fourier_f}, \eqref{renormalized_identity_R}, and \eqref{equation:I_j_q_2_q_12}, applying the dominated convergence theorem yields that
\begin{align*}
    &\langle \varphi,\mathcal{S}_{\mathfrak{i},\mathfrak{i}'}(\mathcal{G}_{z} - \mathcal{G}_{w})\mathcal{S}^{*}_{\mathfrak{i},\mathfrak{i}'}\psi\rangle_{L^{2}_{\textup{sym}'}(\mathbb{R}^{2N-2})} = \frac{1}{4\pi }\int_{\mathbb{R}^{2N-2}} dq_{2-N} \overline{\widehat{\varphi}(q_{2-N})} \widehat{\psi}(q_{2-N})  \sum_{j=0}^{2} I_{j}(q_{2},q_{\overline{12}}) \\
    &= \frac{1}{4\pi }\int_{\mathbb{R}^{2N-2}} dq_{2-N} \overline{\widehat{\varphi}(q_{2-N}) }\widehat{\psi}(q_{2-N})  I_{0}(q_{2-N})  + o(1)= \frac{-1}{4\pi}\langle \varphi, (\mathcal{J}_{z} - \mathcal{J}_{w}) \psi\rangle_{L^{2}_{\textup{sym}'}(\mathbb{R}^{2N-2})}+o(1).
\end{align*}
Here $o(1)$ tends to zero as $R\to \infty$. The proof is complete.
\end{proof}

\subsubsection{\textbf{Invertibility of the principal operator}} \label{section:invertibility}
Let us now  use the identity for the difference operator to show the invertibility of $\mathcal{A}_{z}^{N;\theta}$ on the half-plane
\begin{align*}
    D_{N}^{\theta} := \{z\in \mathbb{C}:\Re(z) < \inf \mathscr{H}_{N}^{\theta;\textup{sym}}\}.
\end{align*}

\begin{lemma} \label{lemma:invertibility_A}
Assume that $\theta \in \mathbb{R}$ and $N\geq 3$. For each $z \in D_{N}^{\theta}$, $\mathcal{A}_{z}^{N;\theta}$ is invertible. 
\end{lemma}

Fix $z \in D_{N}^{\theta}$ and choose $w \in \mathscr{O}'^{\theta}_{N}$. By Lemma \ref{lemma:A_operator_norm} and Theorem \ref{theorem:resolvent_formula}, we know that
\begin{align} \label{equation:resolvent_w}
    (\mathcal{A}_{w}^{N;\theta})^{-1} \in B(L^{2}_{\textup{sym}'}(\mathbb{R}^{2N-2})), \quad (\mathscr{H}_{N}^{\theta;\textup{sym}}-w)^{-1} = \mathcal{G}_{w}+\frac{2}{N(N-1)} \widetilde{\mathcal{S}}_{w} (\mathcal{A}_{w}^{N;\theta})^{-1} \mathcal{S}_{w}.
\end{align}
Thus, Lemma \ref{proposition:Hilbert_identity} shows that
\begin{align*}
    \mathcal{A}_{z}^{N;\theta} 
    &= \mathcal{A}_{w}^{N;\theta} - \frac{2}{N(N-1)}(z-w) \mathcal{S}_{z} \widetilde{\mathcal{S}}_{w} = \mathcal{B}_{z,w}^{\theta}\mathcal{A}_{w}^{N;\theta}.
\end{align*}
Here, by Lemma \ref{lemma:estimate_from_GQT}, $\mathcal{B}_{z,w}^{\theta}$ given below is a bounded operator on $L^{2}_{\textup{sym}'}(\mathbb{R}^{2N-2})$:
\begin{align*}
    \mathcal{B}_{z,w}^{\theta} :=  I' - \frac{2}{N(N-1)}(z-w) \mathcal{S}_{z} \widetilde{\mathcal{S}}_{w} (\mathcal{A}_{w}^{N;\theta})^{-1}.
\end{align*}
Due to the invertibility of $\mathcal{A}_{w}^{N;\theta}$, to prove Lemma \ref{lemma:invertibility_A}, it is enough to prove the invertibility of $\mathcal{B}_{z,w}^{\theta}$. The main idea is to apply the algebraic property in Lemma \ref{lemma:algebric_property}. Hence, it is enough to prove the invertibility of $\widetilde{\mathcal{B}}_{z,w}^{\theta}$ , where $\widetilde{\mathcal{B}}_{z,w}^{\theta}$ is a bounded operator on $L^{2}_{\textup{sym}}(\mathbb{R}^{2N})$ given by
\begin{align*}
    \widetilde{\mathcal{B}}_{z,w}^{\theta} := I - \frac{2}{N(N-1)}(z-w) \widetilde{\mathcal{S}}_{w} (\mathcal{A}_{w}^{N;\theta})^{-1} \mathcal{S}_{z}.
\end{align*}
\begin{lemma} \label{lemma:algebric_property}
Let $H_{1}$ and $H_{2}$ be two Hilbert spaces and let
\begin{align*}
    \mathcal{U} \in B(H_{1},H_{2}), \quad \mathcal{V}\in B(H_{2},H_{1}).
\end{align*}
Then  $I_{2}-\mathcal{U}\mathcal{V}$ is invertible if and only if $I_{1}-\mathcal{V}\mathcal{U}$ is invertible.
\end{lemma}
\begin{proof}
One easily shows that the invertibility of $I_{1}-\mathcal{V}\mathcal{U}$ implies
\begin{align*}
    (I_{2} - \mathcal{U} \mathcal{V})^{-1} = I_{2} + \mathcal{U} (I_{1}-\mathcal{V}\mathcal{U})^{-1} \mathcal{V}.
\end{align*}
\end{proof}

\begin{proof} [Proof of Lemma \ref{lemma:invertibility_A}]
By the resolvent formula \eqref{equation:resolvent_w} and the free resolvent identity \eqref{free_resolvent_identity}, the invertibility of $\widetilde{\mathcal{B}}_{z,w}^{\theta}$ now follows from the identity
\begin{align*}
    (\mathscr{H}_{N}^{\theta;\textup{sym}}-z)(\mathscr{H}_{N}^{\theta;\textup{sym}}-w)^{-1}
    &= I-(z-w) \bkt{\mathcal{G}_{w}+\frac{2}{N(N-1)} \widetilde{\mathcal{S}}_{w} (\mathcal{A}_{w}^{N;\theta})^{-1} \mathcal{S}_{w}}\\
    &= \mathcal{G}_{z}^{-1} \mathcal{G}_{w}-\frac{2(z-w)}{N(N-1)} \widetilde{\mathcal{S}}_{w} (\mathcal{A}_{w}^{N;\theta})^{-1} \mathcal{S}_{w} = \widetilde{\mathcal{B}}_{z,w}^{\theta} \mathcal{G}_{z}^{-1} \mathcal{G}_{w},
\end{align*}
where we have used $\mathcal{S}_{w} \mathcal{G}_{w}^{-1} \mathcal{G}_{z} = \mathcal{S}_{z}$ in the last equality.
Indeed, since $z\in D_{N}^{\theta}$, both  $(\mathscr{H}_{N}^{\theta;\textup{sym}}-z)(\mathscr{H}_{N}^{\theta;\textup{sym}}-w)^{-1}$ and $\mathcal{G}_{z}^{-1} \mathcal{G}_{w}$ are invertible on $L^{2}_{\textup{sym}}(\mathbb{R}^{2N})$. The proof of Lemma \ref{lemma:invertibility_A}
 is complete.
\end{proof}
 
\subsubsection{\textbf{Bounds for the inverse of the principal operator}} \label{section:bounded_inverse}

The following lemma gives a quantitative bound for the operator norm of $(\mathcal{A}_{z}^{N;\theta})^{-1}$ on $\mathscr{G}_{N}^{\theta}$.

\begin{lemma} \label{lemma:A_inverse_operator_norm_estimate}
Suppose that $\theta \in \mathbb{R}$ and $N\geq 3$. For all $z \in \mathscr{G}_{N}^{\theta}$, $(\mathcal{A}_{z}^{N;\theta})^{-1} \in B(L^{2}_{\textup{sym}'}(\mathbb{R}^{2N-2}))$ and 
\begin{align} \label{equation:inverse_A_operator_bound_above_real_line}
    ||(\mathcal{A}_{z}^{N;\theta})^{-1}||_{B(L^{2}_{\textup{sym}'}(\mathbb{R}^{2N-2}))}  \leq   \frac{2}{E_{N}^{\theta}(\Re(z))} < \infty.
\end{align}
\end{lemma}

We first establish the boundedness of the inverse throughout the  half-plane $D_{N}^{\theta}$.

\begin{lemma} \label{lemma:inverse_A_bounded}
Assume that $\theta \in \mathbb{R}$ and $N\geq 3$. For $z \in D_{N}^{\theta}$, $(\mathcal{A}_{z}^{N;\theta})^{-1} \in B(L^{2}_{\textup{sym}'}(\mathbb{R}^{2N-2}))$.
\end{lemma}
\begin{proof}
Fix $z\in D_{N}^{\theta}$. By Lemma \ref{lemma:invertibility_A}, the inverse
$(\mathcal{A}_{z}^{N;\theta})^{-1}$ is defined on the whole space $L^{2}_{\textup{sym}'}(\mathbb{R}^{2N-2})$.
By the closed graph theorem, it therefore suffices to show that $(\mathcal{A}_{z}^{N;\theta})^{-1}$ is closed, or equivalently, that $\mathcal{A}_{z}^{N;\theta}$ is closed. By Lemma \ref{lemma:estimate_from_GQT} and
\eqref{recall_principal_operator}, $\mathcal{A}_{z}^{N;\theta}$ is a bounded perturbation of
$\mathcal{J}_{z}$. Since the Fourier multiplier $\mathcal{J}_{z}$ is closed on the natural domain \eqref{definition:domain_J}, $\mathcal{A}_{z}^{N;\theta}$ is also closed on the same domain \eqref{definition:domain_A}. The proof of Lemma \ref{lemma:inverse_A_bounded} is complete.
\end{proof}

To prove Lemma \ref{lemma:A_inverse_operator_norm_estimate}, it is enough to prove a real version of \eqref{equation:inverse_A_operator_bound_above_real_line}.
\begin{lemma} \label{lemma:A_inverse_operator_norm_estimate_real}
Assume that $\theta \in \mathbb{R}$ and $N\geq 3$. For each $z \in D_{N}^{\theta} \cap \mathbb{R}$, 
\begin{align} \label{equation:A_op_real_z}
    ||(\mathcal{A}_{z}^{N;\theta})^{-1}||_{B(L^{2}_{\textup{sym}'}(\mathbb{R}^{2N-2}))}  =   \frac{1}{E_{N}^{\theta}(z)} < \infty.
\end{align}
\end{lemma}

Let us now conclude Lemma \ref{lemma:A_inverse_operator_norm_estimate} by using Lemma \ref{lemma:A_inverse_operator_norm_estimate_real}.
\begin{proof} [Proof of Lemma \ref{lemma:A_inverse_operator_norm_estimate}]
By Lemma \ref{lemma:estimate_from_GQT}, there exists an absolute constant $C_{\eqref{definition:constant_set_O}} > 0$ such that
\begin{align} \label{definition:constant_set_O}
    &\sup_{z,w \in D_{N}^{\theta}}\frac{2}{N(N-1)} ||\mathcal{S}_{z} ||_{B(L^{2}_{\textup{sym}}(\mathbb{R}^{2N}),L^{2}_{\textup{sym}'}(\mathbb{R}^{2N-2}))} ||\widetilde{\mathcal{S}}_{w} ||_{B(L^{2}_{\textup{sym}'}(\mathbb{R}^{2N-2}),L^{2}_{\textup{sym}}(\mathbb{R}^{2N}))}  \nonumber \\
    &\leq \frac{C_{\eqref{definition:constant_set_O}} N^{2}}{(-\inf \mathscr{H}_{N}^{\theta;\textup{sym}})} < \infty.
\end{align}
The finiteness of the right-hand side uses the strict negativity
\begin{align} \label{negative_0}
    \inf\mathscr{H}_{N}^{\theta;\textup{sym}}<0
    \qquad\forall N\geq 3,
\end{align}
which follows by comparing moments of different orders.
Indeed, Hölder’s inequality gives
\begin{align*} 
    \mathbb{E}[\mathscr{Z}_{t}^{\theta}(G_{1},G_{1})^{N}]^{\frac{1}{N}} \leq \mathbb{E}[\mathscr{Z}_{t}^{\theta}(G_{1},G_{1})^{M}]^{\frac{1}{M}}  \quad \forall 2\leq N\leq M.
\end{align*}
Taking logarithms, dividing by $t$, and letting $t\to\infty$ together with Lemma \ref{lemma:spectrum_theorem_1} shows that the ground-state energy per particle is nonincreasing:
\begin{align} \label{negative_1}
    \frac{1}{M}\inf \mathscr{H}_{M}^{\theta;\textup{sym}} \leq \frac{1}{N}\inf \mathscr{H}_{N}^{\theta;\textup{sym}} \quad \forall 2\leq N\leq M.
\end{align}
By \cite[Section~2.1, Remark~1 following Theorem~1]{JDimock_2004}, the two-particle ground-state energy is strictly negative:
\begin{align*}
    \inf\mathscr{H}_{2}^{\theta;\textup{sym}}<0.
\end{align*}
Combining this with \eqref{negative_1} yields \eqref{negative_0}.

Since $z \in \mathscr{G}_{N}^{\theta}$, it follows that
\begin{align*}
    |\Im(z)| < \frac{(-\inf \mathscr{H}_{N}^{\theta;\textup{sym}})E_{N}^{\theta}(\Re(z)) }{2 C_{\eqref{definition:constant_set_O}} N^{2}}.
\end{align*}
In addition, since $z \in D_{N}^{\theta}$, applying the identity \eqref{equation:resolvent_identity_principal_operator} for difference operator gives
\begin{align*}
    (\mathcal{A}_{z}^{N;\theta})^{-1}-(\mathcal{A}_{\Re(z)}^{N;\theta})^{-1} = 
    - \frac{2}{N(N-1)}(\Re(z)-z) (\mathcal{A}_{z}^{N;\theta})^{-1}\mathcal{S}_{\Re(z)} \widetilde{\mathcal{S}}_{z} (\mathcal{A}_{\Re(z)}^{N;\theta})^{-1}.
\end{align*}
As a consequence, 
\begin{align} \label{estimate:inverse_A}
    ||(\mathcal{A}_{z}^{N;\theta})^{-1}||_{B(L^{2}_{\textup{sym}'}(\mathbb{R}^{2N-2}))} 
    &\leq \frac{1}{E_{N}^{\theta}(\Re(z))} + \frac{C_{\eqref{definition:constant_set_O}} N^{2}}{(-\inf \mathscr{H}_{N}^{\theta;\textup{sym}})}\frac{|\Im(z)|}{E_{N}^{\theta}(\Re(z))} ||(\mathcal{A}_{z}^{N;\theta})^{-1}||_{B(L^{2}_{\textup{sym}'}(\mathbb{R}^{2N-2}))} \nonumber \\
    &\leq \frac{1}{E_{N}^{\theta}(\Re(z))}+\frac{1}{2} ||(\mathcal{A}_{z}^{N;\theta})^{-1}||_{B(L^{2}_{\textup{sym}'}(\mathbb{R}^{2N-2}))}.
\end{align}
Lemma \ref{lemma:inverse_A_bounded} shows that $||(\mathcal{A}_{z}^{N;\theta})^{-1}||_{B(L^{2}_{\textup{sym}'}(\mathbb{R}^{2N-2}))}$ is finite. Thus, \eqref{estimate:inverse_A} yields \eqref{equation:inverse_A_operator_bound_above_real_line}. The proof is complete. 

\end{proof}

Let us turn to the proof of Lemma \ref{lemma:A_inverse_operator_norm_estimate_real}. We first establish the  positivity of $E_{N}^{\theta}(z)$.
\begin{lemma} \label{lemma:positivity_A}
Suppose that $\theta\in \mathbb{R}$ and $N \geq 3$. Then
\begin{align*}
    E_{N}^{\theta}(z) > 0 \quad \forall z \in D_{N}^{\theta} \cap \mathbb{R}.
\end{align*}
\end{lemma}

Our proof relies on the following variational representation. 
\begin{lemma} \label{lemma:A_is_self_adjoint}
Assume that $\theta \in \mathbb{R}$ and $N\geq 3$. For each $z<0$, $\mathcal{A}_{z}^{N;\theta}$ is self-adjoint. In particular, we have
\begin{align}  \label{equation:EN_and_inf_A}
    E_{N}^{\theta}(z) = \inf \mathcal{A}_{z}^{N;\theta}.
\end{align}
\end{lemma}
\begin{proof}
Since $z < 0$, the free resolvent $\mathcal{G}_{z}$ is symmetric. Thus, Lemma \ref{lemma:estimate_from_GQT} shows that the second part in \eqref{recall_principal_operator} is bounded and symmetric. In addition, since $\mathcal{J}_{z}$ is a symmetric Fourier multiplier, $\mathcal{J}_{z}$ is self-adjoint on the natural domain \eqref{definition:domain_J}, and hence, applying the Kato-Rellich theorem shows that $\mathcal{A}_{z}^{N;\theta}$ is self-adjoint on the same domain. The identity \eqref{equation:EN_and_inf_A} now follows from the variational representation \cite[Theorem 2.19]{teschl2009mathematical}.
\end{proof}

By using the identity \eqref{equation:resolvent_identity_principal_operator} together with the $z$-independent domain \eqref{definition:domain_A}, one easily shows that Lemma \ref{lemma:estimate_from_GQT} yields the continuity of $E_{N}^{\theta}(z)$.
\begin{lemma} \label{lemma:E_is_Lip}
Assume that $\theta \in \mathbb{R}$ and $N\geq 3$. There exists $C>0$ such that 
\begin{align} \label{equation:E_is_Lip}
    |E_{N}^{\theta}(z) - E_{N}^{\theta}(w)| \leq C \frac{|z-w|}{\sqrt{(-z)(-w)}} \quad \forall z,w <0.
\end{align}
\end{lemma}

Now we are ready to prove Lemma \ref{lemma:positivity_A}.
\begin{proof} [Proof of Lemma \ref{lemma:positivity_A}]

To prove this, observe that, for sufficiently negative $z$, $E_{N}^{\theta}(z) > 0$. Indeed, by using Lemma \ref{lemma:estimate_from_GQT}, there exists an absolute constant $C>0$ such that
\begin{align} \label{EN_positive}
    E_{N}^{\theta}(z) \geq \frac{\log(-z)-\beta_{\star}(\theta)}{4\pi} - CN^{2} \quad \forall z < 0.
\end{align}

Now suppose that $E_{N}^{\theta}(z_{0}) \leq 0$ for some $z_{0} \in D_{N}^{\theta} \cap \mathbb{R}$. By  \eqref{EN_positive}, $E_{N}^{\theta}$ is positive for all sufficiently negative $z$. Then, by the continuity in Lemma \ref{lemma:E_is_Lip}, there exists $z'_{0} \in D_{N}^{\theta} \cap \mathbb{R}$ such that $E_{N}^{\theta}(z_{0}') = 0$. In particular, by Lemma \ref{lemma:A_is_self_adjoint}, we obtain
\begin{align}  \label{equation:inf_A_is_zero}
    \inf \mathcal{A}_{z_{0}'}^{N;\theta} = E_{N}^{\theta}(z_{0}') = 0.
\end{align}

Since $\sigma(\mathcal{A}_{z_{0}'}^{N;\theta})$ is closed,  \eqref{equation:inf_A_is_zero} implies that $0 \in \sigma(\mathcal{A}_{z_{0}'}^{N;\theta})$. This is a contradiction. Indeed, by the invertibility in Lemma \ref{lemma:invertibility_A} and the boundedness in Lemma \ref{lemma:inverse_A_bounded}, we know that $0 \in \rho(\mathcal{A}_{z_{0}'}^{N;\theta})$. We now obtain the positivity of $E_{N}^{\theta}$ on $D_{N}^{\theta}\cap \mathbb{R}$.
\end{proof}

Lemma \ref{lemma:A_inverse_operator_norm_estimate_real} now follows from Lemma \ref{lemma:positivity_A}.
\begin{proof} [Proof of Lemma \ref{lemma:A_inverse_operator_norm_estimate_real}]
Fix $z \in D_{N}^{\theta} \cap \mathbb{R}$. Similarly, by Lemma \ref{lemma:invertibility_A} and Lemma \ref{lemma:inverse_A_bounded}, $0 \in \rho(\mathcal{A}_{z}^{N;\theta})$. Thus, we conclude
\begin{align*}
    ||(\mathcal{A}_{z}^{N;\theta})^{-1}||_{B(L^{2}_{\textup{sym}'}(\mathbb{R}^{2N-2}))} = \frac{1}{d(0,\sigma(\mathcal{A}_{z}^{N;\theta}))} = \frac{1}{\inf \sigma(\mathcal{A}_{z}^{N;\theta})} = \frac{1}{E_{N}^{\theta}(z)},
\end{align*}
which is finite by using Lemma \ref{lemma:positivity_A}. The proof is complete.
\end{proof}

\subsubsection{\textbf{Analytic continuation of the inverse of the principal operator}}

Recall that $\mathscr{O}_{N}^{\theta} = \mathscr{G}_{N}^{\theta}\cup \mathscr{O}'^{\theta}_{N}$. By Lemma \ref{lemma:A_operator_norm}, we know that $(\mathcal{A}_{\cdot}^{N;\theta})^{-1}$ is holomorphic on $\mathscr{O}'^{\theta}_{N}$, and hence, it suffices to consider its holomorphicity on the open set $\mathscr{G}_{N}^{\theta}$.
\begin{lemma} \label{lemma:Analytic_continuation}
Suppose that $\theta \in \mathbb{R}$ and $N\geq 3$. Then 
\begin{align} \label{equation:differential_equation_for_A_-1}
    \frac{d}{dz} (\mathcal{A}_{z}^{N;\theta})^{-1}
    =\frac{2}{N(N-1)} (\mathcal{A}_{z}^{N;\theta})^{-1}\mathcal{S}_{z} \widetilde{\mathcal{S}}_{z} (\mathcal{A}_{z}^{N;\theta})^{-1} \quad \forall z \in \mathscr{G}_{N}^{\theta}.
\end{align}
In particular, $(\mathcal{A}_{\cdot}^{N;\theta})^{-1}$ is holomorphic on $\mathscr{G}_{N}^{\theta}$.
\end{lemma}
\begin{proof}

By the identity in Lemma \ref{proposition:Hilbert_identity} for the difference and the invertibility in  Lemma \ref{lemma:invertibility_A}, we have
\begin{align} \label{inverse_equation}
    (\mathcal{A}_{w}^{N;\theta})^{-1} - (\mathcal{A}_{z}^{N;\theta})^{-1}
    = \frac{2}{N(N-1)} (w-z) (\mathcal{A}_{z}^{N;\theta})^{-1}\mathcal{S}_{z} \widetilde{\mathcal{S}}_{w} (\mathcal{A}_{w}^{N;\theta})^{-1} \quad \forall z,w \in \mathscr{G}_{N}^{\theta}.
\end{align}
By Lemma \ref{lemma:S_is_holomorphic}, both $\mathcal{S}_{\cdot}$ and $ \widetilde{\mathcal{S}}_{\cdot}$ are holomorphic on $\mathscr{G}_{N}^{\theta}$. Thus, the equation \eqref{inverse_equation} together with the bound in \eqref{equation:inverse_A_operator_bound_above_real_line} then yields 
the continuity of $(\mathcal{A}_{\cdot}^{N;\theta})^{-1}$ on $\mathscr{G}_{N}^{\theta}$. In particular, taking  $w \to z$, $(\mathcal{A}_{\cdot}^{N;\theta})^{-1}$ satisfies the differential equation \eqref{equation:differential_equation_for_A_-1}. The proof is complete. 
\end{proof}

\subsubsection{\textbf{Proof of Theorem \ref{theorem:resolvent_extension}}}

Applying Lemma \ref{lemma:E_is_Lip} and Lemma \ref{lemma:positivity_A}  gives the continuity of
\begin{align*}
    z \in D_{N}^{\theta} \mapsto \frac{|\Im(z)|}{E_{N}^{\theta}(\Re(z))}.
\end{align*}
This together with \eqref{negative_0} shows that $\mathscr{G}^{\theta}_{N}$ is open and connected. Indeed, by positivity of $E_{N}^{\theta}$, $z\in\mathscr{G}^{\theta}_{N}$ can be joined to $\Re(z)$ by a vertical segment contained in $\mathscr{G}^{\theta}_{N}$. In particular, since $(-\infty,\inf\mathscr{H}_{N}^{\theta;\textup{sym}}) \subseteq \mathscr{G}^{\theta}_{N}$, any two points in $\mathscr{G}^{\theta}_{N}$ can be connected by these line segments. Thus, since $\mathscr{G}^{\theta}_{N}\cap \mathscr{O}'^{\theta}_{N} \neq \emptyset$, $\mathscr{O}^{\theta}_{N} = \mathscr{G}^{\theta}_{N}\cup \mathscr{O}'^{\theta}_{N}$ is a domain. By Lemma \ref{lemma:Analytic_continuation} and  Lemma \ref{lemma:S_is_holomorphic}, we know that $\Psi_{N}^{\theta}$ is holomorphic on $\mathscr{G}_{N}^{\theta}$. In addition, Lemma \ref{lemma:A_operator_norm} gives the holomorphicity  of $\Psi_{N}^{\theta}$ on $\mathscr{O}'^{\theta}_{N}$. Therefore, we conclude the holomorphicity of $\Psi_{N}^{\theta}$ on $\mathscr{O}_{N}^{\theta}$. As a consequence, the resolvent formula \eqref{equation:resolvent_representation_2} follows from \eqref{equation:resolvent_representation} and the identity theorem of holomorphic functions. The proof of Theorem \ref{theorem:resolvent_extension} is complete.

\section{Derivation of the ground-state energy formula}
\label{section:proof_formula_ground_state_energy}

In this section, we prove Proposition \ref{theorem:formula_ground_state_energy}
using the extended resolvent formula
\eqref{equation:resolvent_representation_2}.
It suffices to establish \eqref{equation:equation_for_ground_state_energy}.
Indeed, write
\begin{align*}
    e_{0} := \inf\mathscr{H}_{N}^{\theta;\textup{sym}}<0,
    \qquad c:=(-e_{0})^{1/2},
\end{align*}
where we have used \eqref{negative_0}. Fix $\psi\in D(\mathcal{A}_{e_{0}}^{N;\theta})$ with $||\psi||_{L^{2}_{\textup{sym}'}(\mathbb{R}^{2N-2})}=1$,
and define its $L^{2}$-normalized dilation by
\begin{align*}
    \varphi_{c}(x)
    := c^{-(N-1)}\psi(x/c),
    \qquad x\in\mathbb{R}^{2N-2}.
\end{align*}
Then $||\varphi_{c}||_{L^{2}_{\textup{sym}'}(\mathbb{R}^{2N-2})}=1$ and $\widehat{\varphi_{c}}(p)=c^{N-1}\widehat{\psi}(cp)$.
A change of variables $q=cp$ gives
\begin{align*}
    \langle \psi, \mathcal{A}_{e_{0}}^{N;\theta} \psi \rangle_{L^{2}_{\textup{sym}'}(\mathbb{R}^{2N-2})}
     = \frac{1}{4\pi} \log (- e_{0}) + \langle \varphi_{c}, \mathcal{A}_{-1}^{N;\theta} \varphi_{c} \rangle_{L^{2}_{\textup{sym}'}(\mathbb{R}^{2N-2})}.
\end{align*}
The dilation $\psi\mapsto\varphi_{c}$ is unitary and maps the common
domain \eqref{definition:domain_A} onto itself.
Taking the infimum over all unit vectors
$\psi\in D(\mathcal{A}_{e_{0}}^{N;\theta})$ therefore yields
\begin{align} \label{equation:identity_in_theorem:formula_ground_state_energy}
    E_{N}^{\theta}(e_{0}) = \frac{1}{4\pi} \log (- e_{0})+E_{N}^{\theta}(-1).
\end{align}
By combining this identity with \eqref{equation:equation_for_ground_state_energy}, we conclude \eqref{equation:formula_ground_state_energy}.

Let us now prove \eqref{equation:equation_for_ground_state_energy}. By the continuity in Lemma \ref{lemma:E_is_Lip} and the positivity in Lemma \ref{lemma:positivity_A}, we know that
\begin{align}  \label{equation:E_N_at_ground_state_energy_nonnegative}
    E_{N}^{\theta}(e_{0}) \geq  0.
\end{align}
Now we prove \eqref{equation:equation_for_ground_state_energy} by contradiction. Suppose that  
\begin{align}  \label{equation:contradiction_E_N_>0}
    E_{N}^{\theta}(e_{0}) \geq \delta > 0.
\end{align}
One easily shows that $\mathcal{S}_{x}\widetilde{\mathcal{S}}_{x} \geq 0$ for all $x \in D_{N}^{\theta} \cap \mathbb{R}$. This fact together with the integral equation \eqref{equation:integral_principal_operator} and \eqref{equation:contradiction_E_N_>0} implies that 
\begin{align}  \label{equation:A_z_larger_A_E}
    \mathcal{A}_{z}^{N;\theta} \geq \mathcal{A}_{e_{0}}^{N;\theta} \geq \delta I' \quad \forall z \in D_{N}^{\theta} \cap \mathbb{R}.
\end{align}
Thus, by using Lemma \ref{lemma:inverse_A_bounded} and Lemma \ref{lemma:A_is_self_adjoint}, we conclude
\begin{align} \label{equation:uniform_bound_A_inverse}
    \sup_{z \in D_{N}^{\theta} \cap \mathbb{R}} ||(\mathcal{A}_{z}^{N;\theta})^{-1}||_{B(L^{2}_{\textup{sym}'}(\mathbb{R}^{2N-2}))} \leq \delta^{-1}.
\end{align}
The extended resolvent formula \eqref{equation:resolvent_representation_2} gives
\begin{align*}
    &\sup_{z \in D_{N}^{\theta} \cap \mathbb{R}}||(\mathscr{H}_{N}^{\theta;\textup{sym}} - z)^{-1}||_{B({L_{\textup{sym}}^{2}(\mathbb{R}^{2N})})}
    \leq \sup_{z \in D_{N}^{\theta} \cap \mathbb{R}} \Norm{\mathcal{G}_{z} + \frac{2}{N(N-1)} \widetilde{\mathcal{S}}_{z} (\mathcal{A}_{z}^{N;\theta})^{-1}\mathcal{S}_{z} }_{B({L_{\textup{sym}}^{2}(\mathbb{R}^{2N})})} \leq  C,
\end{align*}
where we have used \eqref{equation:uniform_bound_A_inverse} and Lemma \ref{lemma:estimate_from_GQT} in the last estimate. This is a contradiction. Indeed, by the spectral theorem, we know that 
\begin{align*}
    \lim_{z \uparrow \inf \mathscr{H}_{N}^{\theta;\textup{sym}}} ||(\mathscr{H}_{N}^{\theta;\textup{sym}} - z)^{-1}||_{B({L_{\textup{sym}}^{2}(\mathbb{R}^{2N})})}
    = \lim_{z \uparrow \inf \mathscr{H}_{N}^{\theta;\textup{sym}}} \frac{1}{|z-\inf \mathscr{H}_{N}^{\theta;\textup{sym}}|} = +\infty.
\end{align*}
The proof of \eqref{equation:equation_for_ground_state_energy} is complete.

\section{Asymptotics of the mean-field exponent} \label{section:mean_field_exponent}

In this section, we prove Theorem \ref{proposition:mean_field_two_side}.
We begin by rewriting the mean-field exponent in a form that separates
the three- and four-particle contributions.

\begin{lemma} \label{lemma:expression_mean_field_exponent}
Assume that $N \geq 4$ and $\theta \in \mathbb{R}$. Then 

\begin{align} \label{equation:expression_mean_field_exponent}
    &E_{N,\theta}^{\textup{M-F}} = \inf_{v\star u \in \mathscr{A}_{\textup{M-F}}(\mathbb{R}^{2N-2})} \frac{1}{4\pi} \langle v\star u,\mathcal{J}_{-1} v\star u \rangle_{L^{2}_{\textup{sym}'}(\mathbb{R}^{2N-2})} - \frac{\beta_{\star}(\theta)}{4\pi} -I^{(3)}_{N}[v,u]-I^{(4)}_{N}[v,u].
\end{align}
\end{lemma}

The three- and four-particle contributions appearing in \eqref{equation:expression_mean_field_exponent} are defined by
\begin{align*}
    I^{(n)}_{N}[v,u] = 
    \begin{cases}
        \displaystyle 2(N-2) K^{(n)}_{N}[v,u], \quad &\text{if } n = 3;\\
        \displaystyle \binom{N-2}{2} K^{(n)}_{N}[v,u], \quad &\text{if } n = 4.
    \end{cases}
\end{align*}
Writing $U:=\widehat{u}$ and $V:=\widehat{v}$, we set
\begin{align*}
    K^{(3)}_{N}[v,u]
    &:= \frac{1}{(2\pi)^{2}} \int_{0}^{\infty} dt \exp(-t) A_{t;u}^{N-3} \int_{\mathbb{R}^{2}} dp \exp\bkt{-\frac{t}{2}|p|^{2}} \abs{\int_{\mathbb{R}^{2}} dq \exp\bkt{-\frac{t}{2}|q|^{2} }\overline{V}(p+q)U(q) }^{2},\\
    K^{(4)}_{N}[v,u]
    &:= \frac{1}{(2\pi)^{2}} \int_{0}^{\infty} dt \exp(-t) A_{t;u}^{N-4} \abs{ \langle V,X_{t;u} \rangle_{L^{2}(\mathbb{R}^{2})} }^{2},
\end{align*}
where
\begin{align*}
    A_{t;u} :=  \int_{\mathbb{R}^{2}} dq \exp\bkt{-\frac{t}{2} |q|^{2}} |U(q)|^{2}, \quad
    X_{t;u}(k) := \int_{\mathbb{R}^{2}} dq\exp\bkt{-\frac{t}{2}(|q|^{2}+|k-q|^{2})} U(q)U(k-q).
\end{align*}

\begin{proof} [Proof of Lemma \ref{lemma:expression_mean_field_exponent}]
For each fixed collision label, there are $2(N-2)$ other labels
sharing exactly one particle with it and $\binom{N-2}{2}$ labels
disjoint from it. By symmetry, we therefore have
\begin{align*}
    &\frac{2}{N(N-1)}\sum_{(\mathfrak{i},\mathfrak{i}') \neq (\mathfrak{j},\mathfrak{j}')} \langle v\star u,\mathcal{S}_{\mathfrak{i},\mathfrak{i}'} \mathcal{G}_{-1} \mathcal{S}_{\mathfrak{j},\mathfrak{j}'}^{*} v\star u \rangle_{L^{2}_{\textup{sym}'}(\mathbb{R}^{2N-2})}\\
    &= 2(N-2) \langle v\star u,\mathcal{S}_{1,2} \mathcal{G}_{-1} \mathcal{S}_{1,3}^{*} v\star u \rangle_{L^{2}_{\textup{sym}'}(\mathbb{R}^{2N-2})} + \binom{N-2}{2} \langle v\star u,\mathcal{S}_{1,2} \mathcal{G}_{-1} \mathcal{S}_{3,4}^{*} v\star u \rangle_{L^{2}_{\textup{sym}'}(\mathbb{R}^{2N-2})}.
\end{align*}
By using the integral representation
\begin{align} \label{equation:inverse_and_exponential_integral}
    \frac{1}{|q|^{2}/2+1} = \int_{0}^{\infty} \exp\bkt{-t-\frac{t}{2}|q|^{2}} dt,
\end{align}
expressing the following inner products in the Fourier domain yields
\begin{align*}
    \langle v\star u,\mathcal{S}_{1,2} \mathcal{G}_{-1} \mathcal{S}_{1,3}^{*} v\star u \rangle_{L^{2}_{\textup{sym}'}(\mathbb{R}^{2N-2})} &= K_{N}^{(3)}[v,u],\\
    \langle v\star u,\mathcal{S}_{1,2} \mathcal{G}_{-1} \mathcal{S}_{3,4}^{*} v\star u \rangle_{L^{2}_{\textup{sym}'}(\mathbb{R}^{2N-2})} &= K_{N}^{(4)}[v,u].
\end{align*}
Substituting these identities into $E_{N,\theta}^{\textup{M-F}}$ gives \eqref{equation:expression_mean_field_exponent}.
\end{proof}

\subsection{Lower bound for the mean-field exponent}

Let us now show the lower bound in Theorem \ref{proposition:mean_field_two_side}.
The positivity of the Fourier multiplier $\mathcal{J}_{-1}$ yields
\begin{align*}
    E_{N,\theta}^{\textup{M-F}} \geq - \frac{\beta_{\star}(\theta)}{4\pi} -\sup_{v\star u \in \mathscr{A}_{\textup{M-F}}(\mathbb{R}^{2N-2})} I_{N}^{(4)}[v,u] - \sup_{v\star u \in \mathscr{A}_{\textup{M-F}}(\mathbb{R}^{2N-2})} I_{N}^{(3)}[v,u].
\end{align*}
Thus, it suffices to show the following estimates.
\begin{lemma} \label{lemma:sup_K_4}
Suppose that $\theta \in \mathbb{R}$. For all sufficiently large $N$,
\begin{align} \label{estimate:sup_K_4}
    \sup_{v\star u \in \mathscr{A}_{\textup{M-F}}(\mathbb{R}^{2N-2})} K_{N}^{(n)}[v,u] \leq     
    \begin{cases}
        \displaystyle \frac{2 C_{\textup{G-N}}}{N}, \quad &\text{if } n = 4;\\
        \displaystyle O(1/\sqrt{N}), \quad &\text{if } n  =3.
    \end{cases}
\end{align}
\end{lemma}
\begin{proof}
We now bound the four-particle term in \eqref{estimate:sup_K_4}. By the Cauchy–Schwarz inequality, we have 
\begin{align} \label{C_S_oper}
    K_{N}^{(4)}[v,u]  &\leq \frac{1}{(2\pi)^{2}} \int_{0}^{\infty} dt  A_{t;u}^{N-4} ||X_{t;u}||_{L^{2}(\mathbb{R}^{2})}^{2},
\end{align}
where we have used $||V||_{L^{2}(\mathbb{R}^{2})}  = 1$.
The convolution formula and Plancherel's theorem give
\begin{align*}
    ||X_{t;u}||_{L^{2}(\mathbb{R}^{2})}^{2} = (2\pi)^{2} ||G_{t}u||_{L^{4}(\mathbb{R}^{2})}^{4}.
\end{align*}
In particular, applying the two-dimensional Gagliardo–Nirenberg inequality \eqref{equation:Gagliardo_Nirenberg} gives
\begin{align} \label{G_N_oper}
    ||G_{t}u||_{L^{4}(\mathbb{R}^{2})}^{4}
    \leq C_{\textup{G-N}} ||G_{t}u||_{L^{2}(\mathbb{R}^{2})}^{2} ||\nabla G_{t}u||_{L^{2}(\mathbb{R}^{2})}^{2} = C_{\textup{G-N}} A_{2t;u} (-2A'_{2t;u}),
\end{align}
where the last equality follows from Plancherel's theorem. Consequently,
\begin{align*}
    K_{N}^{(4)}[v,u]  \leq -2C_{\textup{G-N}}\int_{0}^{\infty} dt   A_{t;u}^{N-4} A_{2t;u} A'_{2t;u} = -C_{\textup{G-N}}\int_{0}^{\infty} dt   A_{t/2;u}^{N-4} A_{t;u} A'_{t;u}.
\end{align*}
In addition, since $||U||_{L^{2}(\mathbb{R}^{2})} = 1$, applying Jensen's inequality gives $A_{t/2;u}^{2} \leq A_{t;u}$. This implies 
\begin{align*}
    K_{N}^{(4)}[v,u]  
    \leq -C_{\textup{G-N}}\int_{0}^{\infty} dt   A_{t;u}^{(N-4)/2} A_{t;u} A'_{t;u}
    = C_{\textup{G-N}}\int_{0}^{1} a^{(N-2)/2} da = \frac{2C_{\textup{G-N}}}{N}.
\end{align*}

Let us now bound  the three-particle term in \eqref{estimate:sup_K_4}. Write
\begin{align*}
    K_{N}^{(3)}[v,u] = \frac{1}{(2\pi)^{2}} \int_{0}^{\infty} dt \exp(-t) A_{t;u}^{N-3} B_{t;u,v}, 
\end{align*}
where we have set
\begin{align*}
    B_{t;u,v} := \int_{\mathbb{R}^{2}} dp \exp\bkt{-\frac{t}{2}|p|^{2}} \abs{\int_{\mathbb{R}^{2}} dq \exp\bkt{-\frac{t}{2}|q|^{2} }\overline{V}(p+q)U(q) }^{2}. 
\end{align*}
To bound $B_{t;u,v}$, by using 
\begin{align*}
    \frac{|p|^{2}}{4} + \frac{|q|^{2}}{2} \geq \frac{|q|^{2}}{4} + \frac{|p+q|^{2}}{8} \quad \forall p,q \in \mathbb{R}^{2},
\end{align*}
we obtain
\begin{align*}
    B_{t;u,v} \leq ||h_{t}*g_{t}||_{L^{2}(\mathbb{R}^{2})}^{2}.
\end{align*}
Here we have set
\begin{align*}
    h_{t}(p) := f_{t}(-p), \quad f_{t}(p) :=  \exp\bkt{-\frac{t}{4}|p|^{2}} |U(p)|, \quad g_{t}(p) :=  \exp\bkt{-\frac{t}{8}|p|^{2}} |V(p)|.
\end{align*}
By the convolution formula, Plancherel's theorem, and the Cauchy–Schwarz inequality,
\begin{align} \label{three_particle_term_0}
    B_{t;u,v} \lesssim || \widehat{h}_{t}||_{L^{4}(\mathbb{R}^{2})}^{2} || \widehat{g_{t}}||_{L^{4}(\mathbb{R}^{2})}^{2}. 
\end{align}
In particular, applying  the two-dimensional Gagliardo–Nirenberg inequality \eqref{equation:Gagliardo_Nirenberg} yields
\begin{align} \label{three_particle_term_1}
    || \widehat{h}_{t}||_{L^{4}(\mathbb{R}^{2})}^{2}
    \lesssim ||h_{t}||_{L^{2}(\mathbb{R}^{2})} \||p| h_{t}\|_{L^{2}(\mathbb{R}^{2})} = \sqrt{A_{t;u} (-2A'_{t;u})}.
\end{align}
Similarly, we have
\begin{align} \label{three_particle_term_2}
    || \widehat{g_{t}}||_{L^{4}(\mathbb{R}^{2})}^{2} \lesssim \sqrt{D_{t;v} (-4D'_{t;v})}, \quad D_{t;v} := \int_{\mathbb{R}^{2}} dp \exp\bkt{-\frac{t}{4}|p|^{2}} |V(p)|^{2}.
\end{align}
As a consequence, combining \eqref{three_particle_term_0}-\eqref{three_particle_term_2} implies
\begin{align*}
    K_{N}^{(3)}[v,u] &\lesssim \int_{0}^{\infty} dt \exp(-t) A_{t;u}^{N-5/2}  (-2A'_{t;u})^{1/2} (D_{t;v} (-4D'_{t;v}))^{1/2} \\
    &\leq \bkt{\int_{0}^{\infty} dt A_{t;u}^{2N-5} (-2A'_{t;u}) }^{1/2} \bkt{\int_{0}^{\infty} dt D_{t;v} (-4D'_{t;v})  }^{1/2} = O(1/\sqrt{N}),
\end{align*}
where we have used Cauchy-Schwarz inequality  in the second inequality. The proof is complete. 

\end{proof}

\subsection{Asymptotic optimizer for the mean-field exponent}

We now prove the upper bound in Theorem \ref{proposition:mean_field_two_side}.
From the representation \eqref{equation:expression_mean_field_exponent}, choosing a specific $v \star u \in \mathscr{A}_{\textup{M-F}}(\mathbb{R}^{2N-2})$, where 
\begin{align} \label{choice_u_and_v}
    v := \frac{Q^{2}}{||Q||^{2}_{L^{4}(\mathbb{R}^{2})}}, \quad u := \frac{Q}{||Q||_{L^{2}(\mathbb{R}^{2})}},
\end{align}
gives the upper bound
\begin{align} \label{equation:upper_bound_mean_field_expoenent_0}
    E_{N,\theta}^{\textup{M-F}} \leq \frac{1}{4\pi} \langle v\star u,\mathcal{J}_{-1} v\star u \rangle_{L^{2}_{\textup{sym}'}(\mathbb{R}^{2N-2})} - \frac{\beta_{\star}(\theta)}{4\pi} -I^{(4)}_{N}[v,u].
\end{align}
Here we have used $K_{N}^{(3)}[v,u] \geq 0$. The function $Q$ is an optimizer  of the Gagliardo-Nirenberg inequality \eqref{equation:Gagliardo_Nirenberg}. We state its properties from \cite[Theorem B]{cmp/1103922134} below. 
\begin{lemma} \label{lemma:G_N_minimizer}
There exists a positive radial function $Q \in H^{1}(\mathbb{R}^{2}) \cap C^{\infty}(\mathbb{R}^{2})$ solving
\begin{align} \label{equation:stationary_NLS}
    \Delta Q - Q + Q^{3} = 0 \quad \text{on} \quad \mathbb{R}^{2}
\end{align}
such that
\begin{align} \label{equation:NLS_energy_functional}
    \frac{||Q||_{L^{2}(\mathbb{R}^{2})}^{2} ||\nabla Q||_{L^{2}(\mathbb{R}^{2})}^{2} }{||Q||_{L^{4}(\mathbb{R}^{2})}^{4}} = \frac{1}{C_{\textup{G-N}}}.
\end{align}
\end{lemma}

To see that $v \star u \in \mathscr{A}_{\textup{M-F}}(\mathbb{R}^{2N-2})$, we state the regularity of $Q$ below.
Applying  the Sobolev embedding \cite[Corollary 9.11]{brezis2011functional} together with $Q \in H^{1}(\mathbb{R}^{2})$ gives  
\begin{align} \label{equation:Q_L_p_0}
    Q  \in L^{p}(\mathbb{R}^{2}) \quad \forall p \in [2,\infty).
\end{align}
This together with the equation \eqref{equation:stationary_NLS} gives $\Delta Q = Q-Q^{3} \in L^{2}(\mathbb{R}^{2})$, which yields
\begin{align} \label{equation:Q_L_p}
    Q  \in H^{2}(\mathbb{R}^{2}).
\end{align}
Therefore, we conclude $(v,u) \in H^{2}(\mathbb{R}^{2}) \times H^{1}(\mathbb{R}^{2})$, which implies $v \star u \in \mathscr{A}_{\textup{M-F}}(\mathbb{R}^{2N-2})$.

The choice \eqref{choice_u_and_v} is suggested by the equality conditions
in the estimates used to prove Lemma \ref{lemma:sup_K_4}. As \eqref{large_t_contribution} below shows, the leading contribution to $K^{(4)}_{N}[v,u]$ comes from the small-$t$ region. Since $X_{0;u}=U*U$, the Cauchy-Schwarz inequality in \eqref{C_S_oper} suggests choosing $V$ proportional to $U*U$, while \eqref{G_N_oper} suggests choosing $u$ to be a normalized Gagliardo-Nirenberg optimizer. The functions in \eqref{choice_u_and_v}
satisfy both requirements.

We begin by estimating the logarithmic term in
\eqref{equation:upper_bound_mean_field_expoenent_0}.
\begin{lemma} \label{lemma:upper_bound_mean_field_expoenent_0_J}
For the choice in \eqref{choice_u_and_v} and all sufficiently large $N$, it holds that
\begin{align*}
    \frac{1}{4\pi} \langle v\star u,\mathcal{J}_{-1} v\star u \rangle_{L^{2}_{\textup{sym}'}(\mathbb{R}^{2N-2})} \leq \log N + O(1).
\end{align*}
\end{lemma}
\begin{proof}
Applying Jensen's inequality gives
\begin{align} \label{J_-1_upper_bound}
    \langle v\star u,\mathcal{J}_{-1} v\star u \rangle_{L^{2}_{\textup{sym}'}(\mathbb{R}^{2N-2})} \leq \log\bkt{ C_{v} +(N-2)C_{u} + 1  }.
\end{align}
Here we have set
\begin{align*}
    C_{v}:=\int_{\mathbb{R}^{2}} dp \frac{|p|^{2}}{2} |V(p)|^{2}, \quad C_{u}:=\int_{\mathbb{R}^{2}} dp |p|^{2} |U(p)|^{2}.
\end{align*}
By Lemma \ref{lemma:G_N_minimizer}, $Q \in H^{1}(\mathbb{R}^{2})$ is non-zero,  and hence,  Plancherel's theorem shows that $C_{u}$ is a finite positive constant. In addition, the Sobolev embedding \cite[Corollary 9.13]{brezis2011functional} yields $Q  \in H^{2}(\mathbb{R}^{2}) \subseteq L^{\infty}(\mathbb{R}^{2})$. This together with $Q \in H^{1}(\mathbb{R}^{2})$ then gives
\begin{align*}
    ||2Q\nabla Q||_{L^{2}(\mathbb{R}^{2})} \leq 2||Q||_{L^{\infty}(\mathbb{R}^{2})} ||\nabla Q||_{L^{2}(\mathbb{R}^{2})} < \infty.
\end{align*}
Together with Plancherel’s theorem, this shows that $C_{v}$ is a finite non-negative number. The proof is complete. 
\end{proof}

It remains to identify the leading-order contribution of the
four-particle term.
\begin{lemma}  \label{lemma:upper_bound_mean_field_expoenent_0_K}
For the choice in \eqref{choice_u_and_v} and all sufficiently large $N$, it holds that
\begin{align} \label{K_4_for_Q}
    K_{N}^{(4)}[v,u] = \frac{2C_{\textup{G-N}}}{N} + o(1/N).
\end{align}
\end{lemma}
\begin{proof}
We first introduce a time cutoff $C>0$, which will be chosen
sufficiently small below. By the Cauchy--Schwarz inequality, Plancherel's theorem, and \eqref{equation:Q_L_p_0}, we get
\begin{align} \label{limit_A_6}
    |\langle V,X_{t;u} \rangle_{L^{2}(\mathbb{R}^{2})} |^{2} \lesssim \|u\|_{L^{4}(\mathbb{R}^{2})}^{4} <\infty\quad \forall t\geq 0.
\end{align}
Since $t\mapsto A_{t;u}$ is non-increasing and $A_{C;u}<1$, it holds that
\begin{align} \label{large_t_contribution}
    &\limsup_{N\to\infty} N \int_{C}^{\infty} dt \exp(-t) A_{t;u}^{N-4} \abs{ \langle V,X_{t;u} \rangle_{L^{2}(\mathbb{R}^{2})} }^{2} \lesssim  \limsup_{N\to\infty} N A_{C;u}^{N-4} \|u\|_{L^{4}(\mathbb{R}^{2})}^{4} = 0.
\end{align}
Thus, the contribution from $t\geq C$ is $o(1/N)$.

Rescaling time in the remaining integral gives
\begin{align} \label{limit_A}
     N\frac{1}{(2\pi)^{2}} \int_{0}^{C} dt \exp(-t) A_{t;u}^{N-4} \abs{ \langle V,X_{t;u} \rangle_{L^{2}(\mathbb{R}^{2})} }^{2} =  \frac{1}{(2\pi)^{2}} \int_{0}^{CN} dt \exp\bkt{\frac{-t}{N}} A_{t/N;u}^{N-4} \abs{ \langle V,X_{t/N;u} \rangle_{L^{2}(\mathbb{R}^{2})} }^{2}.
\end{align}
It suffices to compute  the right-hand side of \eqref{limit_A}. To do this, we now estimate $A_{t/N;u}$. 
For $N\geq 5$, write
\begin{align} \label{limit_A_1}
    A_{t/N;u} = 1-\frac{1}{N-4} (N-4)\int_{\mathbb{R}^{2}} dp \bkt{1-\exp\bkt{-\frac{t}{N2}|p|^{2}}} |U(p)|^{2} =: 1-\frac{x_{N}(t)}{N-4}.
\end{align}
By using  $1-\exp(-x) \leq x$ for all $x\geq 0$,  the constant $C$ can be chosen sufficiently small such that
\begin{align} \label{limit_A_2}
    x_{N}(t) \leq \frac{N-4}{N} \frac{t}{2} \int_{\mathbb{R}^{2}} dp |p|^{2} |U(p)|^{2}
    \leq \frac{C N}{2} \int_{\mathbb{R}^{2}} dp |p|^{2} |U(p)|^{2} < \frac{N}{4} \quad \forall 0<t<CN.
\end{align}
One easily sees that
\begin{align} \label{limit_A_3}
    \bkt{1-\frac{x}{N}}^{N} \leq \exp(-x) \quad \forall 0\leq x<N/2, \quad N\geq 2. 
\end{align}
Thus, combining \eqref{limit_A_1}-\eqref{limit_A_3} gives
\begin{align} \label{limit_A_4}
    A_{t/N;u}^{N-4} = \bkt{1-\frac{x_{N}(t)}{N-4}}^{N-4} \leq \exp(-x_{N}(t)) \quad \forall 0<t<CN, \quad N\geq 9.
\end{align}
In particular, since $1-\exp(-x) \geq x-x^{2}/2$ for all $x\geq 0$, choosing $C$ to be sufficiently small yields that for all sufficiently large $N$ and each $0<t<CN$, 
\begin{align} \label{limit_A_5}
    x_{N}(t) \geq t\frac{N-4}{2N} \int_{\mathbb{R}^{2}} dp |p|^{2} |U(p)|^{2} - t^{2} \frac{N-4}{4N^{2}} \int_{\mathbb{R}^{2}} dp |p|^{4} |U(p)|^{2} \geq \frac{t}{8}\int_{\mathbb{R}^{2}} dp |p|^{2} |U(p)|^{2}.
\end{align}
Here, due to Lemma \ref{lemma:G_N_minimizer} and \eqref{equation:Q_L_p}, both integrals in \eqref{limit_A_5} are finite and positive. Therefore, combining \eqref{limit_A_4} and \eqref{limit_A_5} shows that for all sufficiently large $N$ and each $0<t<CN$, 
\begin{align} \label{limit_A_7}
    A_{t/N;u}^{N-4} \leq \exp\bkt{-\frac{t}{8} \int_{\mathbb{R}^{2}} dp |p|^{2} |U(p)|^{2}}. 
\end{align}

By using \eqref{limit_A_6} and \eqref{limit_A_7}, applying the dominated convergence theorem yields
\begin{align} \label{limit_A_8}
    &\lim_{N\to\infty} \text{(R.H.S.) of \eqref{limit_A}} = \frac{1}{(2\pi)^{2}} \int_{0}^{\infty} dt \lim_{N\to\infty}  A_{t/N;u}^{N-4} \abs{ \langle V,X_{t/N;u} \rangle_{L^{2}(\mathbb{R}^{2})} }^{2}.
\end{align}
We now identify the pointwise limit. For each fixed $t\geq 0$, the dominated convergence theorem gives
\begin{align*}
    \lim_{N\to\infty} x_{N}(t) = \frac{t}{2} \int_{\mathbb{R}^{2}} dp |p|^{2} |U(p)|^{2},
\end{align*}
which implies
\begin{align} \label{limit_A_9}
    \lim_{N\to\infty} A_{t/N;u}^{N-4}
    = \exp\bkt{-\frac{t}{2} ||\nabla u||_{L^{2}(\mathbb{R}^{2})}^{2} }.
\end{align}
By using Plancherel’s theorem, one easily shows that
\begin{align} \label{limit_A_10}
    \lim_{N\to\infty}\frac{1}{(2\pi)^{2}} | \langle V,X_{t/N;u} \rangle_{L^{2}(\mathbb{R}^{2})} |^{2}
     = \frac{||Q||_{L^{4}(\mathbb{R}^{2})}^{4}}{||Q||_{L^{2}(\mathbb{R}^{2})}^{4}}.
\end{align}
As a consequence, combining \eqref{limit_A_8}-\eqref{limit_A_10} and \eqref{equation:NLS_energy_functional} gives
\begin{align*}
    \lim_{N\to\infty} \text{(R.H.S.) of \eqref{limit_A}}
    = \frac{2||Q||_{L^{4}(\mathbb{R}^{2})}^{4}}{||\nabla u||_{L^{2}(\mathbb{R}^{2})}^{2}||Q||_{L^{2}(\mathbb{R}^{2})}^{4}}  = 2C_{\textup{G-N}}.
\end{align*}
The proof is complete. 
\end{proof}

The upper bound in Theorem \ref{proposition:mean_field_two_side} now follows by combining Lemma \ref{lemma:upper_bound_mean_field_expoenent_0_J} and Lemma \ref{lemma:upper_bound_mean_field_expoenent_0_K}.

\appendix

\section{A beta-function estimate for small parameters} \label{section:beta_property_1}

In this section, we prove Lemma \ref{lemma:beta_property_1}. We begin by proving a general upper bound for $B(p,k)$. Fix $0<p<1$. By using the log-convexity of  $\Gamma$, we know that
\begin{align*}
    \Gamma(k) k= \Gamma(k+1) \leq \Gamma(k+p)^{p} \Gamma(k+p+1)^{1-p} = (k+p)^{1-p} \Gamma(k+p),
\end{align*}
which implies that
\begin{align*}
    B(p,k) \leq \Gamma(p) \frac{(k+p)^{1-p}}{k} = \frac{\Gamma(p)}{k^{p}} \bkt{1+\frac{p}{k}}^{1-p}.
\end{align*}
Thus, for $k\geq p$, 
\begin{align*}
    B(p,k) \leq \frac{\Gamma(p)}{k^{p}} 2^{1-p} \leq 2 \frac{\Gamma(p)}{k^{p}},
\end{align*}
and for $p\geq k$, 
\begin{align*}
    B(p,k) \leq \frac{\Gamma(p)}{k^{p}} \bkt{2\frac{p}{k}}^{1-p} \leq \frac{2}{k} \frac{\Gamma(p+1)}{p^{p}}.
\end{align*}
Note that both  $\Gamma(p+1)$ and $p^{-p}$ are bounded on $0<p<1$. The proof is complete.

\section*{Acknowledgement}

We are grateful to Prof. Li-Cheng Tsai for introducing us to the resolvent
formula \eqref{equation:resolvent_representation}, discussing its behavior
as $N\to\infty$ (see the discussion below
Theorem \ref{proposition:mean_field_two_side}), and explaining
Rajeev's conjecture to us. We appreciate Prof. Kyeongsik Nam's explanations of the
key ideas in \cite{ganguly2025sharp}, which were helpful in the early stages of developing our proofs.
We also thank Prof. Nikos Zygouras
and Colin Piernot for helpful discussions. The first author acknowledges support from EPSRC grant EP/V026100/1,
the Swiss National Science Foundation project MINT (10000849),
and NCCR SwissMAP.

\section*{Declaration of AI use}

We have used AI for language editing, literature searches, and a final check of selected
calculations and arguments. All ideas and proofs in this work were developed by the authors.

\printbibliography

\end{document}